\documentclass[11pt,leqno]{article}

\usepackage[english]{babel}

\usepackage{amsmath}

\usepackage{amssymb}

\usepackage{xypic}

\usepackage[all]{xy}

\usepackage{latexsym}

\usepackage{amsfonts}

\usepackage{epsf}

\usepackage[dvips]{graphicx}

\usepackage[latin1]{inputenc}

\usepackage[T1]{fontenc}

\numberwithin{equation}{section}

\newcommand{\A}{\mathcal{A}}

\newcommand{\B}{\mathcal{B}}

\newcommand{\C}{\mathcal{C}}

\newcommand{\D}{\mathcal{D}}

\newcommand{\K}{\mathcal{K}}

\newcommand{\M}{\mathcal{M}}

\newcommand{\T}{\mathcal{T}}

\newcommand{\W}{\mathcal{W}}

\renewcommand{\mod}{\mathrm{Mod}}

\newcommand{\ind}{\mathrm{Ind}}

\newcommand{\cov}{\mathrm{Cov}}

\newcommand{\coh}{\mathrm{Coh}}

\newcommand{\op}{\mathrm{Op}}

\newcommand{\rc}{\mathbb{R}\textrm{-}\mathrm{c}}

\newcommand{\cc}{\mathbb{C}\textrm{-}\mathrm{c}}

\newcommand{\CC}{\mathbb{C}}

\newcommand{\R}{\mathbb{R}}

\newcommand{\Z}{\mathbb{Z}}

\newcommand{\N}{\mathbb{N}}

\newcommand{\OO}{\mathcal{O}}

\newcommand{\NN}{\mathcal{N}}

\newcommand{\RR}{\mathcal{R}}

\newcommand{\E}{\mathcal{E}}

\newcommand{\F}{\mathcal{F}}

\newcommand{\G}{\mathcal{G}}

\newcommand{\I}{\mathrm{I}}

\newcommand{\II}{\mathcal{I}}

\newcommand{\ltens}{\overset{L}{\otimes}}

\newcommand{\iso}{\stackrel{\sim}{\to}}

\newcommand{\dbtxr}{\mathcal{D}\mathit{b}^t_{X_\mathbb{R}}}

\newcommand{\dbt}{\mathcal{D}\mathit{b}^t}

\newcommand{\db}{\mathcal{D}\mathit{b}}

\newcommand{\ot}{\mathcal{O}^t}

\newcommand{\ol}{\mathcal{O}^\lambda}

\newcommand{\dd}[3]{\mathcal{D}_{{#1} \stackrel{#2}{\to} {#3}}}

\newcommand{\ddual}[3]{\mathcal{D}_{{#1} \stackrel{#2}{\gets} {#3}}}

\newcommand{\ddxy}{\mathcal{D}_{\xmenoy}}

\newcommand{\ddyx}{\mathcal{D}_{\ymenox}}

\newcommand{\xmenoy}{{X}\rightarrow{Y}}

\newcommand{\ymenox}{{Y}\leftarrow{X}}

\newcommand{\ddmn}{\mathcal{D}_{M \to N}}

\newcommand{\ddnm}{\mathcal{D}_{N \gets M}}

\newcommand{\tho}{\mathit{T}\mathcal{H}\mathit{om}}

\newcommand{\muh}{\mu\mathit{hom}}

\newcommand{\mucirc}{\underset{\mu}{\circ}}

\newcommand{\OW}{\OO^\mathrm{w}}
\newcommand{\OWX}{\OO^\mathrm{w}_X}
\newcommand{\OWY}{\OO^\mathrm{w}_Y}
\newcommand{\CWXR}{\C^{{\infty ,\mathrm{w}}}_{X_\R}}
\newcommand{\CWYR}{\C^{{\infty ,\mathrm{w}}}_{Y_\R}}
\newcommand{\CWM}{\C^{{\infty ,\mathrm{w}}}_M}
\newcommand{\CWN}{\C^{{\infty ,\mathrm{w}}}_N}
\newcommand{\CW}{\C^{{\infty ,\mathrm{w}}}}
\newcommand{\wtens}{\overset{\mathrm{w}}{\otimes}}

\newcommand{\vchar}{\mathrm{Char}}

\newcommand{\sol}{\mathcal{S}ol}

\newcommand{\ri}{\mathit{R}\mathcal{I}\mathit{hom}}

\newcommand{\rh}{\mathit{R}\mathcal{H}\mathit{om}}

\newcommand{\ho}{\mathcal{H}\mathit{om}}

\newcommand{\ih}{\mathcal{I}\mathit{hom}}

\newcommand{\Ho}{\mathrm{Hom}}

\newcommand{\Rh}{\mathrm{RHom}}

\newcommand{\id}{\mathrm{id}}

\renewcommand{\dim}{\textbf{Proof.}}

\newcommand{\qed}{\nopagebreak \phantom{} \hfill $\Box$ \\}

\newcommand{\formX}{\omega^{\otimes-1}_{\Delta|X \times X}}

\newcommand{\pT}{\overset{\hspace{0.1cm}}{\dot{T}}}

\newcommand{\pE}{\overset{\hspace{0.1cm}}{\dot{E}}}

\newcommand{\pU}{\overset{\hspace{0.1cm}}{\dot{U}}}

\newcommand{\pV}{\overset{\hspace{0.1cm}}{\dot{V}}}

\newcommand{\supp}{\mathrm{supp}}

\newcommand{\pittau}{\stackrel{\mathbf{.}}{\tau}}

\newcommand{\ppi}{\dot{\pi}}

\renewcommand{\H}{ \mathcal{H}}

\newcommand{\RP}{\mathbb{R}^{{\scriptscriptstyle{+}}}}

\newcommand{\imin}[1]{#1^{-1}}

\newcommand{\lind}[1]{\underset{#1}{\underrightarrow{\lim}}}

\newcommand{\Lind}{\underrightarrow{\lim}}  

\newcommand{\indl}[1]{\underset{#1}{``\underrightarrow{\mathrm{lim}}\mbox{''}}}

\newcommand{\Lpro}{\underleftarrow{\lim}}

\newcommand{\lpro}[1]{\underset{#1}{\underleftarrow{\lim}}}

\newcommand{\exs}[3]{0 \to {#1} \to {#2} \to {#3} \to 0}

\newcommand{\lexs}[3]{0 \to {#1} \to {#2} \to {#3}}

\newcommand{\dt}[3]{{#1} \to {#2} \to {#3} \stackrel{+}{\to}}

\newtheorem{teo}{Theorem}[subsection]

\newtheorem{df}[teo]{Definition}

\newtheorem{cor}[teo]{Corollary}

\newtheorem{oss}[teo]{Remark}

\newtheorem{prop}[teo]{Proposition}

\newtheorem{lem}[teo]{Lemma}

\newtheorem{es}[teo]{Example}

\newtheorem{nt}[teo]{Notations}

\author{Luca Prelli}
\title{\bf{SPECIALIZATION AND MICROLOCALIZATION OF SUBANALYTIC SHEAVES}}
\date{}

\usepackage{latexsym}
\begin{document}

\maketitle

\begin{abstract}
In this paper we define the specialization and microlocalization functors for
subanalytic sheaves. Then we specialize and microlocalize the sheaves of
tempered and Whitney holomorphic functions generalizing some classical constructions.
\end{abstract}

\tableofcontents

\addcontentsline{toc}{section}{\textbf{Introduction}}

\section*{Introduction}

After the fundamental works of Sato on hyperfunctions and
microfunctions and the developement of algebraic analysis, the
methods of cohomological theory of sheaves became very useful for
studying systems of PDE on real or complex analytic manifolds.
Unfortunately sheaf theory is not well suited to study objects
which are not defined by local properties. Since the study of
solutions of a system of PDE in these spaces is very important
(Riemann-Hilbert correspondence, Laplace transform, etc.), many
ways have been explored by the specialists to overcome this
problem. First Kashiwara in \cite{Ka84} defined the functor of
tempered cohomology to solve the Riemann-Hilbert problem for
holonomic systems. Its dual, the functor of formal cohomology was
introduced by Kashiwara and Schapira in \cite{KS96}. Their
microlocal analogous were defined by Andronikof in \cite{An94} and
Colin in \cite{Co98}. These results present a problem: they are
construction ``ad hoc'', and they are not contained in a unifying
theory based on the six Grothendieck operations. In \cite{KS01}
Kashiwara and Schapira introduced the notion of ind-sheaves, and
defined the six Grothendieck operations in this framework. Then
they defined the subanalytic site (a site whose open sets are
subanalytic and the coverings are locally finite) and they proved
the equivalence between subanalytic sheaves and
ind-$\R$-constructible sheaves. They obtained the formalism of the
six Grothendieck operations by including subanalytic sheaves into
the category of ind-sheaves. In \cite{Pr1} a direct,
self-contained and elementary construction of the six Grothendieck
operations for subanalytic sheaves is introduced without using the
more sophisticated and much more difficult theory of ind-sheaves.
In the beginning of this paper we extend some classical
constructions for sheaves, as the Fourier-Sato transform and the
functors of specialization and microlocalization. We introduce
first the category of conic subanalytic sheaves on an analytic
manifold. Then we extend the Fourier-Sato transform to the
category of conic subanalytic sheaves on a vector bundle. At this
point we can start studying subanalytic sheaves from a microlocal
point of view by introducing the functors of specialization and
microlocalization along a submanifold of a real analytic manifold.
Roughly speaking, starting from a sheaf on a real analytic
manifold, we construct conic sheaves on the tangent and the
cotangent bundle respectively. We also show that the functor of
microlocalization is related with the
functor of ind-microlocalization defined in \cite{KSIW}.
Then, applying specialization
(resp. microlocalization) to the subanalytic sheaves of tempered
and Whitney holomorphic functions, we generalize tempered and
formal specialization (resp. microlocalization). In
this way we get a unifying description of Andronikov's and
Colin's ``ad hoc'' constructions. This is also a ``simplification'',
since the definitions of specialization and microlocalization for
subanalytic sheaves are more intuitive for people which
are familiar with the classical sheaf theory of \cite{KS90},
although working with the subanalytic site is delicate and requires
results which are related to the geometry of subanalytic subsets.\\

In more details the contents of this work are as follows.\\

In {\bf Section 1} we recall the results on subanalytic sheaves of \cite{KS01} and \cite{Pr1}.\\

In {\bf Section 2} we construct the category of conic sheaves on a
subanalytic site endowed with an action of $\RP$.\\

In {\bf Section 3} we consider a vector bundle $E$ over a real
analytic manifold and its dual $E^*$ endowed with the natural
action of $\RP$. We define the Fourier-Sato transform
which gives an equivalence between conic subanalytic sheaves on $E$ and conic subanalytic sheaves on $E^*$.\\

Then we define the functor $\nu^{sa}_M$ of specialization along a sub\-manifold $M$ of a real analytic
manifold $X$ ({\bf Section 4}) and its Fourier-Sato transform, the functor
$\mu^{sa}_M$ of microlocalization ({\bf Section 5}). We introduce the functor
$\muh^{sa}$ for subanalytic sheaves and we give an estimate of its support using the notion of microsupport of \cite{KS03}. Then we study its relation with the functor of ind-microlocalization of \cite{KSIW}.\\

We apply these results in {\bf Section 6}. We study the connection
between specialization and microlocalization for subanalytic sheaves and the
classical ones. Specialization of subanalytic sheaves generalize tempered and formal specialization of \cite{An94} and \cite{Co01}, in particular when we specialize Whitney holomorphic functions we obtain the sheaves of functions asymptotically developable of \cite{Ma91} and \cite{Si90}. Moreover, thanks to the functor of microlocalization, we are able to generalize
tempered and formal microlocalization introduced by Andronikof in
\cite{An94} and Colin in \cite{Co98} respectively.\\

{\bf Section 7} is dedicated to the study of the microlocalization of tempered and Whitney holomorphic functions. We prove that the microlocalization of $\ot$ and $\OW$ have (in cohomology) a natural structure of $\E$-module and that locally they are invariant under contact transformations.\\

In {\bf Section 8} we study the Cauchy-Kowaleskaya-Kashiwara theorem with growth conditions. The proof is not a rephrasing of Kashiwara's
proof, which consisted in a reduction to the case of one operator: we first
prove a propagation result (already known for $\ot_X$)
and we use this propagation result to reduce the statements to the commutation
of duality with non-characteristic inverse image. \\

We end this work with a short {\bf Appendix} in which we recall the definitions and we collect
some properties of subanalytic subsets and ind-sheaves, and then we study the inverse
image of the subanalytic sheaves of tempered and Whitney holomorphic functions.\\

\noindent \textbf{Acknowledgments.} We would like to thank Prof.
Pierre Schapira who encouraged us to develop specialization and
microlocalizaiton of subanalytic sheaves, and Prof. A. D'Agnolo
for his many useful remarks.

\section{Review on sheaves on subanalytic sites}\label{review}

In the following $X$ will be a real analytic manifold and $k$ a
field. Reference are made to \cite{KS} for a complete exposition
on sheaves on Grothendieck topologies, to \cite{KS01} and
\cite{Pr1} for an introduction to sheaves on subanalytic sites. We refer to
\cite{BM88} and \cite{Lo93} for the theory of subanalytic sets.

\subsection{Sheaves on subanalytic sites}

Let us recall some results of \cite{KS01} and \cite{Pr1}.

Denote by $\op(X_{sa})$ the category of open subanalytic subsets of
$X$. One endows $\op(X_{sa})$ with the following topology: $S
\subset \op(X_{sa})$ is a covering of $U \in \op(X_{sa})$ if for
any compact $K$ of $X$ there exists a finite subset $S_0\subset S$
such that $K \cap \bigcup_{V \in S_0}V=K \cap U$. We will call
$X_{sa}$ the subanalytic site, and for $U \in \op(X_{sa})$ we
denote by $U_{X_{sa}}$ the category $\op(X_{sa}) \cap U$ with the
topology
induced by $X_{sa}$. We denote by $\op^c(X_{sa})$ the subcategory of
$\op(X_{sa})$ consisting of relatively compact open subanalytic subsets.\\

Let $\mod(k_{X_{sa}})$ denote the category of sheaves on $X_{sa}$.
Then $\mod(k_{X_{sa}})$ is a Grothendieck category, i.e. it admits
a generator and small inductive limits, and small filtrant
inductive limits are exact. In particular as a Grothendieck
category, $\mod(k_{X_{sa}})$ has enough injective objects.\\

Let $\mod_{\rc}(k_X)$  be the abelian category of
$\R$-constructible sheaves on $X$, and consider its subcategory
$\mod^{c}_{\rc}(k_X)$ consisting of sheaves whose support is
compact.\\

We denote by $\rho: X \to X_{sa}$ the natural morphism of sites.
We have functors


\begin{equation*}
\xymatrix{\mod(k_X)   \ar@ <5pt> [rr]^{\rho_*} \ar@
<-5pt>[rr]_{\rho_!} &&
  \mod(k_{X_{sa}}) \ar@ <0pt> [ll]|{\imin \rho}. }
\end{equation*}

The functors $\imin \rho$ and $\rho_*$ are the functors of inverse
image and direct image respectively. The sheaf $\rho_!F$ is the
sheaf associated to the presheaf $\op(X_{sa}) \ni U \mapsto
F(\overline{U})$. In particular, for $U \in \op(X)$ one has
$\rho_!k_U \simeq \lind {V \subset \subset U} \rho_*k_V$, where $V
\in \op(X_{sa})$. Let us summarize the properties of these
functors:

\begin{itemize}

\item the functor $\rho_*$ is fully faithful and left exact, the restriction
of $\rho_*$ to $\mod_{\rc}(k_X)$ is exact,
\item the functor $\imin \rho$ is exact,
\item the functor $\rho_!$ is fully faithful and exact,
\item $(\imin \rho,\rho_*)$ and $(\rho_!,\imin \rho)$ are pairs of
adjoint functors.
\end{itemize}

\begin{nt} Since the functor $\rho_*$ is fully faithful and exact on
$\mod_{\rc}(k_X)$, we can identify $\mod_{\rc}(k_X)$ with its
image in $\mod(k_{X_{sa}})$. When there is no risk of confusion we
will write $F$ instead of $\rho_*F$, for $F \in \mod_{\rc}(k_X)$.
\end{nt}


Let $F \in \mod(k_{X_{sa}})$. There exists a filtrant inductive system $\{F_i\}_{i\in I}$ in $\mod_{\rc}^c(k_X)$ such that $F \simeq \lind i \rho_* F_i$.\\

Let $X,Y$ be two real analytic manifolds, and let $f:X \to Y$ be a
real analytic map. We have a commutative diagram
\begin{equation}
\xymatrix{X \ar[d]^\rho \ar[r]^f & Y \ar[d]^\rho \\
X_{sa} \ar[r]^f & Y_{sa}}
\end{equation}

We get external operations $\imin f$ and $f_*$, which are always
defined for sheaves on Grothendieck topologies. For subanalytic
sheaves we can also define the functor of proper direct image
\begin{eqnarray*}
f_{!!}:\mod(k_{X_{sa}}) & \to & \mod(k_{Y_{sa}}) \\
F & \mapsto & \lind U f_*F_U \simeq \lind K f_*\Gamma_KF
\end{eqnarray*}
where $U$ ranges trough the family of relatively compact open
subanalytic subsets of $X$ and $K$ ranges trough the family of
subanalytic compact subsets of $X$. The notation $f_{!!}$ follows
from the fact that $f_{!!} \circ \rho_* \not\simeq \rho_* \circ
f_!$ in general. If $f$ is proper on $\supp(F)$ then $f_*F \simeq
f_{!!}F$, in this case $f_{!!}$ commutes with $\rho_*$. While
functors $\imin f$ and $\otimes$ are exact, the functors $\ho$,
$f_*$ and $f_{!!}$ are left exact and admits right derived
functors.

To derive these functors we use the category of quasi-injective
objects. An object $F \in \mod(k_{X_{sa}})$ is
quasi-injective if for $U,V \in \op^c(X_{sa})$ with $V \subset U$ the restriction morphism $\Gamma(U;F) \to \Gamma(V;F)$ is surjective or, equivalently, if the functor $\Ho_{k_{X_{sa}}}(\cdot,F)$ is
exact on $\mod^c_{\rc}(k_X)$. Quasi injective objects are
injective with respect to the functors $f_*,f_{!!}$ and, if $G \in
\mod_{\rc}(k_X)$, with respect to the functors
$\Ho_{k_{X_{sa}}}(G,\cdot),\ho(G,\cdot)$.

The functor $Rf_{!!}$ admits a right adjoint, denoted by $f^!$,
and we get the usual isomorphisms between Grothendieck operations
(projection formula, base change formula, K\"unneth formula, etc.)
in the framework of subanalytic sites.

Let $Z$ be a subanalytic locally closed subset of $X$. As in
classical sheaf theory we define
\begin{eqnarray*}
\Gamma_Z: \mod(k_{X_{sa}}) & \to & \mod(k_{X_{sa}})\\
F & \mapsto & \ho(\rho_*k_Z,F) \\
(\cdot)_Z: \mod(k_{X_{sa}}) & \to & \mod(k_{X_{sa}}) \\
F & \mapsto & F \otimes \rho_*k_Z.
\end{eqnarray*}

Finally we recall the properties of the six Grothendieck
operations and their relations with the functors $\imin \rho$,
$R\rho_*$ and $\rho_!$. We refer to \cite{Pr1} for a detailed
exposition.
\begin{itemize}
\item the functor $R^k\ho(\rho_*F,\cdot)$ commutes with filtrant $\Lind$ if $F \in \mod_{\rc}(k_X)$,
\item the functors $R^kf_{!!}$ and $H^kf^!$ commute with filtrant $\Lind$,
\item the functor $\imin \rho$ commutes with $\otimes$, $\imin f$ and
$Rf_{!!}$,
\item the functor $R\rho_*$ commutes with $\rh$, $Rf_*$ and
$f^!$,
\item the functor $\rho_!$ commutes with $\otimes$ and $\imin f$,
\item the restriction of $\otimes$ and $\imin f$ to the category
of $\R$-constructible sheaves commute with $\rho_*$,
\item if $f$ is a topological submersion (i.e. it is locally
isomorphic to a projection $Y \times \R^n \to Y$), then $f^!\simeq
\imin f \otimes f^!k_Y$ commutes with $\imin \rho$ and $Rf_{!!}$
commutes with $\rho_!$.
\end{itemize}
Moreover the functors $Rf_*$, $Rf_{!!}$ and $\rh(F,\cdot)$ with $F
\in \mod_{\rc}(k_X)$ have finite cohomological dimension.

\subsection{Modules over a $k_{X_{sa}}$-algebra}

A sheaf of $k_{X_{sa}}$-algebras (or a $k_{X_{sa}}$-algebra, for
short) is an object $\RR \in \mod(k_{X_{sa}})$ such that
$\Gamma(U;\RR)$ is a $k$-algebra for each $U \in \op(X_{sa})$ and the restriction maps are algebra morphisms. A sheaf of (left) $\RR$-modules is a sheaf $F$ such
that $\Gamma(U;F)$ has a structure of (left)
$\Gamma(U;\RR)$-module for each $U \in \op(X_{sa})$.
Let us denote by $\mod(\RR)$ the
category of sheaves of (left) $\RR$-modules. The category
$\mod(\RR)$ is a Grothendieck category and  the forgetful functor $for:\mod(\RR) \to
\mod(k_{X_{sa}})$ is exact.\\

The functors
\begin{eqnarray*}
\ho_\RR & : &  \mod(\RR)^{op} \times \mod(\RR) \to \mod(k_{X_{sa}}),\\
\otimes_\RR & : &  \mod(\RR^{op}) \times \mod(\RR) \to
\mod(k_{X_{sa}})
\end{eqnarray*}
are well defined. Remark that in the case of $\RR$-modules the functor
$\otimes_\RR$ is only right exact and commutes with $\Lind$.

Let $X,Y$ be two real analytic manifolds, and let $f:X \to Y$ be a
morphism of real analytic manifolds. Let $\RR$ be a
$k_{Y_{sa}}$-algebra. The functors $\imin f$, $Rf_*$ and $Rf_{!!}$
induce functors
\begin{eqnarray*}
\imin f & : & \mod(\RR) \to \mod(\imin f
\RR ),\\
f_* & : &  \mod(\imin f
\RR ) \to \mod(\RR),\\
f_{!!} & : &  \mod(\imin f \RR ) \to \mod(\RR).
\end{eqnarray*}

Now we consider the derived category of sheaves of $\RR$-modules.

\begin{df} An object $F \in \mod(\RR)$ is flat if the functor $\mod(\RR^{op}) \ni G
\to G \otimes_\RR F$ is exact.
\end{df}

Thanks to flat objects we can find a left
derived functor $\otimes^L_\RR$ of the tensor product
$\otimes_\RR$.

\begin{df} An object $F \in \mod(\RR)$ is quasi-injective if its image via the forgetful functor
is quasi-injective in $\mod(k_{X_{sa}})$.
\end{df}

Let $X,Y$ be two real analytic manifolds, and let $f:X \to Y$ be a
real analytic map. Let $\RR$ be a $k_{Y_{sa}}$-algebra. One can
prove that quasi-injective objects are injective with respect to
the functors $f_*$ and $f_{!!}$. The functors $Rf_*$ and $Rf_{!!}$
are well defined and projection formula, base change formula remain
valid for $\RR$-modules. Moreover we have

\begin{teo} The functor $Rf_{!!}:D^+(\imin f \RR) \to D^+(\RR)$
admits a right adjoint. We denote by $f^!:D^+(\RR) \to D^+(\imin f
\RR)$ the adjoint functor.
\end{teo}

\section{Conic sheaves on subanalytic sites}

In this section we study the category of conic sheaves on a
subanalytic site. Reference are made to \cite{KS90} for the
classical theory of conic sheaves.

\subsection{Conic sheaves on topological spaces}

Let $k$ be a field. Let $X$ be a real analytic manifold endowed
with a subanalytic action $\mu$ of $\RP$. In other words we have a
subanalytic map
$$\mu: X \times \RP \to X,$$
which satisfies, for each $t_1,t_2 \in \RP$:
$$
  \begin{cases}
    \mu(x,t_1t_2)=\mu(\mu(x,t_1),t_2), \\
    \mu(x,1)=x.
  \end{cases}
$$
Note that $\mu$ is open. Indeed let $U \in \op(X)$ and $W
\in \op(\RP)$. Then $\mu(U,W)=\bigcup_{t \in
W}\mu(U,t)$, and $\mu(\cdot,t):X \to X$ is a homeomorphism
(with inverse $\mu(\cdot,\imin t)$). We have a diagram
$$\xymatrix{X \ar[r]^{\hspace{-5mm}j} & X \times \RP \ar@ <2pt>
[r]^{\hspace{0.3cm}\mu} \ar@ <-2pt> [r]_{\hspace{0.3cm}p}& X,}$$
where $j(x)=(x,1)$ and $p$ denotes the projection. We have $\mu
\circ j=p \circ j=\id$.

\begin{df} (i) Let $S$ be a subset of $X$. We set $\RP S=\mu(S,\RP).$
If $U \in \op(X)$, then $\RP U \in \op(X)$ since $\mu$ is open.

(ii) Let $S$ be a subset of $X$. We say that $S$ is conic if
$S=\RP S$. In other words, $S$ is invariant by the action of
$\mu$.

(iii) An  orbit of $\mu$ is the set $\RP x$ with $x \in X$.
\end{df}

\begin{df} We say that a subset $S$ of $X$ is $\RP$-connected if
$S \cap \RP x$ is connected for each $x \in S$.
\end{df}

\begin{df} A sheaf $F \in \mod(k_X)$ is conic if $\imin \mu F \simeq \imin p F$.
\begin{itemize}
\item[(i)]We denote by $\mod_{\RP}(k_X)$ the subcategory of
$\mod(k_X)$ consisting of conic sheaves.
\item[(ii)] We denote by $D^b_{\RP}(k_{X_{sa}})$,
the subcategory of $D^b(k_X)$ consisting of objects $F$
such that $H^j(F)$ belongs to $\mod_{\RP}(k_X)$ for all $j
\in \Z$.
\end{itemize}
\end{df}

\begin{prop} Let $U \in \op(X)$ be $\RP$-connected and let $F \in D^b_{\RP}(k_X)$. Then
$$
\mathrm{R}\Gamma(\RP U;F) \iso \mathrm{R}\Gamma(U;F).
$$
\end{prop}

Let us assume the following hypothesis
\begin{equation}\label{hypclass}
\begin{cases}
 \begin{array}{cc}\text{(i) every point $x \in X$ has a fundamental neighborhood}\\
 \text{\ \ system consisting of $\RP$-connected open subsets;}\\
 \end{array}\\
 \text{(ii) for any $x \in X$\ \ the set $\RP x$ is contractible.
}
  \end{cases}
\end{equation}
In this situation (see \cite{Be74}) either $\RP x \simeq \R$ or $\RP x = x$.\\

Denote by $X_{\RP}$ the topological space $X$ endowed with the
conic topology,
i.e. $U \in \op(X_{\RP})$ if it is open for the topology of $X$ and invariant by the action of $\RP$.\\

Let us consider the natural map $\eta: X \to X_{\RP}$. The restriction of $\eta_*$ induces an exact functor denoted by
$\widetilde{\eta}_*$ and we obtain a diagram
\begin{equation}\label{rho}
\xymatrix{\mod_{\RP}(k_X)  \ar[d] \ar@ <2pt>
[rr]^{\widetilde{\eta}_*} &&
  \mod(k_{X_{\RP}}) \ar@ <2pt> [ll]^{\imin \eta} \\
\mod(k_X) \ar[urr]_{\eta_*} && }
\end{equation}

Let $F \in D^b_{\RP}(k_X)$. Let $\varphi$ be the natural map from
$\mathrm{R}\Gamma(\RP U;F)$ to $\mathrm{R}\Gamma(U;\imin \eta F)$ defined by
\begin{equation}\label{varphiclass}
\begin{array}{ccc}
\mathrm{R}\Gamma(\RP U;F) & \to & \mathrm{R}\Gamma(\RP U; R\eta_* \imin \eta F) \\
  & \simeq & \mathrm{R}\Gamma(\RP U; \imin \eta F) \\
 & \to & \mathrm{R}\Gamma(U; \imin \eta F).
\end{array}
\end{equation}

\begin{prop}\label{*-1U} Let $F$ be a sheaf over $X_{\RP}$.
Let $U$ be an open set of $X$ and assume that $U$ is
$\RP$-connected. Then the morphism $\varphi$ defined by
\eqref{varphiclass} is an isomorphism.
\end{prop}

\begin{teo}\label{ex*class} The functors $R\eta_*$ and $\imin \eta$ in \eqref{rho} induce equivalences of derived ca\-te\-go\-ries
$$\xymatrix{D^b_{\RP}(k_X) \ar@ <2pt>
[r]^{\hspace{0cm}R\eta_*} & D^b(k_{X_{\RP}}) \ar@ <2pt> [l]^{\hspace{0cm}\imin\eta}}$$
inverse to each others.
\end{teo}

We need to introduce the subcategory of coherent conic sheaves.

\begin{df}\label{2.2.2} Let $U \in \op(X_{\RP})$. Then $U$ is said to be
relatively quasi-compact if, for any covering $\{U_i\}_{i \in I}$
of $X_{\RP}$, there exists $J \subset I$ finite such that $U
\subset \bigcup_{i \in J} U_i.$ We write $U \subset\subset X_{\RP}$.

We  will denote by $\op^c(X_{\RP})$ the subcategory of
$\op(X_{\RP})$ consisting of relatively quasi-compact open
subsets.
\end{df}

 One can check easily that if $U \in \op^c(X)$, then $\RP
U \in \op^c(X_{\RP})$.



\begin{df} Let $F \in \mod(k_{X_{\RP}})$ and consider the family $\op(X_{sa,\RP})$.
\begin{itemize}
\item[(i)] $F$ is $X_{sa,\RP}$-finite if there exists an epimorphism
$G \twoheadrightarrow F$, with $G \simeq \oplus_{i \in I}
k_{U_i}$,
 $I$ finite and $U_i \in \op^c(X_{sa,\RP})$.
\item[(ii)] $F$ is $X_{sa,\RP}$-pseudo-coherent if for any morphism
$\psi:G \to F$, where $G$ is $X_{sa,\RP}$-finite, $\ker \psi$ is
$X_{sa,\RP}$-finite.
\item[(iii)] $F$ is $X_{sa,\RP}$-coherent if it is both $X_{sa,\RP}$-finite
and $X_{sa,\RP}$-pseudo-coherent.
\end{itemize}
We will denote by $\coh(X_{sa,\RP})$ the subcategory of
$\mod(k_{X_{\RP}})$ consisting of $X_{sa,\RP}$-coherent objects.
\end{df}

\subsection{Conic sheaves on subanalytic sites}

\begin{df} A sheaf of $k$-modules $F$ on $X_{sa}$ is conic if the restriction morphism
$\Gamma(\RP U;F) \to \Gamma(U;F)$ is an isomorphism for each
$\RP$-connected $U \in \op^c(X_{sa})$ with $\RP U \in
\op(X_{sa})$.
\begin{itemize}
\item[(i)]We denote by $\mod_{\RP}(k_{X_{sa}})$ the subcategory of
$\mod(k_{X_{sa}})$ consisting of conic sheaves.
\item[(ii)] We denote by $D^b_{\RP}(k_{X_{sa}})$,
the subcategory of $D^b(k_{X_{sa}})$ consisting of objects $F$
such that $H^j(F)$ belongs to $\mod_{\RP}(k_{X_{sa}})$ for all $j
\in \Z$.
\end{itemize}
\end{df}

\begin{oss}\label{noequivariant} Let $X$ be a real analytic manifold endowed with a subanalytic
action $\mu$ of $\RP$ and consider the following diagram

$$\xymatrix{X \times \RP \ar@ <2pt>
[r]^{\hspace{0.3cm}\mu} \ar@ <-2pt> [r]_{\hspace{0.3cm}p}& X,}$$
where  $p$ denotes the projection. As in classical sheaf theory
one can define the subcategory $\mod^\mu(k_{X_{sa}})$ of
$\mod(k_{X_{sa}})$ consisting of sheaves satisfying $\imin \mu F
\simeq \imin p F$. The categories $\mod^\mu(k_{X_{sa}})$ and
$\mod_{\RP}(k_{X_{sa}})$ are not equivalent in general.

Indeed, let $X=\R$, set $X^+=\{x \in \R;\, x>0\}$ and let $\mu$ be
the natural action of $\RP$ (i.e. $\mu(x,t)=tx$). Let us consider
the sheaf $\rho_!k_{X^+} \in \mod(k_{X_{sa}})$. Then
$$\imin \mu \rho_! k_{X^+} \simeq \rho_! \imin \mu k_{X^+} \simeq \rho_! \imin p k_{X^+}
\simeq \imin p \rho_! k_{X^+}.$$ Let $V=\{x \in \R;\, 1<x<2\}$ and
set $W_m=\{x \in \R;\, {1\over m}<x<m\}$, where $m \in \N
\setminus \{0\}$. Recall that $\rho_!k_{X^+} \simeq \lind {U
\subset \subset X^+} \rho_*k_U \simeq \lind m \rho_* k_{W_m}$. We
have
$$\Gamma(V;\rho_!k_{X^+}) \simeq \lind m\Gamma(V;k_{W_m}) \simeq k,$$
since $V \subset W_m$ for $m \geq 2$. On the other hand, let
$V_n^+=\{x \in \R;\, 0<x<n\}$, where $n \in \N$. Since $\RP V=
X^+$ we have
$$\Gamma(X^+;\rho_!k_{X^+}) \simeq \lpro n \Gamma(V_n^+;\rho_!k_{X^+})
\simeq \lpro n \lind m \Gamma(V_n^+;k_{W_m}),$$ (in the second
isomorphism we used the fact that $V_n^+ \in \op^c(X_{sa})$ for
each $n$) and $\Gamma(V_n^+;k_{W_m})=0$ for each $n,m \in \N$.
Hence $\Gamma(V;\rho_!k_{X^+}) \not\simeq \Gamma(\RP
V;\rho_!k_{X^+})$.
\end{oss}


\begin{df}We denote by $\op(X_{sa,\RP})$ the full subcategory of $\op(X_{sa})$ consisting of conic subanalytic subsets, i.e. $U \in
\op(X_{sa,\RP})$ if $U \in \op(X_{sa})$
and it is invariant by the action of $\RP$.

We denote by $X_{sa,\RP}$ the category
$\op(X_{sa,\RP})$ endowed with the topology induced by $X_{sa}$.
\end{df}

We denote by $\rho_{\RP}:X_{\RP} \to X_{sa,\RP}$ the natural
morphism of sites.

Let  $\eta: X \to X_{\RP}$ and $\eta_{sa}: X_{sa} \to X_{sa,\RP}$
be the natural morphisms of sites. We have a commutative diagram
of sites
\begin{equation}\label{etarhoRP}
\xymatrix{X \ar[r]^{\rho} \ar[d]^{\eta} & X_{sa} \ar[d]^{\eta_{sa}} \\
X_{\RP} \ar[r]^{\rho_{\RP}} & X_{sa,\RP}.}
\end{equation}

\begin{lem}\label{etarho} Let $F \in \coh(X_{sa,\RP})$. Then
$\imin \eta_{sa}
\rho_{\RP*} F \simeq \rho_*\imin \eta F.$
\end{lem}
\dim\ \ Since all these functors are exact on $\coh(X_{sa,\RP})$,
we may reduce to the case $F=k_U$ with $U \in \op^c(X_{sa,\RP})$.
Then we have
$$\imin \eta_{sa} \rho_{\RP*} k_U \simeq \imin \eta_{sa} k_U \simeq
k_U,$$ on the other hand we have
$$\rho_*\imin \eta k_U \simeq
\rho_*k_U \simeq k_U$$ and the result follows. \nopagebreak
\newline \qed

Replacing $\op^c(X_{sa})$ with $\op^c(X_{sa,\RP})$, we can adapt the
results of \cite{KS01}, \cite{Pr1} and we get the following results.

\begin{teo}\label{RPind} (i) Let $G \in \coh(X_{sa,\RP})$ and let
$\{F_i\}$ be a filtrant inductive system in $\mod(k_{X_{sa,\RP}})$.
Then we have an isomorphism
$$\lind i \Ho_{k_{X_{sa,\RP}}}(\rho_{\RP *}G,F_i) \iso
\Ho_{k_{X_{sa,\RP}}}(\rho_{\RP *}G,\lind i F_i).$$
Moreover the functor of direct image $\rho_{\RP*}$ associated to
the morphism $\rho_{\RP}$ in \eqref{etarhoRP} is fully faithful
and exact on $\coh(X_{sa,\RP})$.

(ii) Let $F \in \mod(k_{X_{sa,\RP}})$. There exists a small filtrant
inductive system $\{F_i\}_{i \in I}$ in $\coh(X_{sa,\RP})$ such
that $F \simeq \lind i \rho_{\RP *}F_i$.
\end{teo}

\begin{nt} Since $\rho_{\RP *}$ is fully faithful and exact
on $\coh(X_{sa,\RP})$, we can identify
$\coh(X_{sa,\RP})$ with its image in $\mod(k_{X_{sa,\RP}})$. When there
is no risk of confusion we will write $F$ instead of $\rho_{\RP *}F$,
for $F \in \coh(X_{sa,\RP})$.\\
\end{nt}

We can also find a left adjoint to the functor $\imin {\rho_{\RP}}$.

\begin{prop} The functor $\imin {\rho_{\RP}}$ admits a left adjoint,
denoted by $\rho_{\RP !}$. It satisfies
\begin{itemize}
\item[(i)] for $F \in \mod(k_{X_{\RP}})$ and $U \in \op^c(X_{sa,\RP})$, $\rho_{\RP !} F$ is the sheaf associated to the presheaf $U
\mapsto \lind {V \supset\supset U} \Gamma(V;F)$, \\
\item[(ii)] for $U \in \op(X_{\RP})$ one has $\rho_{\RP !}k_U \simeq \lind{V\subset\subset U, V \in \op^c(X_{sa,\RP})
}k_V$.
\end{itemize}
\end{prop}

\begin{oss} One can check that $\imin {\eta_{sa}} \circ \rho_{\RP *} \simeq \rho_* \circ \imin \eta$ and $\imin \rho \circ \imin {\eta_{sa}} \simeq \imin \eta \circ \imin {\rho_{\RP}}$. Remark that $\rho_! \circ \imin \eta \not\simeq \imin {\eta_{sa}} \circ \rho_{\RP !}$. In fact with the notations of Remark \ref{noequivariant} we have $\rho_!\imin \eta k_{X^+} \simeq \lind m \rho_* k_{W_m}$. On the other hand, $\imin {\eta_{sa}}\rho_{\RP !}k_{X^+} \simeq k_{X^+}$.
\end{oss}

\noindent Let us extend the notion of quasi-injective object of \cite{KS01}, \cite{Pr1} to $\mod(k_{X_{sa,\RP}})$.

\begin{df}\label{dfqinj} An object $F \in \mod(k_{X_{sa,\RP}})$ is quasi-injective if the
functor $\Ho_{k_{X_{sa,\RP}}}(\cdot,F)$ is exact in
$\coh(X_{sa,\RP})$ or, equivalently (see Theorem 8.7.2 of
\cite{KS}) if for each $U,V \in \op^c(X_{sa,\RP})$ with $V \subset U$
the restriction morphism $\Gamma(U;F) \to \Gamma(V;F)$ is
surjective.
\end{df}

The category of quasi-injective object
is cogenerating since it contains injective objects. Moreover it is
stable by filtrant $\Lind$ and $\prod$. We have the following result

\begin{teo}\label{cohinj} The family of quasi-injective sheaves is injective with respect to
the functor $\Ho_{k_{X_{sa,\RP}}}(G,\cdot)$ for each $G \in
\coh(X_{sa,\RP})$.
\end{teo}

In particular

\begin{prop} The family of quasi-injective sheaves is injective with respect to the functor $\Gamma(U;\cdot)$ for any $U \in \op(X_{sa,\RP})$.
\end{prop}

\subsection{An equivalence of categories}\label{consa}

Let $X$ be a real analytic manifold endowed with an action $\mu$
of $\RP$. In the following we shall assume the hypothesis below:
\begin{equation}\label{hypsa}
  \begin{cases}
 \text{(i) every $U \in \op^c_{sa}(X)$ has a finite covering consisting }\\
 \text{\ \ of $\RP$-connected subanalytic open subsets,}\\
 \text{(ii) for any $U \in \op^c_{sa}(X)$ we have $\RP U \in \op(X_{sa})$,}\\
 \text{(iii) for any $x \in X$\ \ the set $\RP x$ is contractible,}\\
 \text{(iv) there exists a covering $\{V_n\}_{n \in \N}$ of $X_{sa}$ such that}\\
 \text{\ \ $V_n$ is $\RP$-connected and $V_n \subset\subset V_{n+1}$ for each $n$}.
  \end{cases}
\end{equation}

Let $U \in \op(X_{sa})$ such that $\RP U$ is still subanalytic.
Let $\varphi$ be the natural map from $\Gamma(\RP U;F)$ to
$\Gamma(U;\imin \eta_{sa} F)$ defined by
\begin{equation}\label{varphisa}
\begin{array}{ccc}
\Gamma(\RP U;F) & \to & \Gamma(\RP U; \eta_{sa*} \imin \eta_{sa} F) \\
  & \simeq & \Gamma(\RP U; \imin \eta_{sa} F) \\
 & \to & \Gamma(U; \imin \eta_{sa} F).
\end{array}
\end{equation}

\begin{prop}\label{*-1Usa} Let $F \in \mod(k_{X_{sa,\RP}})$.
Let $U \in \op(X_{sa})$
and assume that $U$ is $\RP$-connected. Then the morphism
$\varphi$ defined by \eqref{varphisa} is an isomorphism.
\end{prop}
\dim\ \ (i) Assume that $U \in \op^c(X_{sa})$ is $\RP$-connected. Let $F \in \mod(k_{X_{sa,\RP}})$,
then $F=\lind i \rho_{\RP*} F_i$, with $F_i \in
\coh(X_{sa,\RP})$. We have the chain of
isomorphisms
\begin{eqnarray*}
\Ho_{k_{X_{sa}}}(k_U,\imin \eta_{sa}\lind i \rho_{\RP*}F_i) &
\simeq & \Ho_{k_{X_{sa}}}(k_U,
\lind i \rho_* \imin \eta F_i)\\
& \simeq & \lind i \Ho_{k_X}(k_U, \imin \eta F_i)\\
& \simeq & \lind i \Ho_{k_{X_{\RP}}}(k_{\RP U}, F_i)\\
& \simeq &  \Ho_{k_{X_{sa,\RP}}}(k_{\RP U},\lind i \rho_{\RP *}
F_i),
\end{eqnarray*}
where the first isomorphism follows since $\imin \eta_{sa} \circ
\rho_{\RP*} \simeq \rho_* \circ \imin \eta$ by Lemma \ref{etarho} and the third one
follows from the equivalence between conic sheaves on $X$ and
sheaves on $X_{\RP}$. In the fourth isomorphism we used the fact
that $\RP U \in \op^c(X_{sa,\RP})$.

(ii) Let $U \in \op(X_{sa})$ be $\RP$-connected. Let $\{V_n\}_{n
\in \N} \in \cov(X_{sa})$ be a covering of $X$ as in \eqref{hypsa}
(iv) and set $U_n=U \cap V_n$. We have
\begin{equation}\label{ML0}
\Gamma(U;\imin \eta_{sa}F) \simeq \lpro n\Gamma(U_n;\imin \eta_{sa}F) \simeq \lpro n\Gamma(\RP U_n;F)
\simeq \Gamma(\RP U;F).
\end{equation}
\nopagebreak
\newline \qed

We can extend Lemma \ref{etarho} to $\mod(k_{X_{\RP}})$.

\begin{lem}\label{etarhomod} Let $F \in \mod(k_{X_{\RP}})$. Then
$\imin \eta_{sa}
\rho_{\RP*} F \simeq \rho_*\imin \eta F.$
\end{lem}
\dim\ \ Let $F \in \mod(k_{X_{\RP}})$ and let $U \in \op^c(X_{sa})$ be $\RP$-connected. Then
$$
\Gamma(U;\rho_*\imin \eta F) \simeq \Gamma(U;\imin \eta F) \simeq \Gamma(\RP U;F),
$$
where the second isomorphism follows from Proposition \ref{*-1U}. On the other hand
$$
\Gamma(U;\imin {\eta_{sa}}\rho_{\RP *}F) \simeq \Gamma(\RP U;\rho_{\RP *}F) \simeq \Gamma(\RP U;F),
$$
where the second isomorphism follows from Proposition \ref{*-1Usa}. Hence by \eqref{hypsa} (i) $\imin \eta_{sa}
\rho_{\RP*} F \simeq \rho_*\imin \eta F.$\\ \qed

Let us consider the category $\mod_{\RP}(k_{X_{sa}})$ of conic
sheaves on $X_{sa}$. The restriction of $\eta_{sa*}$ induces a
functor denoted by $\widetilde{\eta}_{sa*}$ and we obtain a
diagram
\begin{equation}\label{sarho}
\xymatrix{\mod_{\RP}(k_{X_{sa}})  \ar[d] \ar@ <2pt>
[rr]^{\widetilde{\eta}_{sa*}} &&
  \mod(k_{X_{sa,\RP}}) \ar@ <2pt> [ll]^{\imin \eta_{sa}} \\
\mod(k_{X_{sa}}) \ar[urr]_{\eta_{sa*}} && }
\end{equation}

\begin{teo}\label{ex*sa} The functors $\widetilde{\eta}_{sa*}$ and $\imin \eta_{sa}$ in \eqref{sarho}
  are equivalences of ca\-te\-go\-ries inverse to each
  others.
\end{teo}
\dim\ \ (i) Let $F \in \mod_{\RP}(k_{X_{sa}})$, and let $U \in
\op^c(X_{sa})$ be $\RP$-connected. We have
$$\Gamma(U;F) \simeq  \Gamma(\RP U;F) \simeq
\Gamma(\RP U;\widetilde{\eta}_{sa*}F) \simeq \Gamma(U;\imin
\eta_{sa} \widetilde{\eta}_{sa*}F).$$ The third isomorphism
follows from Proposition \ref{*-1Usa}. Then \eqref{hypsa} (i)
implies $\imin \eta_{sa} \widetilde{\eta}_{sa*} \simeq \id$.

(ii) For any $U \in \op^c(X_{sa,\RP})$ we have:
$$\Gamma(U;\eta_{sa*}\imin \eta_{sa} F) \simeq \Gamma(U;\imin \eta_{sa} F) \simeq \Gamma(U;F)$$ where
the second isomorphisms follows from Proposition \ref{*-1Usa}. This implies $\eta_{sa*} \imin \eta_{sa} \simeq \id$.\\
\qed

\begin{nt} Since $\imin {\eta_{sa}}$ is fully faithful and exact
we will often identify $\coh(X_{sa,\RP})$ with its image in
$\mod_{\RP}(k_{X_{sa}})$. Hence, for $F \in \coh(X_{sa,\RP})$ we
shall often write $F$ instead of $\imin {\eta_{sa}}F$.
\end{nt}

Thanks to Theorem \ref{RPind} we can give another description of the category of conic sheaves.

\begin{teo}\label{indRP} Let $F \in \mod_{\RP}(k_{X_{sa}})$. Then there exists a small filtrant system $\{F_i\}$ in $\coh(X_{sa,\RP})$ such that  $F \simeq \lind i \rho_* \imin \eta F_i$.
\end{teo}

\begin{oss} Let
$F \in \coh(X_{sa,\RP})$. The functor of inverse image commutes
with $\Lind$ and
$$\imin \mu \rho_* \imin \eta F \simeq  \rho_* \imin \mu \imin \eta F
\simeq \rho_* \imin p \imin \eta F \simeq \imin p \rho_* \imin
\eta F.$$ Hence $F \in \coh(X_{sa,\RP})$ implies $F \in
\mod^\mu(k_{X_{sa}})$, where $\mod^\mu(k_{X_{sa}})$ is the
category introduced in Remark \ref{noequivariant}. Since
$\mod^\mu(k_{X_{sa}})$ is stable by filtrant $\Lind$ we have that
$F$ belongs to $\mod^\mu(k_{X_{sa}})$.
Hence $\mod_{\RP}(k_{X_{sa}})$ is a full subcategory of
$\mod^\mu(k_{X_{sa}})$ but $\mod_{\RP}(k_{X_{sa}}) \not\simeq
\mod^\mu(k_{X_{sa}})$ in general. We have the chain of fully
faithful functors
$$\coh(X_{sa,\RP}) \hookrightarrow \mod_{\RP}(k_{X_{sa}})
\hookrightarrow \mod^\mu(k_{X_{sa}}).$$
\end{oss}

\subsection{Derived category}

Assume \eqref{hypsa}. Injective and quasi-injective objects of $\mod(k_{X_{sa}})$ are not contained in $\mod_{\RP}(k_{sa})$. 
For this reason we are going to introduce a subcategory which is useful when we try to find acyclic resolutions.


\begin{lem} Assume that $X$ satisfies \eqref{hypsa}. Then the following property is satisfied:
\begin{equation}\label{Dhypsa}
  \begin{cases}
    \text{ each finite covering of an $\RP$-connected $U \in \op^c(X_{sa})$ }\\
    \text{ has a finite refinement $\{V_i\}_{i=1}^n$ such that each ordered}\\
    \text{ union $\bigcup_{i=1}^jV_i$ is $\RP$-connected for each $j \in \{1,\ldots,n\}$}.
  \end{cases}
\end{equation}
\end{lem}
\dim\ \ Let $U \in \op^c(X_{sa})$ be $\RP$-connected. Then each finite covering of $U$ admits a finite refinement consisting of $\RP$-connected open subanalytic subsets. Let $\{U_i\}_{i=1}^n$ be a finite covering of $U$, $U_i \in \op^c(X_{sa})$ $\RP$-connected for each $i$. We will construct a refinement satisfying \eqref{Dhypsa}.

For $k=1,\dots,n$ set $V_{k11}:=U_k$ and $V_{k1i}:=U_{\sigma(i)} \cap \RP(U_k \cap U_{\sigma(i)})$ for $i=2,\dots,n$ and $\sigma(i)=i-1$ if $i\leq k$, $\sigma(i)=i$ if $i>k$. Then set $U_{k2}:=\bigcup_{i=1}^nV_{k1i}$ and $V_{k2i}:=U_{\sigma(i)} \cap \RP(U_{k2} \cap U_{\sigma(i)})$. For $j=1,\dots, n$ define recursively $U_{kj}=\bigcup_{\ell=1}^j\bigcup_{i=1}^nV_{k\ell i}$ and $V_{kji}=U_{\sigma(i)} \cap \RP(U_{kj} \cap U_{\sigma(i)})$. Remark that 
$\bigcup_{p=1}^j\bigcup_{\ell=1}^n\bigcup_{i=1}^nV_{p\ell i}= \bigcup_{p=1}^j\RP U_p \cap U$. By Lemma \ref{RPconn} below all the sets $V_{kji}$ are $\RP$-connected and $\{V_{kji}\}_{i,k,j}$ is a refinement of $\{U_i\}_i$ satisfying \eqref{Dhypsa} (with the lexicographic order). \\
\qed

\begin{lem}\label{RPconn} Assume that $X$ satisfies \eqref{hypsa} (iii). Let $U,V,W$ be open and $\RP$-connected. Then $U \cup (V \cap \RP(U \cap V)) \cup (W \cap \RP(U \cap W))$ is $\RP$-connected.
\end{lem}
\dim\ \ In what follows, when we write $\RP x$ we suppose that $\RP x \simeq \R$. If $\RP x = x$ everything becomes obvious.

(i) First remark that $U \cap V$ (resp. $U \cap W$, $V \cap W$) is $\RP$-connected. Indeed, let $x_1 \in U \cap \RP x$, $x_2 \in V \cap \RP x$ for some $x \in X$. Then $x_1=\mu(x,a)$, $x_2=\mu(x,b)$.  Every path in $\RP x$ connecting $x_1$ and $x_2$ contains $\mu(x,[a,b])$. Since $U$ and $V$ are $\RP$-connected then $U \cap V \supset \mu(x,[a,b])$.

(ii) Now let us prove that $U \cup (V \cap \RP(U \cap V))$ is $\RP$-connected. Let $x_1,x_2 \in U \cup (V \cap \RP(U \cap V)) \cap \RP x$ for some $x \in X$. Then $x_1=\mu(x,a)$, $x_2=\mu(x,b)$. We want to prove that $\mu(x,[a,b]) \subset U \cup (V \cap \RP(U \cap V))$. If $x_1,x_2 \in U$ it follows since $U$ is $\RP$-connected and if $x_1,x_2 \in V \cap \RP(U \cap V)$ it follows from (i). So we may assume that $x_1 \in U$ and $x_2 \in V \cap \RP(U \cap V)$. Since $U$ is $\RP$-connected and $x_2 \in \RP x_1$, there exists $y=\mu(x,c) \in U \cap V$. Then $\mu(x,[a,c]) \subset U$. In the same way $\mu(x,[b,c]) \subset V \cap \RP(U \cap V)$ and hence $\mu(x,[a,c] \cup [b,c]) \subset U \cup (V \cap \RP(U \cap V))$.

(iii) Let us show that $U \cup (V \cap \RP(U \cap V)) \cup (W \cap \RP(U \cap W))$ is $\RP$-connected. Let $x_1,x_2 \in U \cup (V \cap \RP(U \cap V)) \cup (W \cap \RP(U \cap W)) \cap \RP x$ for some $x \in X$. Then $x_1=\mu(x,a)$, $x_2=\mu(x,b)$. We want to prove that $\mu(x,[a,b]) \subset U \cup (V \cap \RP(U \cap V)) \cup (W \cap \RP(U \cap W))$. By (i) and (ii) we may reduce to the case $x_1 \in V$, $x_2 \in W$. As in (ii), there exist $y_1=\mu(x,c) \in U \cap V$ and $y_2=\mu(x,d) \in U \cap W$. Then $\mu(x,[c,d]) \in U$, $\mu(x,[a,c]) \subset V \cap \RP(U \cap V)$ and $\mu(x,[b,d]) \subset W \cap \RP(U \cap W)$. Hence $\mu(x,[c,d] \cup [a,c] \cup [b,d]) \in U \cup (V \cap \RP(U \cap V)) \cup (W \cap \RP(U \cap W))$ and the result follows.\\
\qed

\begin{df} A sheaf $F \in \mod(k_{X_{sa}})$ is $\RP$-quasi-injective if for each $\RP$-connected $U \in \op^c(X_{sa})$ the restriction morphism $\Gamma(X;F) \to \Gamma(U;F)$ is surjective.
\end{df}


Remark that the functor $\imin {\eta_{sa}}$ sends quasi-injective objects of $\mod(k_{X_{sa,\RP}})$ to $\RP$-quasi-injective objects since $\Gamma(U;\imin {\eta_{sa}}F) \simeq \Gamma(\RP U;F)$ if $U \in \op^c(X_{sa})$ is $\RP$-connected.
Moreover the category of $\RP$-injective objects is cogenerating since injective objects are cogenerating in $\mod(k_{X_{sa}})$.

\begin{prop} \label{conflU} Let  $\exs{F'}{F}{F''}$ be an exact sequence in
$\mod(k_{X_{sa}})$ and assume that $F'$ is $\RP$-quasi-injective. Let
$U \in \op^c(X_{sa})$ be $\RP$-connected. Then the sequence
$$\exs{\Gamma(U;F')}{\Gamma(U;F)}{\Gamma(U;F'')}$$
is exact.
\end{prop}
\dim\ \ Let $s'' \in \Gamma(U;F'')$, and let $\{V_i\}_{i=1}^n$ be
a finite covering of $U$ satisfying \eqref{Dhypsa} and such that
there exists $s_i \in \Gamma(V_i;F)$ whose image is $s''|_{V_i}$.
For $n \geq 2$ on $V_1 \cap V_2$ $s_1-s_2$ defines a section of
$\Gamma(V_1 \cap V_2;F')$ which extends to $s' \in \Gamma(X;F')$.
Replace $s_1$ with $s_1-s'$. We may suppose that $s_1=s_2$ on $V_1
\cap V_2$. Then there exists $t \in \Gamma(V_1 \cup V_2)$ such
that $t|_{V_i}=s_i$, $i=1,2$. Thus the induction proceeds.
\nopagebreak \newline \qed

\begin{prop}\label{RPqinjinj} Let  $\exs{F'}{F}{F''}$ be an exact sequence in
$\mod(k_{X_{sa}})$ and assume that $F'$ is $\RP$-quasi-injective. Let
$U \in \op(X_{sa})$ be $\RP$-connected. Then the sequence
$$\exs{\Gamma(U;F')}{\Gamma(U;F)}{\Gamma(U;F'')}$$
is exact.
\end{prop}
\dim\ \ By \eqref{hypsa} (iv) there exists a covering $\{V_n\}_{n
\in \N}$ of $X_{sa}$ such that $V_n$ is $\RP$-connected and $V_n
\subset\subset V_{n+1}$ for each $n$. For each $n$ the sequence
$$\exs{\Gamma(U\cap V_n;F')}{\Gamma(U\cap V_n;F)}{\Gamma(U\cap V_n;F'')}$$
is exact and the morphism $\Gamma(U\cap V_{n+1};F') \to
\Gamma(U\cap V_n;F)$ is surjective for each $n$ since $F'$ is
$\RP$-quasi-injective. Then by the Mittag-Leffler property (see
Proposition 1.12.3 of \cite{KS90}) the sequence
$$\exs{\lpro n\Gamma(U\cap V_n;F')}{\lpro n\Gamma(U\cap V_n;F)}{\lpro n\Gamma(U\cap V_n;F'')}$$
is exact. Since $\lpro n\Gamma(U \cap V_n;G) \simeq \Gamma(U;G)$
for each $G \in \mod(k_{X_{sa}})$ the result follows. \nopagebreak
\newline \qed

\begin{prop} Let $F',F$ be $\RP$-quasi-injective and consider the exact
sequence  $\exs{F'}{F}{F''}$ in $\mod(k_{X_{sa}})$. Then $F''$ is
$\RP$-quasi-injective.
\end{prop}
\dim\ \ Let $U \in \op^c(X_{sa})$ be $\RP$-connected and let
us consider the diagram below
$$ \xymatrix{\Gamma(X;F) \ar[d]^\alpha \ar[r] & \Gamma(X;F'') \ar[d]^\gamma \\
\Gamma(U;F) \ar[r]^\beta & \Gamma(U;F'').} $$  The morphism
$\alpha$ is surjective since $F$ is $\RP$-quasi-injective and $\beta$
is surjective by Proposition \ref{conflU}. Then $\gamma$ is
surjective. \nopagebreak \newline \qed

It follows from the preceding results that

\begin{prop} \label{RPqinjinj} $\RP$-quasi-injective objects are injective with respect to the functor $\Gamma(U;\cdot)$, with $U \in \op(X_{sa})$ and $\RP$-connected.
\end{prop}

\begin{cor} $\RP$-quasi-injective objects are $\eta_{sa*}$-injective.
\end{cor}
\dim\ \ Let
$\exs{F'}{F}{F''}$ be an exact sequence in
$\mod_{\RP}(k_{X_{sa}})$ and assume that $F'$ is $\RP$-quasi-injective. By Proposition \ref{RPqinjinj} the sequence
$$\exs{\Gamma(U;F')}{\Gamma(U;F)}{\Gamma(U;F'')}$$
for any $U \in \op(X_{sa,\RP})$. This implies that the sequence
$$
\exs{\eta_{sa*}F'}{\eta_{sa*}F}{\eta_{sa*}F''}
$$
is exact. \\ \qed

\begin{teo}\label{Rex*sa} The categories $D^b(k_{X_{sa,\RP}})$ and
$D^b_{\RP}(k_{X_{sa}})$ are equivalent.
\end{teo}
\dim\ \ In order to prove this statement, it is enough to show
that $\imin \eta_{sa}$ is fully faithful. Let $F \in D^b(k_{X_{sa,\RP}})$ and let $F'$ be an injective complex quasi-isomorphic to $F$. Since $\imin {\eta_{sa}}$ sends injective objects to $\RP$-quasi-injective objects which are $\eta_{sa*}$-injective we have
$R\eta_{sa*} \imin {\eta_{sa}} F \simeq \eta_{sa*} \imin {\eta_{sa}} F'
\simeq F' \simeq  F$.
This implies $R\eta_{sa*} \imin \eta_{sa} \simeq \id$, hence $\imin {\eta_{sa}}$ is fully faithful.\\
\qed

Hence for each
$F \in D^b_{\RP}(k_{X_{sa}})$ we have  $F \simeq \imin
{\eta_{sa}}F'$ with $F' \in D^b(k_{X_{sa,\RP}})$.

\begin{prop}\label{preceding} Let $F \in D^b(k_{X_{sa,\RP}})$ and let $U \in \op(X_{sa})$ be $\RP$-connected. There is an isomorphism
$\mathrm{R}\Gamma(\RP U;F) \iso \mathrm{R}\Gamma(U; \imin {\eta_{sa}}F)$.
\end{prop}
\dim\ \ Let $F'$ be a complex of injective objects quasi-isomorphic to $F$. Since $\imin {\eta_{sa}}$ sends injective objects to $\RP$-quasi-injective objects we have
\begin{eqnarray*}
\mathrm{R}\Gamma(\RP U;F) & \simeq & \Gamma(\RP U;F') \\
& \iso & \Gamma(U; \imin {\eta_{sa}}F') \\
& \simeq & \mathrm{R}\Gamma(U; \imin {\eta_{sa}}F),
\end{eqnarray*}
where the second isomorphism follows from Proposition \ref{*-1Usa}.\\
\qed

We extend Lemma \ref{etarhomod} to $D^b(k_{X_{\RP}})$.

\begin{lem}\label{etarhodb} Let $F \in D^b(k_{X_{\RP}})$. Then
$\imin \eta_{sa}
R\rho_{\RP*} F \simeq R\rho_*\imin \eta F.$
\end{lem}
\dim\ \ (i) Let $F \in \mod(k_{X_{\RP}})$ be injective. Then for each $\RP$-connected $U \in \op^c(X_{sa})$ $\mathrm{R}\Gamma(U;R\rho_*\imin \eta F) \simeq \mathrm{R}\Gamma(U;\imin \eta F) \simeq \mathrm{R}\Gamma(\RP U;F)$ is concentrated in degree zero. Hence $\imin \eta F$ is $R\rho_*$-acyclic by \eqref{hypsa} (i).

(ii) Let $F \in D^b_{\RP}(k_X)$ and let $F'$ be a complex of injective objects quasi-isomorphic to $F$. Then $\imin \eta_{sa}
R\rho_{\RP*} F \simeq \imin \eta_{sa}
\rho_{\RP*} F' \simeq \rho_*\imin \eta F' \simeq R\rho_*\imin \eta F,$ where the second isomorphism follows from Lemma \ref{etarhomod} and the third one follows from (i).\\ \qed

\subsection{Operations}

Let $X$ be a real analytic manifold endowed with an analytic
action of $\RP$. In this section we study the operations in the
category of conic sheaves on $X_{sa}$.

\begin{prop} Let $F \in \mod_{\RP}(k_X)$ and let $G \in
\mod_{\RP}(k_{X_{sa}})$.
\begin{itemize}
\item[(i)] we have $\rho_*F \in \mod_{\RP}(k_{X_{sa}})$,
\item[(ii)] we have $\imin \rho G \in \mod_{\RP}(k_X)$.
\end{itemize}
\end{prop}
\dim\ \ (i) Let $U \in \op^c(X_{sa})$ be $\RP$-connected. We have
the chain of isomorphisms
$$ \Gamma(U;\rho_*F) \simeq \Gamma(U;F) \simeq \Gamma(\RP U;F)
\simeq \Gamma(\RP U;\rho_*F).$$

(ii) We have $G=\lind j \rho_*G_j$, with $G_j \in
\coh(X_{sa,\RP})$. Then
$$\imin \rho \lind j \rho_*G_j \simeq \lind j \imin \rho \rho_*G_j
\simeq \lind j G_j$$ and $\lind j G_j$ belongs to
$\mod_{\RP}(k_X)$. \nopagebreak \newline \qed

\begin{prop} Let $F,G \in \mod_{\RP}(k_{X_{sa}})$. Then
\begin{itemize}
\item[(i)] we have $F \otimes G \in \mod_{\RP}(k_{X_{sa}})$,
\item[(ii)] we have $\ho(F,G) \in \mod_{\RP}(k_{X_{sa}})$.
\end{itemize}
\end{prop}
\dim\ \ We have $F=\lind i \rho_*F_i$ and $G=\lind j \rho_*G_j$,
with $F_i,G_j \in \coh(X_{sa,\RP})$.

(i) we have $F \otimes G \simeq \lind {i,j} \rho_* (F_i \otimes
G_j)$ and $F_i \otimes G_j$ belongs to $\coh(X_{sa,\RP})$ for each
$i,j$.

(ii) we have $\ho(F,G) \simeq \lpro i \lind j \rho_*\ho(F_i,G_j)$
and $\ho(F_i,G_j)$ is conic for each $i,j$.
\nopagebreak \newline \qed

Let $f:X \to Y$ be a conic morphism of real analytic manifolds.

\begin{prop} Let $F \in \mod_{\RP}(k_{X_{sa}})$ and let $G \in
\mod_{\RP}(k_{Y_{sa}})$.
\begin{itemize}
\item[(i)] we have $f_*F \in \mod_{\RP}(k_{Y_{sa}})$,
\item[(ii)] we have $\imin f G \in \mod_{\RP}(k_{X_{sa}})$.
\end{itemize}
\end{prop}
\dim\ \ (i) Let $U \in \op^c(Y_{sa})$ be $\RP$-connected. Since
$f$ commutes with the action of $\RP$, the set $\imin f(U)$ is
$\RP$ connected. We have the chain of isomorphisms
$$\Gamma(\imin f(U);F) \simeq \Gamma(\RP \imin f(U);F)
\simeq \Gamma(\imin f(\RP U);F).$$ Hence $\Gamma(U;f_*F) \simeq
\Gamma(\RP U;f_*F)$.

 (ii) We have $G=\lind j \rho_*G_j$, with $G_j \in
\coh(Y_{sa,\RP})$. Then
$$\imin f \lind j \rho_*G_j \simeq \lind j \imin f \rho_*G_j
\simeq \lind j \rho_* \imin f G_j$$ and $\imin f G_j$ is conic
for each $j$. \nopagebreak
\newline \qed

Now let us consider the operations in the derived category of
conic subanalytic sheaves.

\begin{prop} Let $F \in \coh(X_{sa,\RP})$ and let $G \in
D^b_{\RP}(k_{X_{sa}})$. Then $\rh(F,G) \in D^b_{\RP}(k_{X_{sa}})$.
\end{prop}
\dim\ \ We may reduce to the case $G \in \mod_{\RP}(k_{X_{sa}})$.
Then $G=\lind j \rho_*G_j$, with $G_j \in \coh(X_{sa,\RP})$. Then by Proposition 2.2.2 of \cite{Pr1}
$$R^k\ho(F,G) \simeq \lind i \rho_*R^k\ho(F,G_j)$$
for each $k \in \Z$ and the result follows since $R^k\ho(F,G_j)$
is conic for each $k \in \Z$ and for each $j$. \nopagebreak
\newline \qed

\begin{prop} Let $F \in D^b_{\RP}(k_{X_{sa}})$. Then $Rf_*F \in
D^b_{\RP}(k_{Y_{sa}})$.
\end{prop}
\dim\ \ Remark that the functor $\imin {\eta_{sa}}$ sends
injective sheaves to $f_*$-acycic sheaves. This is a consequence
of the fact that $\imin {\eta_{sa}}$ sends injective sheaves to
$\Gamma(U;\cdot)$-acyclic sheaves for each $\RP$-connected $U \in
\op(X_{sa})$. There exists $F' \in D^b(k_{X_{sa,\RP}})$ such that
$F \simeq \imin {\eta_{sa}}F'$. Let $I^\bullet$ be a bounded
injective resolution of $F'$. Then $\imin {\eta_{sa}}I^j$ is conic
and $f_*$-acyclic for each $j$. We have $Rf_*F \simeq f_*\imin
{\eta_{sa}}I^\bullet$ and $f_*\imin {\eta_{sa}}I^j$ is conic for
each $j$. \nopagebreak
\newline \qed

\begin{prop} Let $G \in D^b_{\RP}(k_{Y_{sa}})$. Then $f^!G \in D^b_{\RP}(k_{X_{sa}})$.
\end{prop}
\dim\ \ We may reduce to the case $G \in \mod_{\RP}(k_{Y_{sa}})$. Then
$G=\lind j \rho_*G_j$, with $G_j \in \coh(Y_{sa,\RP})$. By Proposition 2.4.5 of \cite{Pr1}  we
have $H^kf^!G \simeq \lind i \rho_*H^kf^!G_j$ for each $k \in \Z$
and the result follows since $H^kf^!G_j$ is conic for each $k \in
\Z$ and for each $j$. \nopagebreak \newline \qed

\begin{oss}\label{f!!noconic} The functor $f_{!!}: \mod(k_{X_{sa}}) \to \mod(k_{Y_{sa}})$
 does not send conic sheaves to conic sheaves in general. In fact,
 let $p:\R^3 \to \R^2$ be the projection. It is a conic map with respect
 to the natural action of $\RP$ on $\R^3$ and $\R^2$. Set
\begin{eqnarray*}
U & = & \{(x,y) \in \R^2;\; (x-1)^2+y^2<1\}, \\
B_n & = & \{(x,y) \in \R^2;\; x^2+y^2<n\}, \\
B_n^+ & = & B_n \cap (\RP \times \R), \\
S & = & \RP(\partial U \times \{1\}).
\end{eqnarray*}
Let us consider the conic sheaf $k_S$. By definition of proper
direct image we have $\Gamma(U;p_{!!}F)=\lind K \Gamma(\imin p(U);
\Gamma_K F)$, where $K$ ranges through the family of subanalytic
compact subsets of $\R^3$. Since $U$ is bounded we have
$$\Gamma(U;p_{!!}k_S) \simeq \lind K \Gamma(U \times \R;\Gamma_Kk_S) \simeq
\lind m \Gamma(U \times \R;\Gamma_{\R^2 \times [-m,m]}k_S) \simeq
k,
$$
where $m \in \N$. On the other hand we have
$$
\Gamma(\RP U;p_{!!}k_S) \simeq \lpro n\Gamma(B_n^+;p_{!!}k_S)
\simeq \lpro n\lind m \Gamma(B_n^+ \times \R;\Gamma_{\R^2 \times
[-m,m]}k_S)=0,
$$
where $m,n \in \N$, since $\Gamma(B_n^+ \times \R;\Gamma_{\R^2
\times [-m,m]}k_S) =0$ for each $m,n$.
\end{oss}

Hence we shall need a new definition of proper direct image for conic sheaves.

\begin{df} We define functor $f_{\RP !!}$ of proper direct image for conic sheaves in the following way
\begin{eqnarray*}
f_{\RP !!}: \mod_{\RP}(k_{X_{sa}}) & \to & \mod_{\RP}(k_{Y_{sa}}) \\
\lind i \rho_* F_i & \mapsto & \lind i \rho_* f_! F_i,
\end{eqnarray*}
where $F_i \in \coh(X_{sa,\RP})$.
\end{df}

Note that if $F \in \coh(X_{sa,\RP})$ then $f_{\RP !!}\rho_* F \simeq \rho_*f_!F \not\simeq f_{!!}\rho_*F$. Moreover this definition is compatible with the classical one. In fact $f_{\RP !!}$ commutes with $\imin \rho$ and we have the following commutative diagram
$$
\xymatrix{\mod_{\RP}(k_X) \ar[d]^{\rho_*} \ar[r]^{f_!} & \mod_{\RP}(k_Y) \\
\mod_{\RP}(k_{X_{sa}}) \ar[r]^{f_{\RP !!}} & \mod_{\RP}(k_{Y_{sa}}). \ar[u]^{\imin \rho}}
$$

\begin{oss} With the notation of Remark \ref{f!!noconic}, we have
$$
\Gamma(\RP U;p_{\RP !!}k_S) \simeq \Gamma(U;p_{\RP !!}k_S) \simeq k .
$$
In fact the restriction of $p$ to $S \cap \{(x,y,z) \in \R^3;\; x>0\}$ is proper.
\end{oss}

Let us see an explicit formula for the sections of $f_{\RP !!}$. Let $U \in \op^c(X_{sa,\RP})$ and let $F=\lind i \rho_*F_i$ with $F_i \in \coh(X_{sa,\RP})$. We have the chain of isomorphisms
\begin{eqnarray*}
\Gamma(U;\lind i \rho_*f_!F_i) & \simeq & \lind i \Gamma(U;f_!F_i) \\
& \simeq & \lind {i,Z,K} \Gamma(\imin f(U);\Gamma_{Z \cap K}F_i) \\
& \simeq & \lind {i,Z',K} \Gamma(\imin f(U);\Gamma_{Z'\cap K}F_i) \\
& \simeq & \lind {Z',K}\Gamma(\imin f(U);\Gamma_{Z'\cap K}\lind i \rho_*F_i).
\end{eqnarray*}
Here $Z$ ranges into the family of closed subanalytic subsets of $\imin f(U)$ such that $f:Z \to U$ is proper, $Z'$ ranges through the family of closed conic subanalytic subsets of $\imin f(U)$ such that $\imin f(y) \cap \RP x=\{{\rm point}\}$ for any $y\in Y$, $x\in X$, and $K \subset\subset X_{\RP}$ are conic and closed. The first isomorphism follows since $U \in \op^c(X_{sa,\RP})$, the third since $F_i$ is conic for each $i$ and the last one since $\imin f(U) \cap K \subset\subset X_{\RP}$.\\

It is easy to prove that projection formula and base change formula for conic sheaves are satisfied. Moreover, $\RP$-quasi-injective objects are acyclic with respect to the functor $f_{\RP !!}$, since they are $\ho(G,\cdot)$-injective for each $G \in \coh(X_{sa,\RP})$.\\

In order to find a right adjoint to $Rf_{\RP !!}$ we follow the method used to find a right adjoint to the functor proper direct image for subanalytic sheaves. We shall skip the details of the proof (which are an adaptation of the results of \cite{Pr1}). The subcategory $\mathcal{J}_{X_{sa,\RP}}$ of $\RP$-quasi-injective objects
and the functor $f_{\RP !!}$ have the following properties:
\begin{equation*}
  \begin{cases}
    \text{(i) $\mathcal{J}_{X_{sa,\RP}}$ is cogenerating}, \\
    \text{(ii) $\mod(k_{X_{sa,\RP}})$ has finite quasi-injective dimension}, \\
    \text{(iii) $\mathcal{J}_{X_{sa,\RP}}$ is $f_{\RP !!}$-injective}, \\
    \text{(iv) $\mathcal{J}_{X_{sa,\RP}}$ is closed by small $\oplus$}, \\
    \text{(v) $f_{\RP !!}$ commutes with small $\oplus$}.
  \end{cases}
  \end{equation*}
As a consequence of the Brown representability theorem (see
\cite{KS}, Corollary 14.3.7 for details) we find a right adjoint
to the functor $Rf_{\RP !!}$, denoted by $f^!_{\RP}$. By adjunction $f^!_{\RP}$ commutes with $R\rho_*$ and as in \cite{Pr1} one can prove that $H^kf^!_{\RP}$ commutes with filtrant $\Lind$. Hence $f^!_{\RP}$ coincides with the restriction of $f^!$ to $D^b_{\RP}(k_{Y_{sa}})$.


\section{Fourier-Sato transform for subanalytic sheaves.}

In this section we construct the Fourier-Sato transform for
subanalytic sheaves. Reference are made to \cite{KS90} for the
classical Fourier-Sato transform.

\subsection{Conic sheaves on vector bundles}

Let $E \stackrel{\tau}{\to} Z$ be a real vector bundle, with
dimension $n$ over a real analytic manifold $Z$. Then $\RP$ acts
naturally on $E$ by multiplication on the fibers. We identify $Z$
with the zero-section of $E$ and denote by $i:Z \hookrightarrow E$
the embedding. We set $\pE=E \setminus Z$ and $\pittau:\pE \to
Z$ denotes the projection.

\begin{lem}\label{bierstone} The category $\op(E_{sa})$ satisfies
\eqref{hypsa}.
\end{lem}
\dim\ \ Let us prove \eqref{hypsa} (i). Let $\{V_i\}_{i \in \N}$
be a locally finite covering of $Z$ with $V_i \in \op^c(Z_{sa})$
such that $\imin \tau (V_i) \simeq \R^m \times \R^n$ and let
$\{U_i\}$ be a refinement of $\{V_i\}$ with $U_i \in
\op^c(Z_{sa})$ and $\overline{U_i} \subset V_i$ for each $i$. Then
$U$ is covered by a finite number of $\imin \tau (U_i)$ and $U
\cap \imin \tau(U_i)$ is relatively compact in $\imin \tau (V_i)$
for each $i$. We may reduce to the case $E \simeq \R^m \times
\R^n$.
 Let us consider the morphism
of manifolds
\begin{eqnarray*}
\varphi : \R^m \times \mathbb{S}^{n-1} \times \R & \to & \R^m \times \R^n\\
(z,\vartheta,r) & \mapsto & (z,r i(\vartheta)),
\end{eqnarray*}
where $i:\mathbb{S}^{n-1} \hookrightarrow \R^n$ denotes the
embedding. Then $\varphi$ is proper and subanalytic. The subset
$\imin \varphi(U)$ is subanalytic and relatively compact in $\R^m
\times  \mathbb{S}^{n-1} \times \R$. \\


(a) By Lemma \ref{goodcover} $\imin\varphi(U \setminus Z)$ admits a finite cover $\{W_j\}_{j \in J}$ such that the intersections of each $W_j$ with the fibers of $\pi:\R^m \times \mathbb{S}^{n-1} \times \R \to \R^m \times \mathbb{S}^{n-1}$ are contractible or empty. Then $\varphi(W_j)$ is an open subanalytic relatively compact $\RP$-connected subset of $\R^m \times \R^n$ for each $j$. In this way we obtain a finite covering of $U \setminus Z$ consisting of $\RP$-connected subanalytic open subsets.



(b) Let $p \in \pi(\imin\varphi(U \cap Z))$. Then $\imin \pi(p) \cap U$ is a disjoint union of intervals. Let us consider the interval $(m(p),M(p))$, $m(p)<M(p) \in \R$ containing $0$. Set $W_Z=\{(p,r) \in U ;\; m(p)<r<M(p)\}$. 
The set $W_Z$
is open subanalytic (it is a consequence of Proposition 1.2, Chapter 6 of \cite{VD98}), contains $\imin\varphi(U \cap Z)$ and its intersections with the fibers of $\pi$ are contractible. Then $\varphi(W_Z)$ is an open $\RP$-connected subanalytic neighborhood of $U \cap Z$ and it is contained in $U$.

By (a) there exists a finite covering $\{\varphi(W_j)\}_{j \in J}$ of $U
\setminus Z$ consisting of $\RP$-connected subanalytic open
subsets, and $\varphi(W_Z) \cup \bigcup_{j \in J}\varphi(W_j)=U$.\\

By Proposition 8.3.8 of \cite{KS90} the category $\op(E_{sa})$
also satisfies \eqref{hypsa} (ii). Moreover \eqref{hypsa} (iii)
and (iv) are clearly satisfied. \nopagebreak \newline \qed

Now let us consider $E$ endowed with the conic topology. In this situation, an object $U \in \op(E_{\RP})$ is the union of $\pU \in \op(\pE_{\RP})$ and $U_Z \in \op(Z)$ such that $\imin \pittau (U_Z) \subset \pU$. If $U,V \in \op(E_{\RP})$, then $U \subset\subset V$ if $U_Z \subset\subset V_Z$ in $Z$ and $\pU \subset\subset \pV$ in $\pE_{\RP}$ (this means that $\pi(\pU) \subset\subset \pi(\pV)$ in $\pE/\RP$, where $\pi:\pE \to \pE/\RP$ denotes the projection). \\


Applying Theorem \ref{Rex*sa} we have the following

\begin{teo} The categories $D^b_{\RP}(k_{E_{sa}})$ and $D^b(k_{E_{sa,\RP}})$ are
equi\-valent.
\end{teo}

Consider the subcategory $\mod^{cb}_{\rc,\RP}(k_E)$ of
$\mod_{\rc,\RP}(k_E)$ consisting of sheaves whose support is
compact on the base (i.e. $\tau(\supp(F))$ is compact in $Z$).

Let us consider the natural map $\eta: E \to E_{\RP}$. The
restriction of $\imin \eta$ to $\coh(E_{sa,\RP})$ gives rise to
a functor

\begin{equation}\label{rhosa}
\imin {\overline{\eta}}:\coh(E_{sa,\RP}) \to
\mod^{cb}_{\rc,\RP}(k_E)
\end{equation}

Since the functor $\imin \eta$ is fully faithful and exact, we
identify $\coh(E_{sa,\RP})$ as a subcategory of
$\mod^{cb}_{\rc,\RP}(k_E)$.

\begin{teo} \label{Rexsa*} The functor $\imin {\overline{\eta}}$ in \eqref{rhosa} is an
equivalence of categories.
\end{teo}
\dim\ \ (i) Let $F \in \mod^{cb}_{\rc,\RP}(k_E)$. Let us show that
$F$ is $E_{sa,\RP}$-finite. We may reduce to the case $E \simeq
\R^m \times \R^n$ and $Z \simeq \R^m \times \{0\}$.
 It is well
known that if $X$ is a real analytic manifold and $G \in
\mod^c_{\rc}(k_X)$, then $G$ is quasi-isomorphic to a bounded
complex of finite sums $\oplus_W k_W$, where $W \in
\op^c_{sa}(X)$.\\

Let us consider the diagram $Z \stackrel{i}{\hookrightarrow} E
\stackrel{\tau}{\to} Z,$ where $i$ is the embedding. We have
$\tau_*F \in \mod^c_{\rc}(k_Z). $ We have an exact sequence
\begin{equation}\label{exZ}
\oplus_{i \in I}  k_{\imin \tau (V_i)} \to \imin \tau\tau_*F \to
0,
\end{equation}
where $I$ is finite and $V_i \in \op_{sa}^c(Z)$.

 Now let us
consider the diagram $S \stackrel{j}{\hookrightarrow} \pE
\stackrel{\pi}{\to} S,$ where $S=\pE/\RP \simeq \R^m \times
\mathbb{S}^{n-1}$ and $\pi$ is the projection. We have $\imin j
F|_{\pE} \in \mod^c_{\rc}(k_S). $ Since $F|_{\pE}$ is conic $\imin
\pi \imin j F|_{\pE} \simeq F|_{\pE}$. We have an exact sequence
\begin{equation}\label{exS}
\oplus_j k_{\imin \pi (U_j)} \to F_{\pE} \to 0.
\end{equation}
where $J$ is finite and  $U_j \in \op^c_{sa}(S)$.

It is easy to check that the morphism $\imin \tau \tau_*F \oplus
F_{\pE} \rightarrow F$ is an epimorphism and we obtain the
result by (\ref{exZ}) and (\ref{exS}).\\

(ii) Let us show that $F$ is $E_{sa,\RP}$-pseudo-coherent. Let
$G=\oplus_{i \in I} k_{W_i}$, with $I$ finite and $W_i \in
\op^c(E_{sa,\RP})$, and consider a morphism $\psi:G \to F$. Since
$F$ and $ G$ are $\R$-constructible and conic, then $\ker \psi$
belongs to $\mod_{\rc,\RP}(k_E)$, and its support is still compact
on the base. \nopagebreak \newline \qed

As a consequence of Theorems \ref{indRP} and \ref{Rexsa*} one has
the following

\begin{teo} Let $F \in \mod_{\RP}(k_{E_{sa}})$. Then there exists a small filtrant system $\{F_i\}$ in $\mod_{\rc,\RP}^{cb}(k_E)$ such that  $F \simeq \lind i \rho_*F_i$.
\end{teo}



We end this section with the following result, which will be
useful in $\S$\,\ref{fouriersato}.

\begin{lem}\label{itausa} Let $F \in D^b_{\RP}(k_{E_{sa}})$. Then:
\begin{itemize}
\item[(i)] $R\tau_*F \simeq \imin i F$.
\item[(ii)] $R\tau_{!!}F \simeq i^!F$.
\end{itemize}
\end{lem}
\dim\ \ (i) The adjunction morphism defines $R\tau_*F \simeq \imin
i \imin \tau R\tau_* F \to \imin i F$. Let $V \in \op^c(Z_{sa})$.
Then
$$\lind {U \supset V} R^k\Gamma(U;F) \simeq
\lind {\substack{U \supset V \\ \tau(U)=V}}R^k\Gamma(U;F) \simeq
R^k\Gamma(\imin \tau(V);F) \simeq R^k\Gamma(V;R\tau_*F),$$ where
$U \in \op(X_{sa})$ and $\RP$-connected. The second isomorphism
follows from Corollary \ref{preceding}.

(ii) The adjunction morphism defines  $i^!F \to
i^!\tau^!R\tau_{!!}F \simeq R\tau_{!!}F$. Let $V \in
\op^c(Z_{sa})$, and let $K$ be a compact subanalytic
$\RP$-connected neighborhood of $V$ in $E$. Then $\imin \tau(V)
\setminus K$ is $\RP$-connected and subanalytic, and $\RP(\imin
\tau(V) \setminus K)=\imin \tau(V) \setminus Z$. By Corollary \ref{preceding}
we have the isomorphism $\mathrm{R}\Gamma(\imin
\tau(V);\mathrm{R}\Gamma_ZF) \simeq \mathrm{R}\Gamma(\imin \tau
(V);\mathrm{R}\Gamma_KF)$.

It follows from the definition of $R\tau_{!!}$ that for any $k \in
\Z$ and $V \in \op^c(Z_{sa})$ we have $R^k\Gamma(V;R\tau_{!!}F)
\simeq \lind K R^k\Gamma(\imin \tau (V);\mathrm{R}\Gamma_KF)$,
where $K$ ranges through the family of compact subanalytic
$\RP$-connected neighborhoods of $V$ in $E$.

On the other hand  for any $k \in \Z$ we have $R^k\Gamma(V;i^!F)
\simeq R^k\Ho(i_*k_V,F) \simeq R^k\Ho(i_*\imin i \imin \tau k_V,F)
\simeq R^k\Gamma(\imin \tau(V),\mathrm{R}\Gamma_ZF)$ and the
result follows. \nopagebreak \newline \qed




\subsection{Fourier-Sato transformation}\label{fouriersato}

Let $E \stackrel{\tau}{\to} Z$ be a real vector bundle, with
dimension $n$ over a real analytic manifold $Z$ and $E^*
\stackrel{\pi}{\to} Z$ its dual. We identify $Z$ as the
zero-section of $E$ and denote $i:Z \hookrightarrow E$ the
embedding, we define similarly $i:Z \hookrightarrow E^*$. We
denote by $p_1$ and $p_2$ the projections from $E \times_Z E^*$:
$$
\xymatrix{& & E \underset{Z}{\times} E^* \ar[lld]_{p_1} \ar[rrd]^{p_2} & &\\
E \ar[rrd]^{\tau} & & && E^* \ar[lld]_{\pi}\\
& & Z && &}
$$
We set
$$P:=\{(x,y) \in E \underset{Z}{\times} E^*;\; \langle x,y\rangle \geq 0\}$$
$$P':=\{(x,y) \in E \underset{Z}{\times} E^*;\; \langle x,y\rangle \leq 0\}$$
and we define the functors
$$\left\{
\begin{array}{l}
\Psi_{P'}=Rp_{1*}\circ \mathrm{R}\Gamma_{P'} \circ p_2^!:D^b_{\RP}(k_{E^*_{sa}}) \to D^b_{\RP}(k_{E_{sa}})\\
\\
\Phi_{P'}=Rp_{2!!}\circ (\cdot)_{P'} \circ p_1^{-1}:D^b_{\RP}(k_{E_{sa}}) \to D^b_{\RP}(k_{E^*_{sa}})\\
\\
\Psi_{P}=Rp_{2*}\circ \mathrm{R}\Gamma_{P} \circ p_1^{-1}:D^b_{\RP}(k_{E_{sa}}) \to D^b_{\RP}(k_{E^*_{sa}})\\
\\
\Phi_{P}=Rp_{1!!}\circ (\cdot)_{P} \circ
p_2^!:D^b_{\RP}(k_{E^*_{sa}}) \to D^b_{\RP}(k_{E_{sa}})
\end{array}
\right.$$

\begin{oss} These functors are well defined, more generally they
send subanalytic sheaves to conic subanalytic sheaves.
\end{oss}





\begin{lem}\label{suppsa} Let $F \in D^b_{\RP}(k_{E_{sa}})$. Then $\supp
  ((\mathrm{R}\Gamma_P(\imin {p_1}F))_{P'})$ is contained in $Z  \times_Z E^*$.
\end{lem}
\dim\ \ We may reduce to the case $F \in \mod_{\RP}(k_{E_{sa}})$.
Then $F=\lind i \rho_*F_i$, with $F_i \in
\mod^{cb}_{\rc,\RP}(k_E)$. We have
\begin{eqnarray*}
H^k(\mathrm{R}\Gamma_P (\imin {p_1}\lind i \rho_*F_i)_{P'}) &
\simeq & \lind i H^k(\mathrm{R}\Gamma_P(\imin {p_1} \rho_* F_i)_{P'}) \\
& \simeq & \lind i \rho_* H^k(\mathrm{R}\Gamma_P(\imin {p_1} F_i)_{P'}) \\
& \simeq & \lind i \rho_* (H^k(\mathrm{R}\Gamma_P (\imin {p_1}
F_i)_{P'}))_{Z  \times_Z E^*},
\end{eqnarray*}
where
the last isomorphism follows from Lemma 3.7.6 of \cite{KS90}.
\nopagebreak
\newline \qed

\begin{lem}\label{ABsa} Let $A$ and $B$ be two closed subanalytic subsets of $E$ such that $A \cup B=E$,
and let $F \in D^b(k_{E_{sa}})$. Then $\mathrm{R}\Gamma_A(F_B)
\simeq (\mathrm{R}\Gamma_AF)_B$.
\end{lem}

\dim\ \ We have a natural arrow $(\Gamma_AF)_B \to \Gamma_A(F_B)$,
and $R(\Gamma_AF)_B \simeq (\mathrm{R}\Gamma_AF)_B$ since
$(\cdot)_B$ is exact. Then we obtain a morphism
$(\mathrm{R}\Gamma_AF)_B \to \mathrm{R}\Gamma_A(F_B)$. It is
enough to prove that for any $k \in \Z$ and for any $F \in
\mod(k_{E_{sa}})$ we have $(R^k\Gamma_AF)_B \iso
R^k\Gamma_A(F_B)$. Since both sides commute with filtrant $\Lind$,
we may assume $F \in \mod^c_{\rc}(k_E)$. Then the result follows
from the corresponding one for classical sheaves. \nopagebreak
\newline \qed

\begin{prop} The two functors $\Phi_{P'},\Psi_P:D^b_{\RP}(k_{E_{sa}}) \to
  D^b_{\RP}(k_{E^*_{sa}})$ are isomorphic.
\end{prop}
\dim\ \ We have the chain of isomorphisms:
\begin{eqnarray*}
\Phi_{P'}  F & = & Rp_{2!!}(\imin {p_1}F)_{P'} \\
 & \simeq & Rp_{2!!}\mathrm{R}\Gamma_P((\imin {p_1}F)_{P'}) \\
 & \simeq & Rp_{2!!}(\mathrm{R}\Gamma_P(\imin {p_1}F))_{P'} \\
 & \simeq & Rp_{2*}(\mathrm{R}\Gamma_P(\imin {p_1}F))_{P'} \\
 & \simeq & Rp_{2*}\mathrm{R}\Gamma_P(\imin {p_1}F).
\end{eqnarray*}
The first isomorphism follows from Lemma \ref{itausa} (ii), the
second one from Lemma \ref{ABsa}, the third one from Lemma
\ref{suppsa} and the last one from Lemma \ref{itausa} (i).
\nopagebreak \newline \qed

\begin{df} Let $F \in D^b_{\RP}(k_{E_{sa}})$.
\begin{itemize}
\item[$(i)$] The Fourier-Sato transform is the functor $$(\cdot)^{\land}:D^b_{\RP}(k_{E_{sa}}) \to D^b_{\RP}(k_{E_{sa}^*})$$
$$F^{\land}=\Phi_{P'}F \simeq \Psi_PF.$$
\item[$(ii)$] The inverse Fourier-Sato transform is the functor $$(\cdot)^{\vee}:D^b_{\RP}(k_{E_{sa}^*}) \to D^b_{\RP}(k_{E_{sa}})$$
$$F^{\vee}=\Psi_{P'}F \simeq \Phi_PF.$$
\end{itemize}
\end{df}

It follows from definition that the functors ${}^\land$ and
${}^\vee$ commute with $R\rho_*$ and $\imin \rho$. We have
quasi-commutative diagrams
$$
\xymatrix{D^b_{\RP}(k_E) \ar[d]^{R\rho_*} \ar@ <2pt> [r]^{\land} & D^b_{\RP}(k_{E^*}) \ar@ <2pt> [l]^{\vee} \ar[d]^{R\rho_*} \\
D^b_{\RP}(k_{E_{sa}}) \ar@ <2pt> [r]^{\land} &
D^b_{\RP}(k_{E_{sa}^*}) \ar@ <2pt> [l]^{\vee}} \ \ \ \
\xymatrix{D^b_{\RP}(k_E) \ar@ <2pt> [r]^{\land} & D^b_{\RP}(k_{E^*}) \ar@ <2pt> [l]^{\vee} \\
D^b_{\RP}(k_{E_{sa}}) \ar@ <2pt> [r]^{\land} \ar[u]^{\imin \rho} &
D^b_{\RP}(k_{E_{sa}^*}) \ar@ <2pt> [l]^{\vee} \ar[u]^{\imin
\rho}.}
$$

This implies that these functors are the extension to conic
subanalytic sheaves of the classical Fourier-Sato and inverse
Fourier-Sato transforms.

\begin{teo}\label{fouriersa} The functors ${}^{\wedge}$ and ${}^{\vee}$ are equivalence of categories, inverse to each others. In particular we have
$$\Ho_{D^b_{\RP}(k_{E_{sa}})}(F,G) \simeq \Ho_{D^b_{\RP}(k_{E_{sa}^*})}(F^\wedge,G^\wedge).$$
\end{teo}
\dim\ \ Let $F \in D^b_{\RP}(k_{E_{sa}})$. The functors
${}^{\wedge}$ and ${}^{\vee}$ are adjoint functors, then we have a
morphism $F \to F^{\land \vee}$. To show that it induces an
isomorphism it is enough to check that $\mathrm{R}\Gamma(U;F) \to
\mathrm{R}\Gamma(U;F^{\land \vee})$ is an isomorphism on a basis
for the topology of $E_{sa}$. Hence we may assume that $U$ is
$\RP$-connected. By Proposition \ref{preceding} we may suppose that
$U$ is an open su\-ba\-na\-ly\-tic cone of $E$. we have the chain
of isomorphisms:
\begin{eqnarray*}
\Rh(k_U,F^{\land \vee}) & = & \Rh(k_U,\Psi_{P'}\Phi_{P'}F)\\
 & \simeq & \Rh(\Phi_{P'}k_U,\Phi_{P'}F)\\
 & \simeq & \Rh(\Phi_{P'}k_U,\Psi_{P}F)\\
 & \simeq & \Rh(\Phi_{P}\Phi_{P'}k_U,F)\\
 & \simeq & \Rh(k_U,F),
\end{eqnarray*}
where the last isomorphism follows from Theorem 3.7.9 of
\cite{KS90} and from the fact that the functors ${}^\land$ and
${}^\vee$ commute with $R\rho_*$. Similarly we can show that for
$G \in D^b_{\RP}(k_{E_{sa}^*})$ we have an isomorphism $G^{\vee
\land} \iso G$. \nopagebreak \newline \qed



\begin{oss} We have seen that the functors ${}^\vee$ and ${}^\wedge$ commute with $\rho_*$ and $\imin \rho$. They do not commute with $\rho_!$ in general since it does not send conic sheaves to conic sheaves. We have the following commutative diagram
$$
\xymatrix{D^b_{\RP}(k_E) \ar[d]^{\rho_{\RP !}} \ar@ <2pt> [r]^{\land} & D^b_{\RP}(k_{E^*}) \ar@ <2pt> [l]^{\vee} \ar[d]^{\rho_{\RP !}} \\
D^b_{\RP}(k_{E_{sa}}) \ar@ <2pt> [r]^{\land} &
D^b_{\RP}(k_{E_{sa}^*}) \ar@ <2pt> [l]^{\vee}}
$$
\end{oss}

Let us study some functorial properties of the Fourier-Sato
transform. Let $Z'$ be another real analytic manifold and let
$f:Z' \to Z$ be a real analytic map. Set $E'=Z' \times_Z E$ and
denote by $f_\tau$ (resp. $f_\pi$) the map from $E'$ to $E$ (resp.
from $E'{}^*$ to $E^*$).

\begin{prop}\label{3.7.13} Let $F \in D^b_{\RP}(k_{E'_{sa}})$. Then:
\begin{eqnarray*}
(Rf_{\tau *}F)^\wedge & \simeq & Rf_{\pi *}(F^\wedge) \\
(Rf_{\tau \RP !!}F)^\wedge & \simeq & Rf_{\pi \RP !!}(F^\wedge).
\end{eqnarray*}
Let $G \in D^b_{\RP}(k_{E_{sa}})$. Then:
\begin{eqnarray*}
(f^!_\tau G)^\wedge & \simeq & f_\pi^!(G^\wedge) \\
(\imin {f_\tau}G)^\wedge & \simeq & \imin {f_\pi}(G^\wedge).
\end{eqnarray*}
\end{prop}
\dim\ \ The result follows adapting Proposition 3.7.13 of
\cite{KS90}.\\ \qed

Let $E_i$, $i=1,2$ be two real vector bundles over $Z$, $f:E_1 \to
E_2$ a morphism of vector bundles. Set
${}^{t}\hspace{-0.8mm}f:E_2^* \to E_1^*$ the dual morphism.

\begin{prop}\label{3.7.14} (i) Let $F \in D^b_{\RP}(k_{E_{1sa}})$. Then:
\begin{eqnarray*}
{}^{t}\hspace{-0.8mm}f^!(F^\vee) & \simeq & (Rf_*F)^\vee \\
{}^{t}\hspace{-0.8mm}f^!(F^\wedge) & \simeq & (Rf_*F)^\wedge
\otimes {}^{t}\hspace{-0.8mm}f^!k_{E_1^*}\\
\imin {{}^{t}\hspace{-0.8mm}f}(F^\wedge) & \simeq & (Rf_{\RP !!}F)^\wedge \\
\imin {{}^{t}\hspace{-0.8mm}f}(F^\vee) & \simeq & (Rf_{\RP !!}F)^\vee
\otimes {}^{t}\hspace{-0.8mm}f^!k_{E_1^*}.
\end{eqnarray*}
%
%
%
%
(ii) Let $G \in D^b_{\RP}(k_{E_{2sa}})$. Then:
\begin{eqnarray*}
(f^!G)^\wedge & \simeq & (R\hspace{0.5mm}{}^{t}\hspace{-0.8mm}f_*G^\wedge) \\
(f^!k_{E_2}\otimes f^!G)^\vee & \simeq &
(R\hspace{0.5mm}{}^{t}\hspace{-0.8mm}f_*G^\vee) \\
(\imin {f}G)^\vee & \simeq & R\hspace{0.5mm}{}^{t}\hspace{-0.8mm}f_{\RP !!}G^{\vee} \\
(f^!k_{E_2}\otimes \imin {f}G)^\wedge & \simeq &
R\hspace{0.5mm}{}^{t}\hspace{-0.8mm}f_{\RP !!}G^{\land} .
\end{eqnarray*}
%
%
\end{prop}
\dim\ \ The result follows adapting Proposition 3.7.14 of
\cite{KS90}.\\ \qed



Let $E_i$, $i=1,2$ be two vector bundles over a real analytic
manifold $Z$. We set for short $E=E_1 \times_Z E_2$ and $E^*=E_1^*
\times_Z E_2^*$. We denote by $\wedge$ the Fourier-Sato transform
on $E_i$, $i=1,2$ and $E$.

\begin{prop}\label{3.7.15}
Let $F_i \in D^b_{\RP}(k_{E_{isa}})$, $i=1,2$. There is an
isomorphism
$$F_1^ \land \underset{Z}{\boxtimes} F_2^\land \simeq \left( F_1
\underset{Z}{\boxtimes} F_2 \right)^\land.$$
\end{prop}
\dim\ \ Let $p_i^j$ and $p_i$ the $i$-th projection defined on
$E_j \times_Z E^*_j$, $j=1,2$ and $E \times E^*$
respectively. Let $P'_j=\{\langle x_j,y_j\rangle\leq 0\} \subset
E_j \times_Z E^*_j$, $j=1,2$ and
$P'=\{\langle(x_1,x_2),(y_1,y_2)\rangle \leq 0\} \subset E
\times_Z E^*$.
The K\"{u}nneth formula gives rise to the isomorphisms:
\begin{eqnarray*}
F_1^ \land \underset{Z}{\boxtimes} F_2^\land & \simeq &
Rp_{2!!}\left(p_1^{1-1}F_1 \underset{Z}{\boxtimes}
p_1^{2-1}F_2\right)_{P'_1 \times_Z P'_2}, \\
\left( F_1 \underset{Z}{\boxtimes} F_2 \right)^\land & \simeq &
Rp_{2!!}\left(p_1^{1-1}F_1 \underset{Z}{\boxtimes}
p_1^{2-1}F_2\right)_{P'}.
\end{eqnarray*}
It is enough to show that for any sheaf $F \in D^+(k_{(E \times_Z
E^*)_{sa}})$ conic with respect to the actions of $\RP$ on $E_j$
and $E^*_j$, $j=1,2$, the morphism $Rp_{2!!}F_{P'} \to
Rp_{2!!}F_{P'_1 \times_Z P'_2}$ induces an isomorphism
$R^kp_{2!!}F_{P'} \iso R^kp_{2!!}F_{P'_1 \times_Z P'_2} $ for any
$k \in \Z$. We may reduce to the case $F$ concentrated in degree
zero. Then as in $\S\,$\ref{consa} one can show that $F=\lind i
\rho_*F_i$, with $F_i$ conic with respect to the actions of $\RP$
on $E_j$ and $E^*_j$, $j=1,2$, $\R$-constructible and with compact
support on the base
for each $i$. We have the chain of isomorphisms
\begin{eqnarray*}
R^kp_{2!!}(\lind i \rho_*F_i)_{P'} & \simeq & \lind i
\rho_*R^kp_{2!}(F_i)_{P'}\\
& \simeq & \lind i \rho_*R^kp_{2!}(F_i)_{P'_1 \times_Z P'_2}\\
& \simeq & R^kp_{2!!}(\lind i \rho_*F_i)_{P'_1 \times_Z P'_2}.
\end{eqnarray*}
$R^kp_{2!!}$ commutes with $\rho_*$ by Lemma \ref{itausa} and the
second isomorphism follows from Proposition 3.7.15 of \cite{KS90}.
\qed

\section{Specialization of subanalytic sheaves}

In this section we define specialization for subanalytic sheaves. We refer
to \cite{KS90} for the classical theory of specialization.

\subsection{Review on normal deformation}\label{rewdef}

Let $X$ be a real $n$-dimensional analytic manifold and let $M$ be
a closed submanifold of codimension $\ell$. As usual we denote by
$T_MX \overset{\tau}{\to} M$ the normal bundle.\\

We follow the notations of \cite{KS90}. We consider the normal
deformation of $X$, i.e. an analytic manifold $\widetilde{X}_M$,
an application $(p,t):\widetilde{X}_M \to X \times \R$, and an
action of $\R \setminus \{0\}$ on $\widetilde{X}_M$
$(\widetilde{x},r) \mapsto \widetilde{x} \cdot r$ satisfying
$$
  \begin{cases}
   \imin p (X \setminus M) & \text{isomorphic to $(X \setminus M) \times (\R \setminus \{0\})$}, \\
   \ \ \imin t (c) & \text{isomorphic to $X$ for each $c \neq 0$}, \\
   \ \ \imin t (0) & \text{isomorphic to $T_MX$}.
  \end{cases}
$$
Let $s:T_MX \hookrightarrow \widetilde{X}_M$ be the inclusion,
$\Omega$ the open subset of $\widetilde{X}_M$ defined by
$\{t>0\}$, $i_\Omega: \Omega \hookrightarrow \widetilde{X}_M$ and
$\widetilde{p}=p \circ i_\Omega$. We get a commutative diagram
$$
\xymatrix{T_MX \ar[r]^s \ar[d]^\tau & \widetilde{X}_M \ar[d]^p &
\Omega \ar[l]_{\ \ i_\Omega} \ar[dl]^{\widetilde{p}} \\
M \ar[r]^{i_M} & X. & }
$$

The morphism $\widetilde{p}$ is smooth and $\Omega$ is
homeomorphic to $X \times \RP$ by the map $(\widetilde{p},t)$.\\

\begin{df} Let $S$ be a subset of $X$. The normal cone to $S$
along $M$, denoted by $C_M(S)$, is the closed conic subset of
$T_MX$ defined by
$$C_M(S)=T_MX \cap \overline{\imin {\widetilde{p}}(S)}.$$
\end{df}


Let us recall the following result of \cite{KS90}.

\begin{lem} \label{4.1.3} Let $V$ be a conic
open subset of $T_MX$.
\begin{itemize}
\item[(i)] Let $W$ be an open neighborhood of $V$ in
$\widetilde{X}_M$ and let $U=\widetilde{p}(W \cap \Omega)$. Then
$V \cap C_M(X \setminus U)=\varnothing$.
\item[(ii)] Let $U$ be an open subset of $X$ such that $V \cap
C_M(X \setminus U)=\varnothing$. Then $\widetilde{p}(U) \cup V$ is
an open neighborhood of $V$ in $\Omega \cup T_MX$.
\end{itemize}
\end{lem}

\begin{lem}\label{4.1.4} Let $V \in
\op(T_MX_{sa})$ be conic. Each
subana\-lytic neighborhood $W$ of $V$ in
$\widetilde{X}_M$ contains a subana\-lytic
neighborhood $\widetilde{W}$ satisfying
\begin{equation} \label{pitilde}
\begin{cases} \text{(i) the fibers of the
map $\widetilde{p}:\widetilde{W} \cap
\Omega \to X$ are connected,} \\
\text{(ii) $\widetilde{p}(\widetilde{W}
\cap \Omega)$ is subanalytic in $X$.}
\end{cases} \end{equation} \end{lem}

\dim\ \ (a) Let $W \in
\op(\widetilde{X}_{Msa})$ be a
neighborhood of $V$. Up to shrink $W$ we
may suppose $V=W \cap T_MX$. Set
$X'=\widetilde{X}_M \setminus (M \times
\R)$, $S=X'/\RP$. Then $\alpha:X' \to S$
is an $\RP$-bundle and $\varphi: S \to X$
is proper. Consider a continuous
subanalytic section of $T_MX \to T_MX/\RP$, extend it to a continuous subanalytic section $\sigma$ of $X' \to S$
and set $$ W'=\bigcup_{x \in
\alpha(W)}W'_x, $$ where $W'_x$
denotes the connected component of
$\imin\alpha(x)\cap W$ containing
$\sigma(x)$. By construction, the fibers
of $W' \to S$ are connected and $W'$ is an
open neighborhood of $V \setminus M$. Let
$$ W''=\bigcup_{x \in M \cap W}W''_x, $$
where $W''_x$ denotes the connected
component of $(\{x\} \times \R) \cap W$
intersecting $M$. Up to shrink $W''$,
$\widetilde{W}=W' \cup W'' \cup (W
\setminus \Omega)$ is an open neighborhood
of $V$ and satisfies (i).

(b) Let us see that $\widetilde{W}$ is
subanalytic and satisfies (ii). We may
reduce to the case $X=\R^n$, $M=\{0\}
\times \R^{n-\ell} \subset X$,
$\widetilde{X}_M=\R^{n+1}$, $T_MX=\R^n
\times \{0\} \subset \widetilde{X}_M$,
$\Omega=\R^n \times \RP \subset
\widetilde{X}_M$. So that $(x',x'',t)
\cdot c \mapsto (cx',x'',\imin c t)$ is
the action of $\RP$ on $\widetilde{X}_M$
and
 $\widetilde{p}:\Omega \to X$ is the map
 $(x',x'',t) \mapsto
(tx',x'')$. In this situation $X'=\R^{n+1}
\setminus (M \times \R) \simeq
\mathbb{S}^{\ell-1} \times \RP \times
\R^{n-\ell} \times \R$. Moreover $X'/\RP
\simeq \mathbb{S}^{\ell-1} \times
\R^{n-\ell} \times \R$, indeed $$ X'=S
\times \RP =\{(ci(\vartheta),x'',s\imin
c),\;(\vartheta,x'',s) \in S,\; c \in \RP\},
$$ where $i:\mathbb{S}^{n-1}
\hookrightarrow \R^n$ denotes the
embedding. The
section $\sigma:S \iso S \times \{1\} \subset X'$ is a
globally subanalytic subset of $X'$. Let
us consider the globally subanalytic (even
semialgebraic) homeomorphism $\psi:\Omega
\to \Omega$ defined by
$\psi(x',x'',t)=(tx',x'',\imin t)$. Then
$\pi \circ \psi=\widetilde{p}$, where
$\pi:\R^n \times \RP \to \R^n$ is the
projection. The set $\psi(W \cap \Omega)$
is still subanalytic.  Let $p \in \R^n$.
Then $\imin \pi(p) \cap \psi(W \cap X'
\cap \Omega)$ is a disjoint union of
intervals. Let us consider the interval
$(m(p),M(p))$, $m(p)<M(p) \in \R \cup \{\pm\infty\}$
intersecting $\psi(\sigma\cap\Omega)$. Then
$\psi(W')=\{(p,r) \in \psi(W \cap X' \cap
\Omega) ;\; m(p)<r<M(p)\}$.
The set $\psi(W'\cap\Omega)$ is open subanalytic (it is a consequence of Proposition 1.2, Chapter 6 of \cite{VD98}).
Moreover, up to shrink $W$ we may suppose
that $\widetilde{p}(W' \cap
\Omega)=\widetilde{p}(\sigma \cap W \cap \Omega)$
which is subanalytic. Indeed, since we are working in a local chart, we may assume that $W \cap \sigma$ is globally subanalytic.  Let $x'' \in
\R^{n-\ell}$. Then $(\{x''\} \times \R)
\cap W'' $ is a disjoint union of (bounded, up to shrink $W$)
intervals. Let us consider the interval
$(m(x''),M(x''))$, $m(x'')<M(x'') \in \R$
containing $0$. Then $W''=\{r \in W \cap
(M \times \R) ;\; m(x'')<r<M(x'')\}$.
The set $W''$ is subanalytic (it is a consequence of Proposition 1.2, Chapter 6 of \cite{VD98}).
Moreover (up to shrink $W''$) $\widetilde{p}(W'' \cap
\Omega)=W'' \cap M$, which is subanalytic. \nopagebreak \newline \qed

\subsection{Specialization of subanalytic sheaves}\label{nu}

\begin{df} The specialization along $M$ is the functor
\begin{eqnarray*}
\nu^{sa}_M:D^b(k_{X_{sa}}) & \to & D^b(k_{T_MX_{sa}}) \\
F & \mapsto & \imin s\mathrm{R}\Gamma_\Omega\imin p F.
\end{eqnarray*}
\end{df}

\begin{teo} \label{4.2.3} Let $F \in D^b(k_{X_{sa}})$.
\begin{itemize}
\item[(i)] $\nu^{sa}_MF \in D^b_{\RP}(k_{T_MX_{sa}})$.
\item[(ii)] Let $V$ be a conic subanalytic open subset of $T_MX$. Then:
$$H^j(V;\nu^{sa}_MF) \simeq \lind U H^j(U;F),$$
where $U$ ranges through the family of $\op(X_{sa})$ such that
$C_M(X \setminus U) \cap V=\varnothing$.
\item[(iii)] One has the isomorphisms
\begin{eqnarray*}
& (\nu^{sa}_MF)|_M \simeq R\tau_*(\nu^{sa}_MF) \simeq F|_M, & \\
& (\mathrm{R}\Gamma_M \nu^{sa}_MF)|_M \simeq R\tau_{!!}\nu^{sa}_MF \simeq
(\mathrm{R}\Gamma_MF)|_M.
\end{eqnarray*}
\end{itemize}
\end{teo}
\dim\ \ (i) We may reduce to the case $F \in \mod(k_{X_{sa}})$.
Hence $F=\lind i \rho_* F_i$ with $F_i \in \mod_{\rc}(k_{X_{sa}})$
for each $i$.
We have $\imin p \lind i \rho_* F_i \simeq \lind i \rho_* \imin p
F_i$ and $\imin p F_i$ is $\R$-constructible and conic for each
$i$. Hence $\imin p F$ is conic. Since the functors
$\mathrm{R}\Gamma_\Omega$ and $\imin s$ send conic sheaves to conic sheaves
we obtain $\imin s \mathrm{R}\Gamma_\Omega \imin pF=\nu^{sa}_MF \in
D^b_{\RP}(k_{T_MX_{sa}})$.

(ii) Let $U \in \op(X_{sa})$ such that $V \cap C_M(X \setminus
U)=\varnothing$. We have the chain of morphisms
\begin{eqnarray*}
\mathrm{R}\Gamma(U;F) & \to & \mathrm{R}\Gamma(\imin p(U);\imin p F) \\
& \to & \mathrm{R}\Gamma(\imin p(U) \cap \Omega;\imin p F) \\
& \to & \mathrm{R}\Gamma(\imin {\widetilde{p}}(U) \cup V;\mathrm{R}\Gamma_\Omega\imin p F) \\
& \to & \mathrm{R}\Gamma(V;\nu^{sa}_MF)
\end{eqnarray*}
where the third arrow exists since $\imin {\widetilde{p}}(U) \cup
V$ is a neighborhood of $V$ in $\overline{\Omega}$ by Lemma
\ref{4.1.3} (ii). Let us show that it is an isomorphism.
Let $V$ be a conic open subanalytic subset of $T_MX$. We have
\begin{eqnarray*}
H^k(V;\nu^{sa}_MF) & \simeq & \lind W H^k(W;\mathrm{R}\Gamma_\Omega\imin p F) \\
& \simeq & \lind W H^k(W \cap \Omega;\imin p F),
\end{eqnarray*}
where $W$ ranges through the family of subanalytic open
neighborhoods of $V$ in $X$. By Lemma \ref{4.1.4} we may assume
that $W$ satisfies \eqref{pitilde}.
Since $\imin p F$ is conic, we have
\begin{eqnarray*}
H^k(W \cap \Omega;\imin p F) & \simeq & H^k(\imin p(p(W \cap
\Omega));\imin p F) \\
& \simeq & H^k(p(W \cap \Omega) \times \{1\};\imin p F) \\
& \simeq & H^k(p(W \cap \Omega);F),
\end{eqnarray*}
where the second isomorphism follows since every subanalytic
neighborhood of $p(W \cap \Omega) \times \{1\}$ contains an
$\RP$-connected subanalytic neighborhood (the proof is similar to
that of Lemma \ref{4.1.4}). By Lemma \ref{4.1.3} (i) we have that
$p(W \cap \Omega)$ ranges through the family of subanalytic open
subsets $U$ of $X$ such that $V \cap C_M(X \setminus
U)=\varnothing$ and we obtain the result.

(iii)  The result follows adapting Theorem 4.2.3 (iv) of
\cite{KS90}. \nopagebreak \newline \qed


\begin{prop}\label{nuinj} Let $F \in \mod(k_{X_{sa}})$ be quasi-injective. Then $\nu^{sa}_MF$ is concentrated in degree zero. Moreover $\nu_M^{sa}F$ is $\RP$-quasi-injective.
\end{prop}
\dim\ \ Since $\nu^{sa}_MF$ is conic, it is enough to prove that $H^j(V;\nu^{sa}_MF)=0$, $j \neq 0$, when $V$ is a conic subanalytic open subset of $T_MX$. By Theorem \ref{4.2.3} we have $H^j(V;\nu^{sa}_MF) \simeq \lind U H^j(U;F)$, where $U$ ranges through the family of $\op(X_{sa})$ such that
$C_M(X \setminus U) \cap V=\varnothing$, and $H^j(U;F)=0$ if $j \neq 0$ since $F$ is quasi-injective. Moreover if $V_1 \subset V_2$ are conic subanalytic open subset of $T_MX$ the morphism
$$\Gamma(V_2;\nu_M^{sa}F) \simeq \lind {U_1} \Gamma(U_2;F)\to \lind {U_2} \Gamma(U_1;F) \simeq \Gamma(V_1;\nu_M^{sa}F),$$
where $U_i$ ranges through the family of $\op(X_{sa})$ such that
$C_M(X \setminus U_i) \cap V_i=\varnothing$, $i=1,2$,
is surjective.
\nopagebreak \newline \qed

 Let us study the relation with the classical functor of
specialization  $\nu_M$.

\begin{prop} Let $F \in D^b(k_X)$. Then $\imin \rho \nu^{sa}_M
R\rho_* F \simeq \nu_MF$.
\end{prop}
\dim\ \
We have to show that for each $x \in T_MX$ we have $H^k(\imin
\rho\nu^{sa}_MR\rho_*F)_x \simeq H^k(\nu_MF)_x$. Hence we have to
prove the isomorphism $\lind {x \in V} H^k(V;\nu^{sa}_MR\rho_*F)
\simeq \lind {x \in V} H^k(V;\nu_MF)$, where $V$ ranges through
the family of open $\RP$-connected relatively compact subanalytic
subset of $T_MX$. This is a consequence of Theorem \ref{4.2.3} and
Theorem 4.2.3 of \cite{KS90}. \\ \qed

\begin{oss} Remark that the functor of specialization does not commute with
$R\rho_*$ and $\imin \rho$ in general.

 In fact let $V \in
\op^c_{sa}(T_MX)$ be $\RP$-connected and let $F \in \mod(k_X)$.
Then $H^k(V;\nu^{sa}_MR\rho_*F) \simeq \lind U H^k(U;F)$, where
$U$ ranges through the family of $\op(X_{sa})$ such that $C_M(X
\setminus U) \cap V=\varnothing$, which is not cofinal to the
family of $\op(X)$ such that $C_M(X \setminus U) \cap
V=\varnothing$.

Now let $V \in \op^c_{sa}(T_MX)$ be $\RP$-connected and let $G \in
\mod(k_{X_{sa}})$. Then $H^k(V;\nu_M\imin \rho G) \simeq \lind U
\lpro {W \subset\subset U} H^k(W;G)$, where $U$ ranges through the
family of $\op(X_{sa})$ such that $C_M(X \setminus U) \cap
V=\varnothing$ and $W \in \op(X_{sa}).$ Then $H^k(\imin \rho\nu^{sa}_MG)_x \neq H^k(\nu_M\imin \rho G)_x$.\\

\end{oss}


Let $f:X \to Y$ be a morphism of manifolds, let $N$ be a closed
submanifold of $Y$ of codimension $k$ and assume $f(M)\subset N$.
We denote by $f'$ the map from $TX$ to $X \times_Y TY$ associated
to $f$ and by $f_\tau:X \times_Y TY \to TY$ the base change. We
denote by $Tf$ the composite map. Similarly, replacing $X,Y,TX,TY$
by $M,N,T_MX,T_NY$ we get the morphisms $f'_M,f_{M\tau}, T_Mf$.\\

We have a commuative diagram, where all the squares are cartesian

$$
\xymatrix{T_MX \ar[d]^{T_Mf} \ar[r]^{s_X} & \widetilde{X}_M
\ar[d]^{\widetilde{f}'} & \Omega_X \ar[l]^{i_{\Omega_X}}
\ar[r]^{\widetilde{p}_X} \ar[d]^{\widetilde{f}} & X \ar[d]^f \\
T_NY  \ar[r]^{s_Y} & \widetilde{Y}_N & \Omega_Y
\ar[l]^{i_{\Omega_Y}} \ar[r]^{\widetilde{p}_Y} & Y.}
$$

Recall that the following diagram is not cartesian in general

\begin{equation}\label{cleand}
\xymatrix{\widetilde{X}_M \ar[d]^{\widetilde{f}'} \ar[r]^{p_X} & X
\ar[d]^f \\ \widetilde{Y}_N \ar[r]^{p_Y} & Y.}
\end{equation}

\begin{df} (i) One says that $f$ is clean with respect to $N$ if
$\imin f(N)$ is a submanifold $M$ of $X$ and the map
${}^{t}\hspace{-0.8mm}f_M':M \times_N T^*_NY \to T^*_MX$ is
surjective.

(ii) One says that $f$ is transversal to $N$ if the map
${}^{t}\hspace{-0.8mm}f'|_{X \times_YT^*_NY}:X \times_YT^*_NY \to
T^*Y$ is injective.
\end{df}

If $f$ is transversal to $N$ and $\imin f(N)=M$, then the square
\eqref{cleand} is cartesian.\\

We will not prove the following results, which can be easily
recovered adapting $\S$\,IV.4.2 of \cite{KS90}, using the
construction we did for subanalytic sheaves.

\begin{prop} Let $F \in D^b(k_{X_{sa}})$.

\begin{itemize}
\item[(i)] There exists a commutative diagram of canonical
morphisms
$$
\xymatrix{R(T_Mf)_{\RP !!}\nu_MF  \ar[d] \ar[r] & \nu_NRf_{\RP !!}F \ar[d] \\
R(T_Mf)_*\nu_MF & \nu_NRf_*F. \ar[l]}
$$
\item[(ii)] Moreover if $\supp F \to Y$ and $C_M(\supp F) \to
T_NY$ are proper, and if $\supp F \cap \imin f(N) \subset M$, then
the above morphisms are isomorphisms.

In particular $f$ is clean with respect to $N$ and $\imin f(N)=M$,
then the above morphisms are isomorphisms.
\end{itemize}
\end{prop}

\begin{prop}\label{4.2.5} Let $G \in D^b(k_{Y_{sa}})$.

\begin{itemize}

\item[(i)] There exists a commutative diagram of canonical
morphisms
$$
\xymatrix{\omega_{T_MX/T_NY} \otimes \imin{(T_Mf)}\nu_NG  \ar[d] \ar[r] & \nu_M(\omega_{X/Y} \otimes \imin f G) \ar[d] \\
(T_Mf)^!\nu_N G & \nu_Mf^!G \ar[l].}
$$
\item[(ii)] Moreover if $f:X \to Y$ and $f|_M:M \to N$ are smooth
the above morphisms are isomorphisms.
\end{itemize}
\end{prop}

Let $X$ and $Y$ be two real analytic manifolds and let $M,N$ be
two closed submanifolds of $X$ and $Y$ respectively.

\begin{prop} Let $F \in D^b(k_{X_{sa}})$ and $G \in
D^b(k_{Y_{sa}})$. There is a natural morphism
$$\nu_MF \boxtimes \nu_NG \to \nu_{M \times N}(F \boxtimes G).$$
\end{prop}

\begin{cor} Let $F,G \in D^b(k_{X_{sa}})$. There is a natural morphism
$$\nu_MF \otimes \nu_MG \to \nu_M(F \otimes G).$$
\end{cor}

\section{Microlocalization of subanalytic sheaves}

With the construction of the Fourier-Sato transform and the specialization we have all
the tools to define the functor of
microlocalization in the framework of subanalytic sites. See \cite{KS90} for the classical theory of microlocalization. Then we introduce the functor
$\muh^{sa}$ for subanalytic sheaves, we study the relations with the notion of microsupport of \cite{KS03} and with the functor of ind-microlocalization of \cite{KSIW}.

\subsection{Microlocalization of subanalytic sheaves}\label{mu}

 Let us denote by $T_M^*X$ the
conormal bundle to $M$ in $X$, i.e. the kernel of the map $M
\times_X T^*X \to T^*M$. We denote by $\pi$ the projection $T_M^*X
\to M$.

\begin{df} Let $F \in D^b(k_{X_{sa}})$. The microlocalization of
$F$ along $M$ is the Fourier-Sato transform of the specialization,
i.e.
$$\mu^{sa}_M F = (\nu^{sa}_M F)^\wedge.$$
\end{df}

\begin{teo} \label{4.3.2} Let $F \in D^b(k_{X_{sa}})$.
\begin{itemize}
\item[(i)] $\mu^{sa}_MF \in D^b_{\RP}(k_{T^*_MX_{sa}})$.
\item[(ii)] Let $V$ be an open convex subanalytic cone of
$T_M^*X$. Then:
$$H^j(V;\mu^{sa}_MF) \simeq \lind {U,Z} H^j_Z(U;F),$$
where $U$ ranges through the family of $\op(X_{sa})$ such that $U
\cap M=\pi(V)$ and $Z$ through the family of closed subanalytic
subsets such that $C_M(Z) \subset V^\circ$, where $V^\circ$
denotes the polar cone.
\item[(iii)] One has the isomorphisms
\begin{eqnarray*}
& (\mu^{sa}_MF)|_M \simeq R\pi_*(\mu^{sa}_MF) \simeq i_M^!F, & \\
& (\mathrm{R}\Gamma_M \mu^{sa}_MF)|_M \simeq R\pi_{!!}\mu^{sa}_MF \simeq
\imin {i_M} F \otimes i_M^!k_X.
\end{eqnarray*}
\end{itemize}
\end{teo}
\dim\ \ The result follows from the functorial properties of the Fourier-Sato transform and Theorem \ref{4.2.3}.\\
\qed

As in classical sheaf theory, we get the Sato's triangle for
subanalytic sheaves:
$$
F|_M \otimes \omega_{M|X} \to {\rm R}\Gamma_MF|_M \to
R\dot{\pi}_*\mu^{sa}_MF \stackrel{+}{\to}
$$
where $\dot{\pi}$ is the restriction of $\pi$ to $T^*_MX \setminus
M$. \\

\begin{prop}\label{muinj} Let $F \in \mod(k_{X_{sa}})$ be quasi-injective. Then $\imin \rho \mu^{sa}_MF$ is concentrated in degree zero.
\end{prop}
\dim\ \ The result follows from Theorem \ref{4.3.2} (ii).
\nopagebreak \newline \qed

\begin{oss} Remark that the functor of microlocalization does not commute with
$R\rho_*$ and $\imin \rho$ since specialization does not. If $F
\in D^b(k_X)$ we have $\imin \rho \mu^{sa}_M R\rho_*F \simeq
\mu_MF$ since the Fourier-Sato transform commutes with $\imin
\rho$ and $\imin \rho \circ \nu^{sa}_M \circ R\rho_* \simeq
\nu_M$.
\end{oss}


Let $f:X \to Y$ be a morphism of manifolds, let $N$ be a closed
submanifold of $Y$ of codimension $k$ and assume $f(M)\subset N$.
The map $Tf$ defines the maps

$$
\xymatrix{T^*X & X \times_Y T^*Y \ar[l]_{{}^{t}\hspace{-0.8mm}f'\ \
} \ar[r]^{\ \ \ \ f_\pi} & T^*Y}
$$
and similarly one can define the maps ${}^{t}\hspace{-0.8mm}f'_M$
and $f_{M\pi}$.

Applying the Fourier-Sato transform to the morphisms of
$\S$\,\ref{nu} we get the following results (see also \cite{KS90}
$\S$\,IV.4.3 for the classical case)

\begin{prop} Let $F \in D^b(k_{X_{sa}})$.

\begin{itemize}
\item[(i)] There exists a commutative diagram of canonical
morphisms
$$
\xymatrix{Rf_{M\pi \RP !!}{}^{t}\hspace{-0.8mm}f'{}_{\hspace{-1.5mm}M}^{-1}\mu_MF  \ar[d] \ar[r] & \mu_NRf_{!!}F \ar[d] \\
Rf_{M\pi *}({}^{t}\hspace{-0.8mm}f'{}_{\hspace{-1.5mm}M}^! \otimes
\omega_{X/Y} \otimes \omega^{\otimes-1}_{M/N}) & \mu_NRf_*F.
\ar[l]}
$$
\item[(ii)] Moreover if $\supp F \to Y$ and $C_M(\supp F) \to
T_NY$ are proper, and if $\supp F \cap \imin f(N) \subset M$, then
the above morphisms are isomorphisms.

In particular $f$ is clean with respect to $N$ and $\imin f(N)=M$,
then the above morphisms are isomorphisms.
\end{itemize}
\end{prop}

\begin{prop} Let $G \in D^b(k_{Y_{sa}})$.

\begin{itemize}

\item[(i)] There exists a commutative diagram of canonical
morphisms
$$
\xymatrix{R\hspace{-0.8mm}f'_{M\RP !!}(\omega_{M/N} \otimes \imin{f_{M\pi}}\mu_NG)  \ar[d] \ar[r] & \mu_M(\omega_{X/Y} \otimes \imin f G) \ar[d] \\
R\hspace{-0.8mm}f'_{M*}f_{M\pi}^!\mu_N G & \mu_Mf^!G \ar[l].}
$$
\item[(ii)] Moreover if $f:X \to Y$ and $f|_M:M \to N$ are smooth
the above morphisms are isomorphisms.
\end{itemize}
\end{prop}

Let $X$ and $Y$ be two real analytic manifolds and let $M,N$ be
two closed submanifolds of $X$ and $Y$ respectively.

\begin{prop} Let $F \in D^b(k_{X_{sa}})$ and $G \in
D^b(k_{Y_{sa}})$. There is a natural morphism
$$\mu_MF \boxtimes \mu_NG \to \mu_{M \times N}(F \boxtimes G).$$
\end{prop}

\begin{cor} Let $M$ be a closed submanifold of $X$ and let
$\gamma:T^*_MX \times_M T^*_MX \to T^*_MX$ be the morphism given
by the addition. There is a natural morphism
$$\mathrm{R}\gamma_{\RP !!}\left(\mu_MF \underset{M}{\boxtimes} \mu_MG\right) \to \mu_M\left(F \boxtimes G\right) \otimes \omega_{M/X}.$$
\end{cor}

\subsection{The functor $\muh^{sa}$}

 We denote by $\Delta$
the diagonal of $X \times X$, and we denote by $\delta$ the
diagonal embedding. The normal deformation of the diagonal in $X
\times X$ can be visualized by the following diagram
\begin{equation}\label{normdef}
\xymatrix{TX \ar[r]^{\hspace{-0.8cm}\sim} & T_{\Delta}(X \times X)
\ar[r]^{\hspace{0.3cm}s} \ar[d]^{\tau_X} & \widetilde{X \times X}
\ar[d]^p & \Omega
\ar[l]_{\hspace{0.5cm}i_\Omega} \ar[dl]^{\widetilde{p}} \\
& \Delta \ar[r]^{\delta} \ar[dr]^{\sim} & X \times X \ar@ <2pt>
[d]^{q_2} \ar@ <-2pt> [d]_{q_1} & \\
& & X. & }
\end{equation}
Set $p_i=q_i \circ p$, $i=1,2$. While $\widetilde{p}$ and $p_i$,
$i=1,2$, are smooth, $p$ is not, and moreover the square is not
cartesian.\\



\begin{df} Let $F \in D^b_{\rc}(k_X)$ and $G \in D^b(k_{X_{sa}})$. We set
$$\muh^{sa}(F,G):=\mu^{sa}_\Delta \rh(\imin q_2 F, q_1^!G) = (\nu^{sa}_\Delta\rh(\imin
{q_2}F,q_1^!G))^\wedge.$$
\end{df}

Let $\pi$ denote the projection from $T^*_\Delta(X \times X)$ to
$\Delta \simeq X$.

\begin{prop}\label{4.4.2} Let $F \in D^b_{\rc}(k_X)$ and $G \in D^b(k_{X_{sa}})$. There is a canonical isomorphism
$\pi_*\muh^{sa}(F,G) \simeq \rh(F,G)$.
\end{prop}
\dim\ \ The result follows adapting Proposition 4.4.2 of
\cite{KS90}. \\ \qed

\begin{oss} The functor $\muh^{sa}$ is well defined also if $F \in D^b(k_{X_{sa}})$. In this case we do not know if $\muh^{sa}(F,G)$ has bounded cohomology or not.
\end{oss}

\begin{oss} Adapting the results of $\S$\,\ref{nu} and $\S$\,\ref{mu} one gets
the functorial properties of $\muh$ for subanalytic sheaves. Since
the proofs are essentially the same as the classical ones we will
skip them and refer to \cite{KS90}.
\end{oss}

Let $\pi:T^*X \to  X$ be the projection and consider the canonical
1-form $\omega$, the restriction to the diagonal of the map
${}^{t}\hspace{-0.4mm}\pi':T^*X \times_X T^*X \to T^*T^*X$. We
have a diagram
$$
\xymatrix{T^*T^*X & T^*X \times_X T^*X
\ar[l]_{{}^{t}\hspace{-0.4mm}\pi'\ \ } \ar[r]^{\ \ \ \ \ \ \pi_\pi} & T^*X \\
& \Delta_{T^*X} \ar[lu]^{\omega} \ar[u]^{\delta_{T^*X}}
\ar[ru]_\sim}
$$

\begin{lem}\label{omegamu} Let $F \in D^b_{\rc}(k_X)$ and $G \in D^b(k_{X_{sa}})$ We have
$$\imin {\omega} \muh^{sa}(\imin \pi F,\imin \pi G) \simeq
\muh^{sa}(F,G).$$
\end{lem}
\dim\ \ We have the isomorphism
${}^{t}\hspace{-0.4mm}\pi'_{!!}\imin{\pi_\pi} \muh^{sa}(F,G) \iso
\muh^{sa}(\imin \pi F,\imin \pi G)$. Hence we get the isomorphisms
\begin{eqnarray*}
\imin {\omega}\muh^{sa}(\imin \pi F,\imin \pi G) & \simeq & \imin
{\omega}{}^{t}\hspace{-0.4mm}\pi'_{!!}\imin{\pi_\pi}
\muh^{sa}(F,G)
\\
& \simeq & \imin {\delta_{T^*X}}\imin {\pi_\pi}\muh^{sa}(F,G) \\
& \simeq & \muh^{sa}(F,G).
\end{eqnarray*}
\qed

\subsection{Microlocalization and microsupport}

In \cite{KS90} the authors prove that the support of $\muh(F,G)$ is contained in the product of the microsupports of $F$ and $G$. We extend this result to the functor $\muh^{sa}$. Let $X$ be a real analytic manifold and let $T^*X \stackrel{\pi}{\to} X$ be the cotangent bundle. We recall the following two equivalent definitions of microsupport of a subanalytic sheaves of \cite{KS03}. For the notion of microsupport for classical sheaves we refer to \cite{KS90}. For the functorial properties of the microsupport of subanalytic sheaves we refer to \cite{Ma08}.

\begin{df} The microsupport of $F \in D^b(k_{X_{sa}})$, denoted by $SS(F)$ is the subset of $T^*X$ defined as follows. Let $p \in T^*X$, then $p \notin SS(F)$ if one of the following equivalent conditions is satisfied.
\begin{itemize}
\item[(i)] There exists a conic neighborhood $U$ of $p$ and a small filtrant system $\{F_i\}$ in $C^{[a,b]}(\mod_{\rc}(k_X))$ with  $SS(F_i) \cap U=\varnothing$ such that $F$ is quasi-isomorphic to $\lind i \rho_* F_i$ in a neighborhood of $\pi(p)$.
\item[(ii)] There exists a conic neighborhood $U$ of $p$ such that for any $G \in D_{\rc}^b(k_X)$ with $\supp(G)\subset\subset \pi(U)$ and such that $SS(G) \subset U \cup T^*_XX$, one has $\Ho_{D^b(k_{X_{sa}})}(G,F)=0$.
\end{itemize}
\end{df}

\begin{oss} In \cite{KS03} microsupport was defined for ind-sheaves. The above definition follows from the equivalence between subanalytic sheaves and ind-$\R$-constructible sheaves (see \cite{Ma08} for details).
\end{oss}

 We need the following result of \cite{Ma08}.

\begin{lem} Let $X,Y$ be two real analytic manifolds and let $q_1,q_2$ be the projections from $X \times Y$ to $X$ and $Y$ respectively. Let $G \in D^b_{\rc}(k_Y)$ and $F \in D^b(k_{X_{sa}})$. Then
\begin{equation} \label{anarita}
SS(\rh(\imin q_1 G,q_2^!F)) \subseteq  SS(F) \times SS(G)^a.
\end{equation}
\end{lem}

Let $M$ be a real closed submanifold of $X$.

\begin{prop}\label{muSS(F)} Let $F \in D^b(k_{X_{sa}})$. Then $\supp (\mu_M^{sa} F) \subseteq SS(F) \cap T^*_MX$.
\end{prop}
\dim\ \ Let $F \in D^b(k_{X_{sa}})$ and let $p \notin SS(F)$. There exists conic neighborhood $U$ of $p$ and a small filtrant system $\{F_i\}$ in $C^{[a,b]}(\mod_{\rc}(k_X))$ with  $SS(F_i) \cap \overline{U}=\varnothing$ such that there exists $W \in \op(X_{sa})$ with $U \subseteq \imin \pi (W)$ and $F_W \simeq \lind i \rho_* F_i$.
We have $H^k\mu_M^{sa} F_W \simeq \lind i \rho_* H^k \mu_M F_{iW}$, hence $(\mu_M^{sa} F)|_U=0$ since $\supp (\mu_M F_i) \subseteq SS(F_i)$. \\ \qed

\begin{cor}\label{muhSS(F)G} Let $G \in D^b_{\rc}(k_X)$, $F \in D^b(k_{X_{sa}})$. Then
$$
\supp (\muh^{sa}(F,G)) \subseteq SS(F) \cap SS(G).
$$
\end{cor}
The result follows from Proposition \ref{muSS(F)} and \eqref{anarita}.\\
\qed

Let $f:X \to Y$ be a morphism of real analytic manifolds and denote by $f_\pi: X \times_Y T^*Y \to T^*Y$ the base change map.

\begin{df} Let $f:X \to Y$ be a morphism of real analytic manifolds and let $F \in D^b(k_{Y_{sa}})$. One says that $f$ is non characteristic for $SS(F)$ if
$$
\imin f_\pi(SS(F)) \cap T^*_XY \subseteq X \times_Y T^*_YY.
$$
If $f$ is a closed embedding $X$ is said to be non characteristic.
\end{df}

\begin{prop}\label{noncharSS(F)} Let $f:X \to Y$ be a morphism of real analytic manifolds and let $F \in D^b(k_{Y_{sa}})$. Assume that $f$ is non characteristic for $SS(F)$. Then the natural morphism
$$
\imin f F \otimes \omega_{X|Y} \to f^!F
$$
is an isomorphism.
\end{prop}
\dim\ \ We may reduce to the case $f$ closed embedding, hence we have to prove the isomorphism $F|_X \otimes \omega_{X|Y} \simeq \mathrm{R}\Gamma_XF|_X$ when $SS(F) \cap T^*_XY \subseteq T^*_YY$. Consider the Sato's triangle
$$
\dt{F|_X \otimes \omega_{X|Y}}{\mathrm{R}\Gamma_XF|_X}{R\dot{\pi}_*\mu^{sa}_XF}.
$$
Since $SS(F) \cap T^*_XY \subseteq T^*_YY$ we have $R\dot{\pi}_*\mu^{sa}_XF=0$ by Proposition \ref{muSS(F)} and the result follows. \\ \qed

\begin{lem}\label{rhD'} Let $F \in D^b_{\rc}(k_X)$ and let $G \in
D^b(k_{X_{sa}})$. Then
$$
D'F \boxtimes  G \iso \rh(\imin {q_1}F,\imin {q_2}G).
$$
\end{lem}
\dim\ \ We may reduce to the case $F=k_U$ with $D'k_U \simeq k_{\overline{U}}$ and $G$ quasi-injective. Set $G=\lind i \rho_* G_i$ with $G_i \in \mod_{\rc}(k_X)$, we have the chain of isomorphisms
$$
(\imin {q_2}G)_{\imin {q_1}(\overline{U})} \simeq \lind i \rho_*(\imin {q_2}G_i)_{\imin {q_1}(\overline{U})} \simeq \lind i \rho_*\Gamma_{\imin {q_1}(U)}(\imin {q_2}G_i) \simeq \Gamma_{\imin {q_1}(U)}(\imin {q_2}G)
$$
where the second isomorphism follows from Proposition 3.4.4 of \cite{KS90}.\\
\qed

\begin{prop}\label{D'rhrcSS} Let $F \in D^b_{\rc}(k_X)$ and let $G \in
D^b(k_{X_{sa}})$. Suppose that $SS(F) \cap SS(G) \subseteq T_X^*X$. Then
$$
D'F \otimes  G \iso \rh(F,G).
$$
\end{prop}
\dim\ \ Let $\delta: \Delta \to X \times X$ be the embedding and
let us consider the Sato's triangle
\begin{eqnarray*}
& \imin \delta\rh(\imin {q_1}F,q_2^!G) \otimes
\omega_{\Delta|X
\times X}  \to  \delta^!\rh(\imin {q_1}F,q_2^!G) & \notag \\
 & \to R\dot{\pi}_*\muh^{sa}(F,G) \stackrel{+}{\to}. &
\end{eqnarray*}
We have $\delta^!\rh(\imin {q_1}F,q_2^!G) \simeq
\rh(F,G)$. Moreover
\begin{eqnarray*}
\imin \delta\rh(\imin {q_1}F,q_2^!G) \otimes
\omega_{\Delta|X
\times X} & \simeq & \imin \delta\rh(\imin {q_1}F,\imin {q_2}G) \\
& \simeq & \imin \delta (D'F \boxtimes G) \\
& \simeq & D'F \otimes G
\end{eqnarray*}
where the second isomorphism follows from Lemma \ref{rhD'}. Then we obtain a distinguished triangle
$$
\dt{D'F \otimes G}{\rh(F,G)}{R\dot{\pi}_*\muh^{sa}(F,G)}
$$
and the result follows since $R\dot{\pi}_*\muh^{sa}(F,G)=0$ by Corollary \ref{muhSS(F)G}.\\ \qed

\subsection{The link with the functor $\mu$ of microlocalization}

We will study the relation between microlocalization for
subanalytic sheaves and the functor $\mu$ of \cite{KSIW}. We will
need first some results. Let $X$ be a real analytic manifold and
consider the normal deformation of $\Delta$ in $X \times X$
visualized by the diagram \eqref{normdef}.
\begin{lem}\label{omegap} Let $G \in D^b(k_{(X \times X)_{sa}})$ and let $F \in D^b_{\rc}(k_{(X \times X)_{sa}})$.
 We have $\mathrm{R}\Gamma_\Omega\imin p \rh(F,G) \simeq \mathrm{R}\Gamma_\Omega \rh(\imin p F,\imin p G)$.
\end{lem}
\dim\ \ We may reduce to the case $F \in \mod_{\rc}(k_{X \times X})$ and $G \in \mod(k_{(X \times
X)_{sa}})$. Hence $G=\lind i \rho_* G_i$ with $G_i \in
\mod_{\rc}(k_{(X \times X)_{sa}})$ for each $i$. Let
$\Omega_{\widetilde{X \times X}_{sa}}$ be the site induced by
$\widetilde{X \times X}_{sa}$ on $\Omega$. We have
\begin{eqnarray*}
H^k \imin {i_{\Omega_{\widetilde{X \times X}_{sa}}}} \imin p
\rh(F,\lind i \rho_* G_i)
 & \simeq & \lind i \rho_* H^k \imin
{i_\Omega} \imin p \rh(F, G_i) \\
 & \simeq & \lind i \rho_* H^k  \rh(\imin
{i_\Omega} \imin p F,\imin
{i_\Omega} \imin p G_i) \\
& \simeq & \lind i \rho_* H^k \imin
{i_\Omega}  \rh(\imin p F,\imin p G_i) \\
& \simeq & H^k \imin {i_{\Omega_{\widetilde{X \times X}_{sa}}}}
\rh(\imin p F,\lind i \rho_* \imin p G_i).
\end{eqnarray*}
The second isomorphism follows because $p \circ i_\Omega$ is
smooth. Composing with $Ri_{\Omega_{\widetilde{X \times
X}_{sa}}\hspace{-2mm}*}$ we obtain the result. \nopagebreak
\newline \qed

\begin{lem}\label{omegabar} Let $G \in D^b(k_{X_{sa}})$, then for $i=1,2$
\begin{eqnarray}
k_{\overline{\Omega} }\otimes \imin {p_i}G & \simeq & \rh(k_\Omega,\imin
{p_i}G) \label{i} \\
k_{\Omega} \otimes \imin {p_i}G & \simeq & \rh(k_{\overline{\Omega}},\imin
{p_i}G). \label{ii}
\end{eqnarray}
\end{lem}
\dim\ \ Let us prove \eqref{i}. Since for $i=1,2$ $p_i$ is smooth, Proposition 3.16 of \cite{Ma08} implies that $SS(\imin {p_i}G) \cap SS(k_\Omega)$ is
contained on the zero section of $T^*(\widetilde{X \times X})$. Then the result is a consequence of Proposition \ref{D'rhrcSS}, and the fact that $D'k_\Omega \simeq k_{\overline{\Omega}}$. The proof of \eqref{ii} is similar.
\\
\qed

Let $\sigma$ be a section of $T^*X \to X$ and consider the
following commutative diagram with cartesian square
\begin{equation}\label{TXT*X}
\xymatrix{TX & \\
T^*X \times_X TX \ar[u]^{\pi_2} \ar[d]_{\pi_1} & TX
\ar[l]^{\ \ \ \ \ \ \sigma'} \ar[d]^{\tau_X} \ar[ul]_{\id} \\
T^*X & X \ar[l]^{\ \ \sigma}.}
\end{equation}
We set
\begin{eqnarray*}
P & := & \{((x,\xi),(x,v)) \in T^*X \times_X TX; \; \langle
\xi,v\rangle \geq 0\} \\
P' & := & \{((x,\xi),(x,v)) \in T^*X \times_X TX; \; \langle
\xi,v\rangle \leq 0\} \\
P_\sigma & := & \{(x,v) \in  TX; \; \langle
\sigma(x),v\rangle \geq 0\}=\sigma'{}^{-1}(P) \\
P'_\sigma & := & \{(x,v) \in  TX; \; \langle \sigma(x),v\rangle
\leq 0\}=\sigma'{}^{-1}(P')
\end{eqnarray*}
The kernel $K_\sigma$ is defined as follows
\begin{equation}\label{microker}
K_\sigma:=Rp_{!!}(k_{\overline{\Omega}} \otimes
\rho_!k_{P_\sigma}) \otimes \rho_!\delta_*\formX.
\end{equation}

\begin{prop}\label{relsigma} (i) Let $F \in D^b_{\rc}(k_X)$ and $G \in
D^b(k_{X_{sa}})$. There is a natural arrow
$$
\varphi:\rh(F,K_\sigma \circ G) \to \imin \sigma \muh^{sa}(F,G),
$$
where $K_\sigma \circ G$ means $Rq_{1!!}(\imin {q_2}G \otimes
K_\sigma)$.

(ii) Let $\rho:X \to X_{sa}$ be the natural functor of sites. Then
$\imin \rho(\varphi)$ is an isomorphism.
\end{prop}
\dim\ \ $\mathrm{(i)}_a$ We have the chain of isomorphisms
\begin{eqnarray*}
\imin \sigma \muh^{sa}(F,G) & \simeq & \imin \sigma
R\pi_{1!!}(\imin
{\pi_2}\nu^{sa}_\Delta\rh(\imin {q_2}F,q_1^!G) \otimes k_{P'}) \\
& \simeq &  R\tau_{X!!}\sigma'{}^{-1}(\imin
{\pi_2}\nu^{sa}_\Delta\rh(\imin {q_2}F,q_1^!G) \otimes k_{P'}) \\
& \simeq &  R\tau_{X!!}(\nu^{sa}_\Delta\rh(\imin {q_2}F,q_1^!G)
\otimes k_{P'_\sigma}).
\end{eqnarray*}
Consider the normal deformation of $\Delta$ in $X \times X$
visualized by the diagram \eqref{normdef}. We have
\begin{eqnarray*}
\nu^{sa}_\Delta\rh(\imin {q_2}F,q_1^!G) & \simeq & \imin s
\mathrm{R}\Gamma_\Omega\imin p \rh(\imin {q_2} F,q_1^!G) \\
& \simeq & \imin s  \rh(\imin {p_2} F,\mathrm{R}\Gamma_\Omega\imin p q_1^!G) \\
& \simeq & \imin s  \rh(\imin {p_2} F,\mathrm{R}\Gamma_\Omega\imin {p_1} G) \otimes \imin \tau_X\omega_X \\
& \simeq & \imin s  \rh(\imin {p_2} F,\mathrm{R}\Gamma_\Omega\imin {p_1} G) \otimes \imin \tau_X\formX \\
& \simeq & \imin s  \rh(\imin {p_2} F,\imin {p_1} G \otimes
k_{\overline{\Omega}}) \otimes \imin \tau_X\formX,
\end{eqnarray*}
 where the second isomorphism
follows from Lemma \ref{omegap} and the last one follows from
Lemma \ref{omegabar}. Hence we get
\begin{eqnarray*}
\imin \sigma \muh^{sa}(F,G) & \simeq & R\tau_{X!!}(\imin s
\rh(\imin {p_2} F,\imin {p_1} G
\otimes k_{\overline{\Omega}}) \otimes \imin \tau_X\formX \otimes k_{P'_\sigma})\\
& \simeq & Rp_{2!!}s_{!!}(\imin s  \rh(\imin {p_2} F,\imin {p_1} G
\otimes k_{\overline{\Omega}}) \otimes \imin \tau_X\formX \otimes k_{P'_\sigma})\\
& \simeq & Rp_{2!!}(\rh(\imin {p_2} F,\imin {p_1} G \otimes
k_{\overline{\Omega}}) \otimes \imin p\delta_*\formX \otimes
k_{P'_\sigma}).
\end{eqnarray*}

$\mathrm{(i)}_b$ On the other hand we have the chain of
isomorphisms
\begin{eqnarray*}
K_\sigma \circ G & \simeq & Rq_{1!!}(\imin {q_2}G \otimes
Rp_{!!}(k_{\overline{\Omega}} \otimes \rho_!k_{P_\sigma}) \otimes
\rho_!\delta_*\formX) \\
& \simeq & Rp_{1!!}(\imin {p_2}G \otimes k_{\overline{\Omega}}
\otimes \rho_!k_{P_\sigma} \otimes
\imin p \rho_!\delta_*\formX) \\
& \simeq & Rp_{2!!}(\imin {p_1}G \otimes k_{\overline{\Omega}}
\otimes \rho_!k_{P'_\sigma} \otimes
\imin p \rho_!\delta_*\formX) \\
& \simeq & Rp_{2!!}(\imin {p_1}G \otimes k_{\overline{\Omega}}
\otimes \rho_!(k_{P'_\sigma} \otimes \imin p \delta_*\formX)).
\end{eqnarray*}
Hence we get
\begin{eqnarray*}
\rh(F,K_\sigma \circ G) & \simeq & \rh(F,Rp_{2!!}(\imin {p_1}G
\otimes k_{\overline{\Omega}} \otimes \rho_!(k_{P'_\sigma} \otimes
\imin p
\delta_*\formX))) \\
& \simeq & Rp_{2!!}\rh(\imin {p_2}F,\imin {p_1}G \otimes
k_{\overline{\Omega}} \otimes \rho_!(k_{P'_\sigma} \otimes \imin p
\delta_*\formX)) \\
& \simeq & Rp_{2!!}(\rh(\imin {p_2}F,\imin {p_1}G \otimes
k_{\overline{\Omega}}) \otimes \rho_!(k_{P'_\sigma} \otimes \imin
p \delta_*\formX)).
\end{eqnarray*}

$\mathrm{(i)}_c$ The adjunction morphism defines a morphism
$\rho_! \to \rho_*$. It induces the morphism
$$\varphi:\rh(F,K_\sigma \circ G) \to \imin \sigma \muh^{sa}(F,G).
$$ \newline

(ii) Composing with $\imin \rho$ we get $\imin \rho \circ \rho_!
\iso \imin \rho \circ \rho_* \simeq \id$. Hence we get
$$\imin \rho(\varphi):\imin \rho\rh(F,K_\sigma \circ G) \iso \imin \rho \imin \sigma \muh^{sa}(F,G).
$$
\nopagebreak \newline \qed

Let $\pi:T^*X \to  X$ be the projection and consider the canonical
1-form $\omega$, the restriction to the diagonal of the map
${}^{t}\hspace{-0.4mm}\pi'$. Replace $X$ with $T^*X$ and $\sigma$
with $\omega$ in \eqref{microker} and consider the microlocal
kernel
$$K_{T^*X}=Rp_{!!}(k_{\overline{\Omega}} \otimes
\rho_!k_{P_{\omega}}) \otimes
\rho_!\delta_*\omega^{\otimes-1}_{\Delta_{T^*X}|T^*X \times T^*X}
$$

\begin{df} The functor of microlocalization of \cite{KSIW} is
defined as
\begin{eqnarray*}
\mu:D^b(k_{X_{sa}}) & \to & D^b(k_{T^*X_{sa}}) \\
F & \mapsto & \mu F=K_{T^*X} \circ \imin \pi F
\end{eqnarray*}
\end{df}

\begin{teo}\label{relmu} (i) Let $F \in D^b_{\rc}(k_X)$ and $G \in
D^b(k_{X_{sa}})$. There is a natural arrow
\begin{equation}\label{morphmu}
\varphi:\rh(\imin \pi F,\mu G) \to \muh^{sa}(F,G).
\end{equation}

(ii) Let $\rho:T^*X \to T^*X_{sa}$ be the natural functor of sites. Then
$\imin \rho(\varphi)$ is an isomorphism.
\end{teo}
\dim\ \ (i) By Lemma \ref{omegamu} and Proposition \ref{relsigma}
(i) we get the morphisms
$$\muh^{sa}(F,G) \simeq \imin {\omega} \muh^{sa}(\imin \pi F,\imin \pi G) \gets \rh(\imin \pi F,\mu
G).$$

(ii) The result follows from Proposition \ref{relsigma} (ii). \\
\qed

\begin{es} The morphism \eqref{morphmu} is not an isomorphism in
general. For example let $F \in \mod(k_{X_{sa}})$. Then
$$\rh(\imin \pi k_X,\mu F) \simeq \rh(k_{T^*X},\mu F) \simeq \mu
F,$$

on the other hand we have

$$\muh^{sa}(k_X,F) \simeq j_*\mu_XF \simeq j_*F,$$
where $j:T^*_XX \hookrightarrow T^*X$ denotes the embedding of the
zero section.
\end{es}

\section{Holomorphic functions with growth conditions}

In this section we show how the functors we defined before generalize classical constructions.

In $\S$\,\ref{asymptotic} we show the relation between specialization of Whitney holomorphic functions with the functor of formal specialization of \cite{Co01}, and the sheaf of asymptotically developable functions of \cite{Ma91} and \cite{Si90}.

In $\S$\,\ref{Andronikof-Colin} we study the microlocalization of
tempered and Whitney holomorphic functions. We establish a relation with the functors of
tempered and formal microlocalization introduced by Andronikof in
\cite{An94} and Colin in \cite{Co98}.

\subsection{Review on temperate and formal cohomology}

From now on, the base field is $\CC$. Let $M$ be a real analytic
manifold. One denotes by $\db_M$ the sheaf of Schwartz's
distributions, by $\C^\infty_M$
the sheaf of $\C^\infty$-functions, and by $\D_M$ the sheaf of finite order
differential operators with analytic coefficients. As usual, for $F \in \mod(\CC_M)$,
we set $D'F=\rh(F,\CC_M)$.\\

\noindent In \cite{Ka84}
the author defined the functor
$$\tho(\cdot,\db_M):\mod_{\rc}(\CC_M) \to \mod(\D_M)$$
in the following way: let $U$ be a subanalytic subset of $M$ and
$Z=M \setminus U$. Then the sheaf $\tho(\CC_U,\db_M)$ is defined
by the exact sequence
$$\exs{\Gamma_Z\db_M}{\db_M}{\tho(\CC_U,\db_M)}.$$
This functor is exact and extends as functor in the derived
category, from $D^b_{\rc}(\CC_M)$ to $D^b(\D_M)$. Moreover the
sheaf $\tho(F,\db_M)$ is soft for any $\R$-constructible sheaf
$F$.

\begin{df} Let $Z$ be a closed subset of $M$. We denote by
$\II^\infty_{M,Z}$ the sheaf of $\C^\infty$-functions on $M$
vanishing up to infinite order on $Z$.
\end{df}

\begin{df} A Whitney function on a closed subset $Z$ of $M$ is an indexed
family $F=(F^k)_{k\in \N^n}$ consisting of continuous functions on
$Z$ such that $\forall m \in \N$, $\forall k \in \N^n$, $|k| \leq
m$, $\forall x \in Z$, $\forall \varepsilon >0$ there exists a
neighborhood $U$ of $x$ such that $\forall y,z \in U \cap Z$

$$\left|F^k(z)-\sum_{|j+k|\leq m}{(z-y)^j \over
j!}F^{j+k}(y)\right| \leq \varepsilon d(y,z)^{m-|k|}.$$ We denote
by $W^\infty_{M,Z}$ the space of  Whitney $\C^\infty$-functions on
$Z$. We denote by $\W^\infty_{M,Z}$ the sheaf $U \mapsto
W_{U,U\cap Z}^\infty$.
\end{df}
In \cite{KS96} the authors defined the functor
$$\cdot \wtens \C^\infty_M:\mod_{\rc}(\CC_M) \to \mod(\D_M)$$
in the following way: let $U$ be a subanalytic open subset of $M$
and $Z=M \setminus U$. Then $\CC_U \wtens
\C^\infty_M=\II^\infty_{M,Z}$, and $\CC_Z \wtens
\C^\infty_M=\W^\infty_{M,Z}$. This functor is exact and extends as
a functor in the derived category, from $D^b_{\rc}(\CC_M)$ to
$D^b(\D_M)$. Moreover the sheaf $F \wtens \C^\infty_M$ is soft for
any $\R$-constructible sheaf $F$.\\

Now let $X$ be a complex manifold, $X_\R$ the underlying real
analytic manifold and $\overline{X}$ the complex conjugate
manifold. The product $X \times \overline{X}$ is a
complexification of $X_\R$ by the diagonal embedding $X_\R
\hookrightarrow X \times \overline{X}$. One denotes by $\OO_X$ the
sheaf of holomorphic functions and by $\D_X$ the sheaf of finite
order differential operators with holomorphic coefficients. For $F
\in D^b_{\rc}(\CC_X)$ one sets
\begin{eqnarray*}
\tho(F,\OO_X) & = & \rh_{\D_{\overline{X}}}(
\OO_{\overline{X}},\tho(F,\db_{X_\R})), \\
F \wtens \OO_X & = & \rh_{\D_{\overline{X}}}(\OO_{\overline{X}},F \wtens \C^\infty_{X_\R}),
\end{eqnarray*}
 and these functors are
called the functors of temperate and formal cohomology respectively.

\subsection{Tempered and Whitney holomorphic functions}

\begin{df} One denotes by $\dbt_M$ the presheaf of tempered distributions on
$M_{sa}$ defined as follows:
$$U \mapsto \Gamma(M;\db_M)/\Gamma_{M\setminus U}(M;\db_M).$$
\end{df}

 As a consequence of the \L ojasievicz's inequalities
\cite{Lo59}, for $U,V \in \op(M_{sa})$ the sequence
$$\exs{\dbt_M(U \cup
V)}{\dbt_M(U)\oplus\dbt_M(V)}{\dbt_M(U \cap V)}$$ is exact. Then
$\dbt_M$ is a sheaf on $M_{sa}$. Moreover, by definition $\dbt_M$ is quasi-injective.\\

\begin{df} One denotes by $\CWM$ the presheaf of Whitney
$\C^\infty$-functions on $M_{sa}$ defined as follows:
$$U \mapsto \Gamma(M;H^0D'\CC_U \wtens \C^\infty_M).$$
\end{df}
 As a consequence of a result of \cite{Ma67}, for $U,V \in \op(M_{sa})$ the sequence
$$\lexs{\CWM(U \cup
V)}{\CWM(U)\oplus\CWM(V)}{\CWM(U \cap V)}$$ is exact. Then $\CWM$
is a sheaf on $M_{sa}$. Moreover if $U \in \op(M_{sa})$ is locally cohomologically trivial (l.c.t. for short), i.e. if $D'\CC_U \simeq \CC_{\overline{U}}$, the morphism $\Gamma(X;\CWM) \to \Gamma(U;\CWM)$ is surjective and $\mathrm{R}\Gamma(U;\CWM)$ is concentrated in degree zero.

We have the
following result (see \cite{KS01}, \cite{Pr1}).

\begin{prop} For each $F \in \mod_{\rc}(\CC_M)$ one has the isomorphisms
\begin{eqnarray*}
\imin \rho \rh(F,\dbt_M) & \simeq & \tho(F,\db_M), \\
\imin \rho\rh(F,\CWM) & \simeq & D'F\wtens\C^\infty_M.
\end{eqnarray*}
\end{prop}

Now let $X$ be a complex manifold, $X_\R$ the underlying real
analytic manifold and $\overline{X}$ the complex conjugate
manifold. One denotes by $\ot_X$ and $\OWX$ the sheaves of tempered and Whitney holomorphic
functions defined as follows:
\begin{eqnarray*}
\ot_X & := & \rh_{\rho_!\D_{\overline{X}}}(\rho_!\OO_{\overline{X}},\dbtxr), \\
\OWX & := & \rh_{\rho_!\D_{\overline{X}}}(\rho_!\OO_{\overline{X}},\CWXR).
\end{eqnarray*}

The relation with the functors of temperate and formal cohomology are given by
the following result (see \cite{KS01}, \cite{Pr1})

\begin{prop} For each $F \in D^b_{\rc}(\CC_X)$ one has the isomorphisms
\begin{eqnarray*}
\tho(F,\OO_X) & \simeq & \imin \rho \rh(F,\ot_X), \\
D'F \wtens \OO_X & \simeq & \imin \rho \rh(F,\OWX).
\end{eqnarray*}
\end{prop}

\subsection{Asymptotic expansions}\label{asymptotic}

Let $M$ be a real analytic manifold. We consider a slight generalization of the sheaf of Whitney $\C^\infty$-functions of \cite{KS01}.

\begin{df} Let $F \in \mod_{\rc}(\CC_M)$ and let $U \in \op(M_{sa})$. We define the presheaf $\CW_{M|F}$ as follows:
$$U \mapsto \Gamma(M;H^0D'\CC_U \otimes F \wtens \C^\infty_M).$$
\end{df}
Let $U,V \in \op(M_{sa})$, and consider the exact sequence
$$\exs{\CC_{U\cap V}}{\CC_U \oplus \CC_V}{\CC_{U\cup V}},$$
applying the functor $\ho(\cdot,\CC_M)=H^0D'(\cdot)$ we obtain
$$\lexs{H^0D'\CC_{U\cap V}}{H^0D'\CC_U \oplus H^0D'\CC_V}{H^0D'\CC_{U\cup V}},$$
applying the exact functors $\cdot \otimes F$, $\cdot\wtens \C^\infty_M$ and taking global sections
we obtain
$$\lexs{\CW_{M|F}(U\cup V)}{\CW_{M|F}(U) \oplus \CW_{M|F}(V)}{\CW_{M|F}(U\cap V)}.$$
This implies that $\CW_{M|F}$ is a sheaf on $M_{sa}$. Moreover if $U \in \op(M_{sa})$ is l.c.t., the morphism $\Gamma(X;\CWM) \to \Gamma(U;\CW_{M|F})$ is surjective and $\mathrm{R}\Gamma(U;\CW_{M|F})$ is concentrated in degree zero. Let $\exs{F}{G}{H}$ be an exact sequence in $\mod_{\rc}(\CC_M)$, we obtain an exact sequence in $\mod(\CC_{M_{sa}})$
\begin{equation}\label{exsFGH}
\exs{\CW_{M|F}}{\CW_{M|G}}{\CW_{M|H}}.
\end{equation}

We can easily extend the sheaf $\CW_{M|F}$ to the case of $F \in D^b_{\rc}(\CC_M)$, taking a finite resolution of $F$ consisting of locally finite sums $\oplus \CC_V$, with $V$ l.c.t. in $\op^c(M_{sa})$. In fact, the sheaves $\CW_{M|\oplus \CC_V}$ form a complex quasi-isomorphic to $\CW_{M|F}$ consisting of acyclic objects with respect to $\Gamma(U;\cdot)$, where $U$ is l.c.t. in $\op^c(M_{sa})$.

As in the case of Whitney $\C^\infty$-functions one can prove that, if $G \in D^b_{\rc}(\CC_M)$ one has
$$\imin \rho \rh(G,\CW_{F|M}) \simeq D'G \otimes F \wtens \C^\infty_M.$$


\begin{es} Setting $F=\CC_M$ we obtain the sheaf of Whitney $\C^\infty$-functions. Let $N$ be a closed analytic submanifold of $M$. Then $\CW_{M|\CC_{M\setminus N}}$ is the sheaf of Whitney $\C^\infty$-functions vanishing on $N$ with all their derivatives.
\end{es}

\begin{nt} Let $Z$ be a locally closed subanalytic subset of $M$. We set for short $\CW_{M|Z}$ instead of $\CW_{M|\CC_Z}$.
\end{nt}

Let $N$ be a closed analytic submanifold of $M$, let $T_NM \stackrel{\tau}{\to} N$ be the normal vector bundle and consider the normal deformation $\widetilde{M}_N$ as in $\S$\,\ref{rewdef}.

Set $F=\CC_{M\setminus N}$, $G=\CC_M$, $H=\CC_N$ in \eqref{exsFGH}. The exact sequence
$$\exs{\CW_{M|M\setminus N}}{\CWM}{\CW_{M|N}}$$
induces an exact sequence
$$\exs{\nu^{sa}_N\CW_{M|M\setminus N}}{\nu^{sa}_N\CWM}{\nu^{sa}_N\CW_{M|N}},$$
in fact let $V$ be a l.c.t. conic subanalytic subset of $\widetilde{M}_N$ and $U \in \op(M_{sa})$ such that
$C_M(X \setminus U) \cap V=\varnothing$, then we can find a l.c.t. $U'\subset U$ satisfying the same property. Moreover it is easy to see that $\nu^{sa}_N\CW_{M|N} \simeq \imin \tau \CW_{M|N}$, hence we get the exact sequence
\begin{equation}\label{asym}
\exs{\nu^{sa}_N\CW_{M|M\setminus N}}{\nu^{sa}_N\CWM}{\imin \tau \CW_{M|N}}.
\end{equation}

Now let us study the relation with the constructions of \cite{Co01}. In that work the author defined the functors of Whitney specialization as follows: let $F \in D^b_{\rc}(k_M)$, then
$$
w\nu_N(F,\C^\infty_M) = \imin s \rh_{\D_{\widetilde{X}}}(\dd {\widetilde{M}_N}{}{M},(\imin p F)_\Omega \wtens \C^\infty_M).
$$
The stalks are given by the following formula: let $v \in T_NM$. Then
$$
H^k(w\nu_N(F,\C^\infty_M))_v \simeq \lind U H^k(X;F_{\overline{U}} \wtens \C^\infty_M),
$$
where $U \in \op(M_{sa})$ l.c.t. such that
$v \notin C_M(X \setminus U)$. By Theorem \ref{4.2.3} (ii) it turns out that
$$
w\nu_N(F,\C^\infty_M) \simeq \imin \rho \nu^{sa}_N\CW_{M|F},
$$
this means that Whitney specialization is obtained by specializing the sheaf $\CW_{M|F}$.\\

Assume that $M \simeq \{(x,y) \in \R^\ell \times \R^{n-\ell}\}$ and $N \simeq \{0\} \times \R^{n-\ell}$. A sector $S$ of $M$ is a subanalytic open subset $S=U \times V$ with $U \in \op(\R^{n-\ell}_{sa})$ and $V=W \cap B(0,\varepsilon)$, where $W \in \op(\R^\ell_{sa,\RP})$ and $B(0,\varepsilon)$ is the open ball of center 0 and radius $\varepsilon$. We say that $S'$ is a subsector of $S$ if $\overline{S'} \setminus N \subset S$ or, with the notations of Definition \ref{2.2.2} if $\RP S' \subset\subset \RP S$ with the conic topology. We write for short $S' < S$.

\begin{df} Let $S$ be an open sector of $M$ and let $f \in \C^\infty_M$. One says that $f$ is asymptotically developable along $M$, if there exists a formal series $\sum_{k \in \N^\ell}a_k(x)y^k$ with $\C^\infty$ coefficients $a_k$ such that,
for all $S' < S$, $m \in \N$, there exists $C>0$ such that
$$
\forall (x,y) \in S', \ \ \left|f(x,y)-\sum_{|k|\leq m}a_k(x)y^k\right| \leq C\|y\|^m.
$$
One denotes by $\sigma_N(S) \subset \C^\infty_M$ the space of functions asymptotically developable along $M$, and by $\sigma_M^\infty=\{f \in \C^\infty_N(S),\; \forall k \in \N^\ell,\; D^kf \in \sigma_M(S)\}$.
\end{df}

We have the following result (see Proposition 2.10 of \cite{Co01}).

\begin{prop} Let $S$ be a sector of $M$. Then $\Gamma(\RP S;\imin \rho \nu^{sa}_N\CWM) \simeq \sigma^\infty_N(S)$ and $\Gamma(\RP S;\imin \rho \nu^{sa}_N\CW_{M|M\setminus N})$ is the subspace of functions asymptotically developable to the identically zero series.
\end{prop}

Applying the functor $\imin \rho$ to the exact sequence \eqref{asym} we obtain the exact sequence
$$
\exs{\imin \rho\nu^{sa}_N\CW_{M|M\setminus N}}{\imin \rho\nu^{sa}_N\CWM}{\imin \rho \imin \tau \CW_{M|N}},
$$
where the surjective arrow is the map which associates to a function its asymptotic expansion.\\

Let $X$ be a complex manifold and let $Z$ be a complex submanifold of $X$. Let $F \in D^b_{\rc}(\CC_X)$. We denote by $\OW_{X|F}$ the sheaf defined as follows:
$$
\OW_{X|F} := \rh_{\rho_!\D_{\overline{X}}}(\rho_!\OO_{\overline{X}},\CW_{X_\R|F}).
 $$
Let $\exs{F}{G}{H}$ be an exact sequence in $\mod_{\rc}(\CC_X)$. Then the exact sequence \eqref{exsFGH} gives rise to the distiguished triangle
\begin{equation}\label{dtFGH}
\dt{\OW_{X|F}}{\OW_{X|G}}{\OW_{X|H}}.
\end{equation}
If we consider the functor of specialization of formal cohomology of \cite{Co01}
$$
w\nu_Z(F,\OO_X) = \rh_{\imin \tau \D_{\overline{X}}}(\imin \tau \OO_{\overline{X}},w\nu_Z(F,\C^\infty_{X_\R})),
$$
we have the isomorphism
$$
w\nu_Z(F,\OO_X) \simeq \imin \rho \nu^{sa}_Z\OW_{X|F}.
$$
Setting $F=\CC_{X\setminus Z}$, $G=\CC_X$, $H=\CC_Z$ in \eqref{dtFGH} and applying the functor of specialization, we have the distinguished triangle
\begin{equation}\label{dtColin}
\dt{\imin \rho \nu^{sa}_Z\OW_{X|X \setminus Z}}{\imin \rho \nu^{sa}_Z\OWX}{\imin \rho \imin \tau \OO_{X|Z}}.
\end{equation}
The sheaves $\imin \rho \nu^{sa}_Z\OWX$ and $\imin \rho \imin \tau \OO_{X|Z}$ are concentrated in degree zero. This follows from the following result of \cite{DS96}: if $U \in \op(X_{sa})$ is convex, then $\mathrm{R}\Gamma(X;\CC_{\overline{U}} \wtens \OO_X)$ is concentrated in degree zero.
Moreover (see \cite{KS96}) the sheaf $\imin \rho \OW_{X|Z}$ is isomorphic to the sheaf $\OO_{X}\widehat{|}_Z$, the formal completion of $\OO_X$ along $Z$. We have an exact sequence
$$
0 \to \imin \rho H^0\nu^{sa}_Z\OW_{X|X \setminus Z} \to \imin \rho \nu^{sa}_Z\OO_X \to \imin \tau \OO_{X}\widehat{|}_Z \to \imin \rho H^1\nu^{sa}_Z\OW_{X|X\setminus Z} \to 0.
$$

Let $\sigma^{hol}_Z(S)$ be the space of holomorphic function asymptotically developable in $S$, having an asymptotic expansion with holomorphic coefficients. We have the following results of \cite{Co01}.

\begin{prop}  Let $S$ be a sector of $M$. Then $\Gamma(\RP S;\imin \rho \nu^{sa}_Z\OWX) \simeq \sigma^{hol}_Z(S)$ and $\Gamma(\RP S;\imin \rho \nu^{sa}_Z\OW_{X|X\setminus Z})$ is the subspace of functions asymptotically developable to the identically zero series.
\end{prop}

\begin{prop} The distinguished triangle \eqref{dtColin} induces an exact sequence outside the zero section
\begin{equation}\label{exspT}
\exs{\imin \rho H^0 \nu^{sa}_Z\OW_{X|Z}|_{\pT_ZX}}{\imin \rho \nu^{sa}_Z\OWX|_{\pT_ZX}}{\imin \pittau \OO_{X}\widehat{|}_Z}.
\end{equation}
On the zero section we have the exact sequence
\begin{equation}\label{exsZ}
\exs{\OO_X|_Z}{\OO_{X}\widehat{|}_Z}{\imin \rho H^1\nu_Z^{sa}\OW_{X|X\setminus Z}|_Z}.
\end{equation}
\end{prop}
Remark that on the exact sequence \eqref{exsZ} we used Theorem \ref{4.2.3} (iii) and the fact that $\imin \rho \OWX \simeq \OO_X$.

\begin{es} Set $X=\CC$ and $Z=\{0\}$. Outside the zero section the sheaves $\imin \rho \nu^{sa}_Z\OW_{X|X \setminus Z}$ and $\imin \rho \nu^{sa}_Z\OWX$ are the well known sheaves $\A_0$ and $\A$ of Malgrange \cite{Ma91} and Sibuya \cite{Si90}. These sheaves were defined in the real blow up of the origin of $\CC$ identified with $\mathbb{S}^1 \times (\RP \cup \{0\})$. Let $\pi$ be the projection on $\CC$. The sequence \eqref{exspT} is a generalization of the exact sequence in $\mod(\CC_{\mathbb{S}^1})$
$$
\exs{\A_0}{\A}{\imin \pi \CC[[z]]},
$$
and the sequence \eqref{exsZ} is a generalization of the exact sequence
$$
\exs{\CC\{z\}}{\CC[[z]]}{H^1(\mathbb{S}^1;\A_0)}.
$$
\end{es}

\subsection{Microlocalization of $\OO^t_X$ and $\OWX$}\label{Andronikof-Colin}

Let $f:X \to Y$ be a smooth morphism of real analytic manifolds. We have
the following results (see the Appendix):
\begin{eqnarray*}
\rh_{\rho_!\D_X}(\rho_!\ddxy,\dbt_X) \simeq \imin f \dbt_Y, \\
\rh_{\rho_!\D_X}(\rho_!\ddxy,\CW_X) \simeq \imin f \CW_Y.
\end{eqnarray*}


Let us consider the normal deformation of the diagonal in $X
\times X$ of diagram \eqref{normdef}. Let $F \in
D^b_{\rc}(\CC_X)$. We recall the definitions of the Andronikof's
functor of microlocalization of tempered distributions
$$
t\muh(F,\db_X):=(\imin s (\D_{\widetilde{X \times
X}\stackrel{p_1}{\gets}X} \otimes_{\D_{\widetilde{X \times X}}}
\tho((\imin {p_2}F)_{\Omega},\db_{\widetilde{X \times
X}})[-1]))^\wedge
$$
and Colin's microlocalization of the Whitney tensor product
$$
F \underset{\mu}{\wtens} \C^\infty_X := (\imin s \rh( \dd {\widetilde{X \times
X}}{p_1}{X},(\imin p_2F)_{\overline{\Omega}} \wtens
\C^\infty_{\widetilde{X \times X}}))^\vee.
$$

\begin{teo}\label{muhomtwr} Let $F \in D^b_{\rc}(\CC_X)$. We have the isomorphisms
\begin{eqnarray}
\imin \rho \muh^{sa}(F,\dbt_X) & \simeq &
t\muh(F,\db_X), \label{Andronikofr} \\
\imin \rho \muh^{sa}(F,\CW_X) & \simeq &
(D'F \underset{\mu}{\wtens} \C^\infty_X)^a, \label{Colinr}
\end{eqnarray}
where $(\cdot)^a$ denotes
the direct image of the antipodal map.
\end{teo}
\dim\ \ Let $G \in D^b(\CC_{X_{sa}})$. We have
\begin{eqnarray*}
\muh^{sa}(F,G) & \simeq & (\imin s
\mathrm{R}\Gamma_\Omega\imin p \rh(\imin
{q_2}F,q_1^!G))^\wedge \\
& \simeq & (\imin s \mathrm{R}\Gamma_\Omega\rh(\imin p\imin
{q_2}F,
\imin pq_1^!G))^\wedge \\
& \simeq & (\imin s \rh((\imin p\imin {q_2}F)_\Omega,
\imin pq_1^!G))^\wedge \\
& \simeq & (\imin s (\rh((\imin p \imin {q_2}F)_\Omega,\imin p
\imin {q_1}G)\otimes \imin pq_1^!\CC_X))^\wedge \\
& \simeq & (\imin s \rh((\imin {p_2}F)_\Omega,\imin {p_1}
G)\otimes \imin s\imin pq_1^!\CC_X)^\wedge,
\end{eqnarray*}
where the second isomorphism follows from Lemma \ref{omegap}.\\

\noindent (i) Let us prove \eqref{Andronikofr}. Setting $G=\dbt_X$ and composing with $\imin \rho$ we have
\begin{eqnarray*}
\lefteqn{\imin \rho \muh^{sa}(F,\dbt_X)    } \\
& \simeq & (\imin s \imin \rho \rh((\imin {p_2}F)_\Omega,\imin
{p_1}\dbt_X) \otimes \imin s \imin
p q_1^!\CC_X)^\wedge  \\
& \simeq &  (\imin s \imin \rho \rh_{\rho_!\D_{\widetilde{X \times
X}}}(\rho_!\D_{\widetilde{X \times
X}\stackrel{p_1}{\to}X},\rh((\imin
{p_2}F)_{\Omega},\dbt_{\widetilde{X \times X}})) \otimes \imin
s\imin
p q_1^!\CC_X)^\wedge \\
& \simeq & (\imin s \rh_{\D_{\widetilde{X \times
X}}}(\D_{\widetilde{X \times X}\stackrel{p_1}{\to}X},\imin \rho
\rh((\imin {p_2}F)_{\Omega},\dbt_{\widetilde{X \times X}})) \otimes
\imin s\imin
p q_1^!\CC_X)^\wedge  \\
& \simeq &  (\imin s \rh_{\D_{\widetilde{X \times
X}}}(\D_{\widetilde{X \times X}\stackrel{p_1}{\to}X},\tho((\imin
{p_2}F)_{\Omega},\db_{\widetilde{X \times X}})) \otimes \imin
s\imin
p q_1^!\CC_X)^\wedge  \\
& \simeq & (\imin s (\ddual {X}{p_1}{\widetilde{X \times X}}
\otimes_{\D_{\widetilde{X \times X}}} \tho((\imin
{p_2}F)_{\Omega},\db_{\widetilde{X \times X}})[-1]))^\wedge
\\
& \simeq & t\muh(F,\db_X).
\end{eqnarray*}

\smallskip

\noindent (ii) Let us prove \eqref{Colinr}. Setting $G=\CW_X$ and composing with $\imin \rho$ we have
\begin{eqnarray*}
\lefteqn{\imin \rho \muh^{sa}(F,\CW_X)} \\
& \simeq & (\imin s\imin \rho
\rh((\imin p_2F)_\Omega,\imin {p_1}\CW_X) \otimes \imin s\imin pq_1^!\CC_X)^\wedge \\
& \simeq & (\imin s\imin \rho
\rh((\imin p_2F)_\Omega,\imin {p_1}\CW_X))^{\vee a} \\
& \simeq & (\imin s\imin \rho
\rh_{\rho_!\D_{\widetilde{X \times
X}}}(\rho_! \dd {\widetilde{X \times X}}{p_1}{X},\rh((\imin p_2F)_\Omega,\CW_{\widetilde{X \times X}})) )^{\vee a} \\
& \simeq & (\imin s
\rh_{\D_{\widetilde{X \times
X}}}( \dd {\widetilde{X \times X}}{p_1}{X},\imin \rho\rh((\imin p_2F)_\Omega,\CW_{\widetilde{X \times X}})) )^{\vee a} \\
& \simeq & (\imin s
\rh_{\D_{\widetilde{X \times
X}}}( \dd {\widetilde{X \times X}}{p_1}{X},D'(\imin p_2F)_\Omega \wtens \C^\infty_{\widetilde{X \times X}}) )^{\vee a} \\
& \simeq & (\imin s \rh_{\D_{\widetilde{X \times
X}}}( \dd {\widetilde{X \times
X}}{p_1}{X},(\imin p_2D'F)_{\overline{\Omega}} \wtens
\C^\infty_{\widetilde{X \times X}}) )^{\vee a}
\end{eqnarray*}
where
the last
isomorphism follows since
$$D'((\imin {p_2}F)_\Omega) \simeq \mathrm{R}\Gamma_\Omega
D'(\imin p_2 F) \simeq \mathrm{R}\Gamma_\Omega \imin p_2D'F 
\simeq (\imin
{p_2}D'F)_{\overline{\Omega}}.$$ Here we used Lemma \ref{omegabar} and the fact that $p_2$ is smooth.
We have
$$
\imin s \rh_{\D_{\widetilde{X \times
X}}}( \dd {\widetilde{X \times X}}{p_1}{X},(\imin
p_2D'F)_{\overline{\Omega}} \wtens \C^\infty_{\widetilde{X \times X}})
)^{\vee a} =  (D'F \underset{\mu}{\wtens} \C^\infty_X)^a
$$
and the result follows.\\
\qed

Let $X$ be a complex manifold and let $F \in D^b_{\rc}(\CC_X)$. In \cite{An94} and \cite{Co98} the authors constructed the functors $t\muh(F,\OO_X)$ of tempered microlocalization and $F \underset{\mu}{\wtens} \OO_X$ of formal microlocalization taking the Dolbeaut resolutions of the real ones.

\begin{teo}\label{muhomtw} Let $F \in D^b_{\rc}(\CC_X)$. We have the isomorphisms
\begin{eqnarray}
\imin \rho \muh^{sa}(F,\ot_X) & \simeq &
t\muh(F,\OO_X), \label{Andronikof} \\
\imin \rho \muh^{sa}(F,\OWX) & \simeq &
(D'F \underset{\mu}{\wtens} \OO_X)^a, \label{Colin}
\end{eqnarray}
where $(\cdot)^a$ denotes
the direct image of the antipodal map.
\end{teo}
\dim\ \ The result follows by taking Dolbeaut resolutions on the left and the right-hand side of \eqref{Andronikofr} and \eqref{Colinr}. Let us see the proof of \eqref{Andronikof}. Let $F \in D^b_{\rc}(\CC_X)$. We have the chain of isomorphisms
\begin{eqnarray*}
\lefteqn{\muh^{sa}(F,\ot_X)} \\
& \simeq & (\imin s \rh((\imin {p_2}F)_\Omega,p_1^!\rh_{\rho_!\D_{\overline{X}}}(\rho_!\OO_{\overline{X}},\dbt_{X_\R})))^\wedge[-1] \\
& \simeq & (\imin s \rh((\imin {p_2}F)_\Omega,\rh_{\imin {p_1}\rho_!\D_{\overline{X}}}(\imin {p_1} \rho_!\OO_{\overline{X}},p_1^!\dbt_{X_\R})))^\wedge[-1] \\
& \simeq & (\imin s \rh_{\imin {p_1}\rho_!\D_{\overline{X}}}(\imin {p_1} \rho_!\OO_{\overline{X}},\rh((\imin {p_2}F)_\Omega,p_1^!\dbt_{X_\R})))^\wedge[-1] \\
& \simeq & (\rh_{\imin s\imin {p_1}\rho_!\D_{\overline{X}}}(\imin s\imin {p_1} \rho_!\OO_{\overline{X}},\imin s\rh((\imin {p_2}F)_\Omega,p_1^!\dbt_{X_\R})))^\wedge[-1] \\
& \simeq & (\rh_{\imin {\tau_X}\rho_!\D_{\overline{X}}}(\imin {\tau_X}\rho_!\OO_{\overline{X}},\imin s\rh((\imin {p_2}F)_\Omega,p_1^!\dbt_{X_\R})))^\wedge[-1] \\
& \simeq & \rh_{\imin \pi\rho_!\D_{\overline{X}}}(\imin \pi\rho_!\OO_{\overline{X}},(\imin s\rh((\imin {p_2}F)_\Omega,p_1^!\dbt_{X_\R}))^\wedge)[-1] \\
& \simeq & \rh_{\rho_!\D_{\overline{X}}}(\rho_!\OO_{\overline{X}},\muh^{sa}(F,\dbt_{X_\R})).
\end{eqnarray*}
Applying the functor $\imin \rho$ we have
\begin{eqnarray*}
\imin \rho \muh^{sa}(F,\ot_X) & \simeq & \imin \rho \rh_{\rho_!\D_{\overline{X}}}(\rho_!\OO_{\overline{X}},\muh^{sa}(F,\dbt_{X_\R})) \\ & \simeq & \rh_{\D_{\overline{X}}}(\OO_{\overline{X}},\imin \rho\muh^{sa}(F,\dbt_{X_\R})) \\
& \simeq & \rh_{\D_{\overline{X}}}(\OO_{\overline{X}},t\muh(F,\db_{X_\R})) \\
& \simeq & t\muh(F,\OO_X).
\end{eqnarray*}
The proof of \eqref{Colin} is similar.\\
\qed

\section{Integral transforms}

In this section we give some applications related to the microlocalization of subanalytic sheaves. We show the existence of a natural action of microlocal operators on tempered and formal microlocalization. We show also the invariance under contact transformation of tempered and formal microlocalization.

\subsection{$\E_X$-modules}

Let $X$ be a complex manifold of complex dimension $d_X$.
Following the notations of \cite{KS90} one sets
$\E^\R_X=H^{d_X}(\mu_\Delta\OO_{X \times X}^{(0,{d_X})})$. It is a
sheaf of rings over $T^*X$ and for each $F \in D^b(\CC_X)$, $j \in
\Z$ the sheaf $H^j\muh(F,\OO_X)$ is naturally endowed with a
structure of left $\E^\R_X$-module. The sheaf $\E^\R_X$ is called
the ring of microlocal operators on $X$. It contains a subring,
denoted by $\E_X$ and called the ring of (finite-order)
microdifferential operators. We will not recall all the properties
of this sheaf and refer to \cite{Sc85} for a detailed study.

In \cite{An94} the author introduced the sheaf $\E^{\R,f}_X$ of
tempered microdifferential operators as
$$\E^{\R,f}_{X} := H^{d_X}(t\muh(\CC_\Delta,\OO_{X \times X})
\underset{\OO_{X \times X}}{\ltens} \OO_{X \times X}^{(0,d_X)}).$$

\begin{prop} One has $\E^{\R,f}_X \simeq
\imin \rho H^{d_X}\mu^{sa}_\Delta\OO_{X \times X}^{t(0,d_X)}$.
\end{prop}
\dim\ \ There is a natural morphism
$$\muh^{sa}(\CC_\Delta,\OO_{X \times X})
\underset{\rho_!\OO_{X \times X}}{\ltens} \rho_!\OO_{X \times
X}^{(0,d_X)} \to \muh^{sa}(\CC_\Delta,\OO_{X \times
X}^{t(0,d_X)}).$$ Applying the exact functor $\imin \rho$ we
obtain the morphism
$$\E^{\R,f}_X \to
\imin \rho H^{d_X}\muh^{sa}(\CC_\Delta,\OO_{X \times
X}^{t(0,d_X)}) \simeq \imin \rho H^{d_X}\mu^{sa}_\Delta\OO_{X
\times X}^{t(0,d_X)} ,$$
which is an isomorphism on the fibers by Theorem \ref{4.3.2}.\\
\qed

Let us recall the following results:

\begin{itemize}
\item the sheaf $t\muh(\CC_\Delta,\OO_{X \times X})$ is
concentrated in degree $d_X$,
\item one has the ring inclusions $\E_X \subset \E^{\R,f}_X
\subset \E^\R_X$.
\end{itemize}

\subsection{Integral transforms}


 Let $X,Y,Z$ be three manifolds. Let
$q_{ij}$ be the $(i,j)$-th projection defined on $X \times Y
\times Z$ and let $p_{ij}$ be the $(i,j)$-th projection defined on
$T^*X \times T^*Y \times T^*Z$. Let $p^a_{ij}$ be the composition
of $p_{ij}$ with the antipodal map $a$ and let $j:X \times Y \times Y \times Z \to X \times Y \times Z$ be the diagonal embedding. Consider the following diagram
\begin{equation}\label{michint}
\xymatrix{T^*(X \times Y) \times T^*(Y \times Z) & T^*X \times T^*Y \times T^*Z \ar `[r] `[ddd]_{p_{13}} [ddd] \ar[l]_{\hspace{8mm} p_{12}^a \times p_{23}^a} \ar[d]^{\id \times p_2 \times a}_\wr & \\
\ar@{} [dr] |{\Box} T^*(X \times Y) \times_Y T^*(Y \times Z)  \ar[u]^{\delta_\pi} \ar[d]_{{}^t\delta'} & T^*X \times T^*_{\Delta_Y}(Y \times Y) \times T^*Z \ar[l] \ar[d] & \\
T^*(X \times Y \times Z) & T^*X \times Y \times T^*Z \ar[l]^{{}^tq'_{13}} \ar[d]^{q_{13\pi}} & \\
& T^*X \times T^*Z. &}
\end{equation}

For $F_1 \in D^b(k_{(X \times Y)_{sa}})$ and $F_2 \in D^b(k_{(Y \times Z)_{sa}})$ set $F_1 \circ F_2=Rq_{13!!}(\imin {q_{12}}F_1 \otimes \imin {q_{23}}F_2)$ and for $G_1 \in D^b(k_{(T^*X \times T^*Y)_{sa}})$ and $G_2 \in D^b(k_{(T^*Y \times T^*Z)_{sa}})$ set $G_1 \stackrel{a}{\circ} G_2=Rp^a_{13 \RP !!}(p_{12}^{a-1}G_1 \otimes p_{23}^{a-1}G_2)$.
We need this proposition which follows from the functorial properties
of $\muh^{sa}$ (it is an adaptation of Proposition 4.4.11 of
\cite{KS90}).

\begin{prop} \label{4.4.11}  Let us consider the sheaves $K_1 \in D^b_{\rc}(k_{X \times
Y}),$ $ F_1 \in D^b(k_{(X \times Y)_{sa}})$, $ K_2 \in
D^b_{\rc}(k_{Y \times Z})$ and $F_2 \in D^b(k_{(Y \times Z)_{sa}})$. Suppose that $q_{13}$ is proper on $\supp(\imin {q_{12}} K_1 \otimes \imin {q_{23}} K_2)$.
There is a morphism
\begin{equation}\label{4.4.14}
\muh^{sa}(K_1,F_1) \stackrel{a}{\circ}
\muh^{sa}(K_2,F_2)
\to \muh^{sa}(K_1 \circ K_2,F_1 \circ F_2).
\end{equation}
\end{prop}

\begin{prop} Let $\lambda=\varnothing,t$. Let $K_1 \in D^b_{\rc}(\CC_{X \times Y})$ and $K_2 \in D^b_{\rc}(\CC_{Y \times Z})$. Suppose that $q_{13}$ is proper on $\supp(\imin {q_{12}} K_1 \otimes \imin {q_{23}} K_2)$. Morphism
\eqref{4.4.14} defines a morphism
\begin{eqnarray}\label{4.4.14ring}
& \muh^{sa}(K_1,\OO_{X \times Y}^{\lambda (0,d_Y)}) \stackrel{a}{\circ}
\muh^{sa}(K_2,\OO_{Y \times Z}^{\lambda (0,d_Z)}) & \\
& \to \muh^{sa}(K_1 \circ K_2,\OO_{X \times Z}^{\lambda (0,d_Z)})[-d_Y]. \notag &
\end{eqnarray}
\end{prop}
\dim\ \ It follows from \eqref{4.4.14} setting $F_1=\OO_{X \times Y}^{\lambda (0,d_Y)}$, $F_2=\OO_{Y \times Z}^{\lambda (0,d_Z)}$ and using the integration morphism $\OO_{X \times Y}^{\lambda (0,d_Y)} \circ \OO_{Y \times Z}^{\lambda (0,d_Z)} \to \OO_{X \times Z}^{\lambda (0,d_Z)}[-d_Y]$. \\
\qed

\begin{cor} Morphism \eqref{4.4.14} induces the ring structures on $\E^\R_X$ and $\E^{\R,f}_X$.
\end{cor}
\dim\ \ It follows applying $\imin \rho$ to \eqref{4.4.14ring} with $X=Y=Z$ and $K_1=K_2=\CC_\Delta[-d_X]$.\\
\qed

\begin{prop} Let $F \in D^b_{\rc}(\CC_X)$. Morphism
\eqref{4.4.14} defines a morphism
\begin{equation}\label{safreccia}
\mu_\Delta^{sa}\OO_{X \times X}^{\lambda(0,d_X)}[d_X] \otimes
\muh^{sa}(F,\ol_X) \to \muh^{sa}(F,\ol_X).
\end{equation}
\end{prop}
\dim\ \ We apply Proposition \ref{4.4.11} with $X=Y$, $Z=\{{\rm point}\}$. We set $(K_1,K_2,F_1,F_2)=(\CC_\Delta[-d_X],F,\OO^{\lambda(0,d_X)}_{X \times X},\ol_X)$.
In this case we have $\CC_\Delta \circ F \simeq F$. We obtain the desired morphism using the integration morphism $\OO^{\lambda (0,d_X)}_{X \times X} \circ \ol_X \to \ol_X[-d_X]$.\\
\qed

Applying the functor $\imin \rho$ to \eqref{safreccia}, we find the
morphisms of \cite{An94} and \cite{KS90} (recall that $\imin \rho\muh^{sa}(F,\ot_X) \simeq t\muh(F,\OO_X)$).

\begin{cor}\label{safreccia1} Morphism \eqref{safreccia} induces  morphisms
\begin{eqnarray}
\E^{\R,f}_X \otimes \imin \rho\muh^{sa}(F,\ot_X) & \to &
\imin \rho\muh^{sa}(F,\ot_X) \label{Androfreccia} \\
\E^{\R}_X \otimes \muh(F,\OO_X) & \to &
\muh(F,\OO_X) \label{Somfreccia}
\end{eqnarray}
which induces a structure of $\E^{\R,f}_X$-module (resp. $\E^\R_X$-module) on the sheaves $H^k\imin \rho
\muh^{sa}(F,\ot_X)$ (resp. $H^k\muh(F,\OO_X)$), for each $k \in \Z$ .
\end{cor}

\

Now we will study the action of $\E^{\R,f}_X$ on formal
microlocalization. We first need to introduce the sheaf of tempered $\C^\infty$-functions.

\begin{df} Let $X$ be a real analytic manifold and let $U \in \op(X)$. Let $f \in \Gamma(U;\C^\infty_X)$. One says that $f$ has poynomial growth at $p \in X$ if for a local coordinate system $(x_1,\ldots,x_n)$ around $p$, there exists a compact neighborhood $K$ of $p$ and $N \in \N$ such that
$$
\sup_{x \in K \cap U}(d(x,K \setminus U))^N|f(x)| < \infty.
$$
One says that $f$ is tempered at $p$ if all its derivatives have polynomial growth at $p$. One says that $f$ is tempered if it is tempered at any point.
\end{df}

\begin{df} One denotes by $\C^{\infty,t}_X$ the presheaf of tempered
$\C^\infty$-functions on $X_{sa}$ defined as follows:
$$U \mapsto \{f \in \Gamma(U;\C^\infty_X),\;f\text{ is tempered}\}.$$
\end{df}
 As a consequence of a result of \cite{KS96}, for $U,V \in \op(X_{sa})$ the sequence
$$\lexs{\C^{\infty,t}_X(U \cup
V)}{\C^{\infty,t}_X(U)\oplus\C^{\infty,t}_X(V)}{\C^{\infty,t}_X(U \cap V)}$$ is exact. Then $\C^{\infty,t}_X$
is a sheaf on $X_{sa}$. Moreover $\mathrm{R}\Gamma(U;\C^{\infty,t}_X)$ is concentrated in degree zero for any $U \in \op(X_{sa})$.

We have the following results (see \cite{KS01}).

\begin{prop} For each $F \in D^b_{\rc}(\CC_X)$ one has the isomorphism
$$
\imin \rho \rh(F,\C^{\infty,t}_X) \simeq \tho(F,\C^\infty_X),
$$
where $\tho(F,\C^\infty_X)$ is the sheaf of \cite{KS96}. When $F=\CC_U$, $U \in \op(X_{sa})$ it is defined by $V \mapsto \C_V^{\infty,t}(U \cap V)$.
\end{prop}

\begin{prop} Let $X$ be a complex manifold, $X_\R$ the underlying real analytic manifold and $\overline{X}$ the conjugate manifold. Then
$$
\ot_X \simeq \rh_{\rho_!\D_{\overline{X}}}(\rho_!\OO_{\overline{X}},\C^{\infty,t}_{X_\R}). $$
\end{prop}

We prove the following result.

\begin{lem} \label{smoothtemp} Let $f:X \to Y$ be a smooth morphism of real analytic manifolds. Then we have the isomorphism
$$
f^{-1}\C^{\infty,t}_Y \iso \rh_{\rho_!\D_X}(\ddxy,\C^{\infty,t}_X).
$$
\end{lem}
\dim\ \ We may reduce to the case of a projection $\pi:Y \times \R \to Y$. We shall prove that the morphism
$$
\partial_t:\C^{\infty,t}_{Y \times \R} \to \C^{\infty,t}_{Y \times \R},
$$
where $t$ denotes the variable in $\R$, is surjective.
Let $U \in \op^c((Y \times \R)_{sa})$, then by Lemma \ref{goodcover} it admits a finite covering $\{U_i\}_{i=1}^N$ such that each $U_i$ is simply connected and the intersections of each $U_i$ with the fibers of $\pi$ are contractible (or empty).
Hence we may reduce to the case that the intersections of $U$ with the fibers of $\pi$ are contractible (or empty). Moreover we can assume that
$$
U = \{(x,t) \in Y \times \R; \; f(x) < t < g(x)\},
$$
where $f,g:\pi(U) \to \R$ are continuous subanalytic maps and $\pi(U)$ is simply connected.
Let us consider $h,k:\pi(U) \to \R$ subanalytic and $\varphi \in \Gamma(\pi(U);\C^\infty_Y)$ such that $f<h<\varphi<k<g$.

Let $s \in \Gamma(U;\C^{\infty,t}_{Y \times \R})$ and define
$$
\widetilde{s}(x,t) = \int_{(x,\varphi(x))}^{(x,t)}s(x,\tau){\rm d}\tau.
$$
Then $ \widetilde{s} \in \Gamma(U;\C^\infty_{Y \times \R})$ and $\partial_t\widetilde{s}=s$. Moreover
$$
|\widetilde{s}(x,t)| \leq |\varphi(x)-t| \sup_{(x,\tau) \in \{x\} \times [\varphi(x),t]}|s(x,\tau)|.
$$
Since $U$ is bounded, there exists $M>0$ such that $|\varphi(x)-t| < M$ for each $(x,t) \in U$. Since $s$ is tempered, for each $x \in \pi(U)$ and each $\tau \in [\varphi(x),t]$ there exist $c_1,r_1>0$ such that
\begin{eqnarray*}
|s(x,\tau)| & \leq & c_1{1 \over d((x,\tau),\partial U)^{r_1}} \\
& \leq & c_1{1 \over \min(d((x,t),\partial U),d((x,h(x)),\partial U),d((x,k(x)),\partial U))^{r_1}}.
\end{eqnarray*}
As a consequence of \L ojaciewicz's inequality (see Theorem 6.4 of \cite{BM88}) there exist $c_2,r_2>0$ such that
$$
d((x,h(x)),\partial U),d((x,k(x)),\partial U) \geq c_2d(x,\partial(\pi(U)))^{r_2} \geq c_2d((x,t),\partial U)^{r_2}.
$$
Hence there exists $c,r>0$ such that
$$
\widetilde{s}(x,t) \leq c{1 \over d((x,t),\partial U)^r}
$$
and the result follows.\\
\qed

\begin{lem}\label{prodtrans} Let $f:X \to Y$ be a smooth morphism of real analytic manifolds. Let $\M,\NN \in \D^b(\D_X)$. There is a natural morphism
$$
\rh_{\D_X}(\ddxy,\M) \overset{L}{\underset{\imin f\A_Y}{\otimes}} \rh_{\D_X}(\ddxy,\NN) \to \rh_{\D_X}(\ddxy,\M \overset{L}{\underset{\A_X}{\otimes}} \NN).
$$
\end{lem}
\dim\ \ 
By Lemma 4.9 of \cite{Ka03} we have
\begin{equation}\label{4.9}
\ddyx \overset{L}{\underset{\A_X}{\otimes}} \ddxy \simeq \ddyx \overset{L}{\underset{\imin f\A_Y}{\otimes}} \imin f\D_Y.
\end{equation}
Then if $\M$ is a $\D_X$-module
\begin{eqnarray*}
(\ddyx \overset{L}{\underset{\D_X}{\otimes}} \M) \overset{L}{\underset{\imin f\A_Y}{\otimes}} \imin f\D_Y
& \simeq & (\ddyx  \overset{L}{\underset{\imin f\A_Y}{\otimes}} \imin f \D_Y) \overset{L}{\underset{\D_X}{\otimes}} \M\\
& \simeq & (\ddyx  \overset{L}{\underset{\A_X}{\otimes}} \ddxy) \overset{L}{\underset{\D_X}{\otimes}} \M\\
& \simeq & \ddyx   \overset{L}{\underset{\D_X}{\otimes}} (\M\overset{L}{\underset{\A_X}{\otimes}} \ddxy).
\end{eqnarray*}
Now when $f$ is smooth $\ddyx \overset{L}{\underset{\D_X}{\otimes}} \cdot \simeq \rh_{\D_X}(\ddxy,\cdot)[d_X-d_Y]$. Then if $\NN$ is another $\D_X$-module
\begin{eqnarray*}
\lefteqn{\rh_{\D_X}(\ddxy,\M) \overset{L}{\underset{\imin f\A_Y}{\otimes}} \rh_{\D_X}(\ddxy,\NN)} \\
& \simeq & \rh_{\D_X}(\ddxy,\M \overset{L}{\underset{\A_X}{\otimes}} \ddxy \overset{L}{\underset{\imin f\D_Y}{\otimes}} \rh_{\D_X}(\ddxy,\NN)) \\
& \to & \rh_{\D_X}(\ddxy,\M \overset{L}{\underset{\A_X}{\otimes}} \NN).
\end{eqnarray*}
\qed

\begin{lem}\label{1pezzo} Let $X$ be a real analytic manifold. Let $F,G \in D^b_{\rc}(\CC_X)$ and let $S$ be a closed subanalytic subset of $X$. There is a morphism
\begin{eqnarray*}
& \imin \rho \rh(F,(\C^{\infty,t}_X)_S) \otimes_{\A_X} \imin \rho \rh(D'((F \otimes G)_S), \CW_X) & \\
& \to \imin \rho \rh(D'(G_S),\CW_X). &
\end{eqnarray*}
\end{lem}
\dim\ \ (i) Let $V_1,V_2 \in \op(X_{sa})$. The sheaf $\imin \rho \Gamma_{V_1}(\C^{\infty,t}_X)_S$ is concentrated in degree zero since $\C^{\infty,t}_X$ is $\Gamma(U;\cdot)$-acyclic for each $U \in \op(X_{sa})$. Moreover the sheaves $\imin \rho \rh(D'\CC_{V_1 \cap V_2 \cap S},\C^{\infty,{\rm w}}_X) \simeq \CC_{V_1 \cap V_2 \cap S} \wtens \C^\infty_X$ and $\imin \rho \rh(D'\CC_{V_2 \cap S},\C^{\infty,{\rm w}}_X) \simeq \CC_{V_2 \cap S} \wtens \C^\infty_X$ are also concentrated in degree zero. There is a morphism
    \begin{equation}\label{V1V2}
    \imin \rho \Gamma_{V_1}(\C^{\infty,t}_X)_S \otimes_{\A_X} \CC_{V_1 \cap V_2 \cap S} \wtens \C^\infty_X \to \CC_{V_2 \cap S} \wtens \C^\infty_X.
\end{equation}
This follows since the multiplication of a function tempered on $V_1$ by a function vanishing with all its derivatives outside $V_1$ is a function vanishing with all its derivatives outside $V_1$.

(ii) By Theorem 1.1 of \cite{KS96} the morphism \eqref{V1V2} extends to a morphism
\begin{equation}\label{V1G}
    \imin \rho \Gamma_{V_1}(\C^{\infty,t}_X)_S \otimes_{\A_X} G_{V_1 \cap S} \wtens \C^\infty_X \to G_S \wtens \C^\infty_X,
\end{equation}
 functorial in $G \in \mod_{\rc}(\CC_X)$. By adjuction this gives a morphism
\begin{equation}\label{V1Gad}
    \imin \rho \Gamma_{V_1}(\C^{\infty,t}_X)_S \to \ho_{\A_X}(G_{V_1 \cap S} \wtens \C^\infty_X,G_S \wtens \C^\infty_X).
\end{equation}
By Theorem 1.1 of \cite{KS96} the morphism \eqref{V1Gad} extends to a morphism
\begin{equation}\label{FGad}
    \imin \rho \ho(F,(\C^{\infty,t}_X)_S) \to \ho_{\A_X}((F \otimes G)_S \wtens \C^\infty_X,G_S \wtens \C^\infty_X).
\end{equation}
functorial in $F \in \mod_{\rc}(\CC_X)$.

(iii) Let $F,G \in D^b_{\rc}(\CC_X)$. We have the following chain of morphisms
\begin{eqnarray*}
\imin \rho \rh(F,(\C^{\infty,t}_X)_S) & \iso & R(\imin \rho \ho(F,(\C^{\infty,t}_X)_S))\\
& \to & R(\ho_{\A_X}((F \otimes G)_S \wtens \C^\infty_X,G_S \wtens \C^\infty_X)) \\
& \to & \rh_{\A_X}((F \otimes G)_S \wtens \C^\infty_X,G_S \wtens \C^\infty_X),
\end{eqnarray*}
where the first isomorphism follows since $\imin \rho$ is exact and $(\cdot)_S$ sends quasi-injective objects to quasi-injective objects, the second arrow follows from \eqref{FGad} and the third one is a canonical morphism of derived functors (see \cite{KS}, Proposition 13.3.13).

By adjunction we obtain the desired morphism.\\ \qed

\begin{lem} \label{2pezzo} Let us consider the normal deformation of the diagonal in $X
\times X$ of diagram \eqref{normdef}. Let $F,G \in
D^b_{\rc}(\CC_X)$. There is a morphism
\begin{eqnarray*}
& \imin \rho \nu^{sa}_\Delta\rh(\imin {q_1}F,\imin {q_2}\C^{\infty,t}_X) \otimes_{\A_X}  \imin \rho \nu^{sa}_\Delta\rh(\imin {q_1}D'(F \otimes G),\imin {q_2}\CW_X) & \\
& \to  \imin \rho \nu^{sa}_\Delta\rh(\imin {q_1}D'G,\imin {q_2}\CW_X). &
\end{eqnarray*}
\end{lem}
\dim\ \ (i) As in the proof of Theorem \ref{muhomtwr}, if $X$ is a real analytic manifold, $K \in D^b_{\rc}(\CC_X)$, $\lambda=t,{\rm w}$, we have
\begin{eqnarray*}
\lefteqn{\imin \rho \nu^{sa}_\Delta\rh(\imin {q_1}K,\imin {q_2}\C^{\infty,\lambda}_X)} \\
& \simeq & \imin \rho \imin s \rh((\imin {p_1}K)_\Omega,\imin {p_2}\C^{\infty,\lambda}_X) \\
& \simeq & \imin \rho \imin s \rh((\imin {p_1}K)_\Omega,\rh_{\rho_!\D_{\widetilde{X \times X}}}(\rho_!\D_{\widetilde{X \times X} \to X},\C^{\infty,\lambda}_{\widetilde{X \times X}})) \\
& \simeq & \imin s \rh_{\D_{\widetilde{X \times X}}}(\D_{\widetilde{X \times X} \to X},\imin \rho\rh((\imin {p_1}K)_\Omega,\C^{\infty,\lambda}_{\widetilde{X \times X}})),
\end{eqnarray*}
where the second isomorphism follows from Lemma \ref{smoothtemp}.

(ii)  By Lemma \ref{omegabar} for $H \in D^b_{\rc}(\CC_X)$ we have $(\imin {p_1}D'H)_\Omega \simeq D'((\imin {p_1}H)_{\overline{\Omega}})$ and $\mathrm{R}\Gamma_\Omega\imin {p_2} \C^{\infty,t}_X \simeq (\imin {p_2}\C^{\infty,t}_X)_{\overline{\Omega}}$.

(iii) By Lemma \ref{prodtrans} with $(X,Y)=(\widetilde{X \times X},X)$, $\M=\rh(\imin {p_1}F,(\C^{\infty,t}_{\widetilde{X \times X}})_{\overline{\Omega}})$, $\NN=\rh(D'((\imin {p_1}(F \otimes G))_{\overline{\Omega}}), \CW_{\widetilde{X \times X}}))$,
we are reduced to find a morphism
\begin{eqnarray*}
& \imin \rho \rh(\imin {p_1}F,(\C^{\infty,t}_{\widetilde{X \times X}})_{\overline{\Omega}}) \otimes_{\A_{\widetilde{X \times X}}} \imin \rho \rh(D'((\imin {p_1}(F \otimes G))_{\overline{\Omega}}), \CW_{\widetilde{X \times X}}) & \\
& \to \imin \rho \rh(D'((\imin {p_1}G)_{\overline{\Omega}}),\CW_{\widetilde{X \times X}}) &
\end{eqnarray*}
which follows replacing
$(X,S,F,G)$ with $(\widetilde{X \times X},\overline{\Omega},\imin {p_1}F,\imin {p_1}G)$ in Lemma \ref{1pezzo}.
\qed

Let us consider the complex case. Let $X$ be a complex manifold.

\begin{lem}\label{prodcplx} Let ${\cal L},\H \in D^b(\D_{X_\R})$. There is a natural morphism
$$
\rh_{\D_{\overline{X}}}(\OO_{\overline{X}},{\cal L}) \overset{L}{\underset{\OO_X}{\otimes}} \rh_{\D_{\overline{X}}}(\OO_{\overline{X}},\H) \to \rh_{\D_{\overline{X}}}(\OO_{\overline{X}},{\cal L} \overset{L}{\underset{\A_{X_\R}}{\otimes}} \H).
$$
\end{lem}
\dim\ \ By definition we have ${\cal L} \overset{L}{\underset{\A_{X_\R}}{\otimes}} \H = \D_{X_\R \to X_\R \times X_\R} \overset{L}{\underset{\A_{X_\R \times X_\R}}{\otimes}} ({\cal L} \overset{\D}{\boxtimes} \M)$. Hence we get
$$
{\cal L} \overset{\D}{\boxtimes} \H \to \rh_{\D_{X_\R}}(\D_{X_\R \to X_\R \times X_\R}, {\cal L} \overset{L}{\underset{\A_{X_\R}}{\otimes}} \H).
$$
There is a chain of morphisms
\begin{eqnarray*}
\lefteqn{\rh_{\D_{\overline{X}}}(\OO_{\overline{X}},{\cal L}) \overset{L}{\underset{\OO_X}{\otimes}} \rh_{\D_{\overline{X}}}(\OO_{\overline{X}},{\cal H})} \\
& \simeq & \D_{X \to X \times X} \overset{L}{\underset{\D_{X \times X}}{\otimes}} (\rh_{\D_{\overline{X}}}(\OO_{\overline{X}},{\cal L}) \overset{\D}{\boxtimes} \rh_{\D_{\overline{X}}}(\OO_{\overline{X}},{\cal H})) \\
& \to & \D_{X \to X \times X} \overset{L}{\underset{\D_{X \times X}}{\otimes}} \rh_{\D_{\overline{X}} \boxtimes \D_{\overline{X}}}(\OO_{\overline{X}} \boxtimes \OO_{\overline{X}},\D_{X \times X} \overset{L}{\underset{\D_X \boxtimes \D_X}{\otimes}} ({\cal L} \boxtimes {\cal H})) \\
& \to & \D_{X \to X \times X} \overset{L}{\underset{\D_{X \times X}}{\otimes}} \rh_{\D_{\overline{X} \times \overline{X}}}(\OO_{\overline{X} \times \overline{X}},\D_{\overline{X} \times \overline{X}} \overset{L}{\underset{\D_{\overline{X}} \boxtimes \D_{\overline{X}}}{\otimes}} \D_{X \times X} \overset{L}{\underset{\D_X \boxtimes \D_X}{\otimes}} ({\cal L} \boxtimes {\cal H})) \\
& \to & \D_{X \to X \times X} \overset{L}{\underset{\D_{X \times X}}{\otimes}} \rh_{\D_{\overline{X} \times \overline{X}}}(\OO_{\overline{X} \times \overline{X}}, {\cal L} \overset{\D}{\boxtimes} \H) \\
& \to & \D_{X \to X \times X} \overset{L}{\underset{\D_{X \times X}}{\otimes}} \rh_{\D_{\overline{X} \times \overline{X}}}(\OO_{\overline{X} \times \overline{X}},\rh_{\D_{X_\R}}(\D_{X_\R \to X_\R \times X_\R}, {\cal L} \overset{L}{\underset{\A_{X_\R}}{\otimes}} \H)) \\
& \simeq & \D_{X \to X \times X} \overset{L}{\underset{\D_{X \times X}}{\otimes}} \rh_{\D_X}(\D_{X \to X \times X}, \rh_{\D_{\overline{X}}}(\OO_{\overline{X}},{\cal L} \overset{L}{\underset{\A_{X_\R}}{\otimes}} \H)) \\
& \to & \rh_{\D_{\overline{X}}}(\OO_{\overline{X}},{\cal L} \overset{L}{\underset{\A_{X_\R}}{\otimes}} \H).
\end{eqnarray*}
\qed

\begin{lem} \label{morfnu} Let us consider the normal deformation of the diagonal in $X
\times X$ of diagram \eqref{normdef}. Let $F,G \in
D^b_{\rc}(\CC_X)$. There is a morphism
\begin{eqnarray*}
& \imin \rho \nu^{sa}_\Delta\rh(\imin {q_1}F,\imin {q_2}\ot_X) \otimes_{\OO_X}  \imin \rho \nu^{sa}_\Delta\rh(\imin {q_1}D'(F \otimes G),\imin {q_2}\OW_X) & \\
& \to  \imin \rho \nu^{sa}_\Delta\rh(\imin {q_1}D'G,\imin {q_2}\OW_X). &
\end{eqnarray*}
\end{lem}
\dim\ \ If $X$ is a complex manifold, $K \in D^b_{\rc}(\CC_X)$, $\lambda=t,{\rm w}$, we have
\begin{eqnarray*}
\lefteqn{\imin \rho \nu^{sa}_\Delta\rh(\imin {q_1}K,\imin {q_2}\OO^\lambda_X)} \\
& \simeq & \imin \rho \imin s \rh((\imin {p_1}K)_\Omega,\imin {p_2}\OO^\lambda_X) \\
& \simeq & \imin \rho \imin s \rh((\imin {p_1}K)_\Omega,\imin {p_2}\rh_{\rho_!\D_{\overline{X}}}(\rho_!\OO_{\overline{X}},\C^{\infty,\lambda}_{X_\R})) \\
& \simeq & \rh_{\D_{\overline{X}}}(\OO_{\overline{X}},\imin \rho \imin s \rh((\imin {p_1}K)_\Omega,\imin {p_2}\C^{\infty,\lambda}_{X_\R})) \\
& \simeq & \rh_{\D_{\overline{X}}}(\OO_{\overline{X}},\imin \rho \nu^{sa}_\Delta\rh(\imin {q_1}K,\imin {q_2}\C^{\infty,\lambda}_{X_\R})).
\end{eqnarray*}
Set
\begin{eqnarray*}
{\cal L} & = & \imin \rho \nu^{sa}_\Delta\rh(\imin {q_1}F,\imin {q_2}\C^{\infty,t}_{X_\R})), \\
\H & = & \imin \rho \nu^{sa}_\Delta\rh(\imin {q_1}D'(F \otimes G),\imin {q_2}\C^{\infty,\mathrm{w}}_{X_\R})).
\end{eqnarray*}
By Lemma \ref{prodcplx} there is a natural morphism
$$
\rh_{\D_{\overline{X}}}(\OO_{\overline{X}},{\cal L}) \overset{L}{\underset{\OO_X}{\otimes}} \rh_{\D_{\overline{X}}}(\OO_{\overline{X}},\H) \to \rh_{\D_{\overline{X}}}(\OO_{\overline{X}},{\cal L} \overset{L}{\underset{\A_{X_\R}}{\otimes}} \H).
$$
Then the result follows from Lemma \ref{2pezzo}. \\
\qed


\begin{lem}\label{intmorph} Let $f:X \to Y$ be a smooth morphism of complex manifolds. Then there is a natural morphism
$$
Rf_{!!}\Omega^{\rm w}_X[d_X] \to \Omega_Y[d_Y].
$$
\end{lem}
\dim\ \ By Theorem \ref{isoOW} we have the isomorphism
$$
f^!\OWY[2d_Y]  \iso \rh_{\rho_!\D_X}(\rho_!\ddxy,\OWX)[2d_X].
$$
We have $\rh_{\rho_!\D_X}(\rho_!\ddxy,\OWX) \simeq \rho_!\ddyx \overset{L}{\underset{\rho_!\D_X}{\otimes}} \OWX[d_Y-d_X]$.
Hence we get $$\rho_!\ddyx \overset{L}{\underset{\rho_!\D_X}{\otimes}} \OWX \simeq f^!\OWY[d_Y-d_X].$$ By adjunction we get $$Rf_{!!}(\rho_!\ddyx \overset{L}{\underset{\rho_!\D_X}{\otimes}} \OWX) \to \OWY[d_Y-d_X].$$
From this we can deduce $$Rf_{!!}\Omega_X^{\rm w} \to Rf_{!!}(\Omega_X^{\rm w} \overset{L}{\underset{\rho_!\D_X}{\otimes}} \rho_!\ddxy) \to \Omega_Y^{\rm w}[d_Y-d_X].$$
\qed

Let us consider the diagram \eqref{michint} with $Z=\{{\rm point}\}$. Set $p_X:T^*X \times T^*Y \to T^*X$, $p_Y:T^*X \times T^*Y \to T^*Y$, $q_X:X \times Y \to X$, $q_Y:X \times Y \to Y$.

\begin{prop}\label{mu10.6} Let $G \in D^b_{\rc}(\CC_X)$ and $K \in D^b_{\rc}(\CC_{X \times Y})$ such that 
$q_Y$ is proper on $\supp(\imin {q_X} G) \cap \supp (K)$. 
Then we have a morphism
\begin{eqnarray}\label{multtw}
& \imin \rho\muh^{sa}(K,\OO^{t(0,d_Y)}_{X \times Y})[d_Y] \stackrel{a}{\circ} \imin \rho \muh^{sa}(D'(K \circ G),\OW_Y) &  \\
& \to \imin \rho\muh^{sa}(D'G,\OW_X). \notag &
\end{eqnarray}
\end{prop}
\dim\ \ We will prove the assertion in several steps. Set
\begin{eqnarray*}
H_1 & = & \imin \rho \nu_\Delta^{sa}\rh(\imin {q_1}K,q_2^!\OO^{t(0,d_Y)}_{X \times Y}) \\
& \simeq & \imin \rho \nu_\Delta^{sa}\rh(\imin {q_1}K,\imin {q_2}\OO^{t(0,d_Y)}_{X \times Y})[2d_{X \times Y}] \\
H_2 & = & \imin \rho \nu_\Delta^{sa}\rh(\imin {q_1}D'(K \circ G),q_2^!\OW_Y) \\
& \simeq & \imin \rho\nu_\Delta^{sa}\rh(\imin {q_1}D'(K \circ G),\imin {q_2}\OW_Y)[2d_Y].
\end{eqnarray*}
Since the Fourier-Sato transform commutes with $\imin \rho$ we have
$$
H_1^\land \simeq \imin \rho \muh^{sa}(K,\OO^{t(0,d_Y)}_{X \times Y}), \ \ H_2^\land \simeq \imin \rho \muh^{sa}(D'(K \circ G),\OW_Y).
$$

(i) By the commutativity of the diagram \eqref{michint} we have an isomorphism
$$
Rp_{X !}^a((H_1^{\wedge })^a \otimes p_Y^{a-1}H_2^\wedge) \simeq Rq_{X\pi !}{}^tq_X'{}^{-1}R{}^t\delta'_{!}\imin{\delta_\pi}(H_1^\wedge \boxtimes H_2^\wedge).
$$

(ii) By Proposition 3.7.15 of \cite{KS90} we have an isomorphism
$$
(H_1)^\wedge \boxtimes (H_2)^\wedge \iso (H_1 \boxtimes H_2)^\wedge,
$$

(iii) Denote by $Tq_Y:T(X \times Y) \to TY$ the tangent map. By Propositions 3.7.13 and 3.7.14 of \cite{KS90} we have the isomorphism
$$
R{}^t\delta'_{!}\imin{\delta_\pi}(H_1 \boxtimes H_2)^\wedge \simeq (H_1 \otimes \imin {Tq_Y}H_2)^\wedge[-2d_Y].
$$

(iv) We have the chain of morphisms
\begin{eqnarray*}
\lefteqn{\imin {Tq_Y} \nu^{sa}_\Delta\rh(\imin {q_1}D'(K\circ G),\imin {q_2}\OW_Y)}\\
& \simeq & \nu^{sa}_\Delta\rh(\imin {q_1}\imin {q_Y}D'(K \circ G),\imin {q_2}\imin {q_Y}\OW_Y) \\
& \simeq & \nu^{sa}_\Delta\rh(\imin {q_1}D'(\imin {q_Y}q_{Y_*}(K \otimes \imin {q_X}G)),\imin {q_2}\imin {q_Y}\OW_Y) \\
& \to & \nu^{sa}_\Delta\rh(\imin {q_1}D'(K \otimes \imin {q_X}G),\imin {q_2}\imin {q_Y}\OW_Y) \\
& \simeq & \nu^{sa}_\Delta\rh(\imin {q_1}D'(K \otimes \imin {q_X}G),\imin {q_2}\rh_{\rho_!\D_{X \times Y}}(\rho_!\dd{X \times Y}{}{Y},\OW_{X \times Y})),
\end{eqnarray*}
where the first isomorphism follows since $q_Y$ is smooth, the second one since $\supp(\imin {q_X} G) \cap \supp (K)$ is
proper over $Y$ and the last one is given by inverse image formula for Whitney holomorphic functions.

(v) We have a morphism
$$
(H_1 \otimes \imin {Tq_Y}H_2)^\wedge[-2d_Y] \to \imin \rho \muh^{sa}(\imin {q_X}D'G,\OO^{\mathrm{w}(0,d_Y)}_{X \times Y}).
$$
To prove the existence of this morphism we use (iv) to prove the morphism
$$
H_1 \otimes \imin {Tq_Y}H_2[-2d_Y] \to \nu^{sa}_\Delta\rh(\imin {q_1}\imin {q_X}D'G,\imin {q_2}\OO^{\mathrm{w}(0,d_Y)}_{X \times Y})[2d_{X \times Y}].
$$
Hence we may reduce to the case of the morphism
\begin{eqnarray*}
& \imin \rho \nu^{sa}_\Delta\rh(\imin {q_1}K,\imin {q_2}\ot_{X \times Y}) \otimes  \imin \rho \nu^{sa}_\Delta\rh(\imin {q_1}D'(K \otimes \imin {q_X}G),\imin {q_2}\OW_{X \times Y}) & \\
& \to  \imin \rho \nu^{sa}_\Delta\rh(\imin {q_1}D'(\imin {q_X}G),\imin {q_2}\OW_{X \times Y}). &
\end{eqnarray*}
This is a consequence of Lemma \ref{morfnu} with
$(X,F,G)$ replaced by $(X \times Y,K,\imin {q_X}G)$.

(vi) We have the chain of morphisms
\begin{eqnarray*}
\lefteqn{Rq_{X\pi !}{}^tq'_X{}^{-1} \imin \rho\muh^{sa}(\imin {q_X}D'G,\OO^{\mathrm{w}(0,d_Y)}_{X \times Y})}\\
& \to & \imin \rho\muh^{sa}(Rq_{X*}\imin {q_X}D'G,Rq_{X!!}\OO^{\mathrm{w}(0,d_Y)}_{X \times Y}) \\
& \to & \imin \rho\muh^{sa}(D'G,\OW_X)[-d_Y],
\end{eqnarray*}
where the second morphism is a consequence of the integration morphism $Rq_{X!!}\OO^{\mathrm{w}(0,d_Y)}_{X \times Y} \to \OW_X[-d_Y]$ defined in Lemma \ref{intmorph} (see also Remark 3.4 of \cite{KS96}) and the fact that $Rq_{X*}\imin {q_X} \simeq \id$.
Composing morphisms (i)-(vi) we get the desired morphism. \\
\qed

\begin{cor}\label{safreccia2} Let $F \in D^b_{\rc}(\CC_X)$. Morphism \eqref{multtw} defines
a morphism
\begin{equation}\label{Colinfreccia}
\E^{\R,f}_X \otimes \imin \rho\muh^{sa}(F,\OW_X) \to
\imin \rho\muh^{sa}(F,\OW_X)
\end{equation}
which induces a structure of $\E_X^{\R,f}$-module on $H^k\imin \rho \muh^{sa}(F,\OW_X)$ for each $k \in \Z$.
\end{cor}
\dim\ \ We apply Proposition \ref{mu10.6} setting $X=Y$ and $(G,K)=(D'F,\CC_\Delta)$.
In this case we have $D'(\CC_\Delta
\circ D'F) \simeq D'D'F \simeq F$.\\
\qed

In this way we find the
morphism of \cite{Co98} (recall that $\imin \rho \muh^{sa}(F,\OW_X) \simeq (D'F \underset{\mu}{\wtens} \OO_X)^a$)
$$
(\E^{\R,f}_X)^a \otimes F \underset{\mu} {\wtens} \OO_X \to F \underset{\mu}{\wtens} \OO_X.
$$

\begin{oss} We would like to see the compatibility between this morphism and the one of Andronikof (\cite{An94}, Proposition 3.3.10). Steps (i) to (iii) of Proposition \ref{mu10.6} are the same. We need the compatibility between the multiplications. We will see the compatibility between
$$
\imin \rho \mathrm{R}\Gamma_Z\ot_X \otimes \imin \rho \rh(F,\ot_X) \to \imin \rho \rh(F,\ot_X)
$$
and
$$
\imin \rho \mathrm{R}\Gamma_Z\ot_X \otimes \imin \rho \rh(F,\OW_X) \to \imin \rho\rh(F,\OW_X)
$$
when $Z \subset X$ is closed subanalytic and $F \in D^b_{\rc}(\CC_X)$.

We reduce to the case of a real analytic manifold and we use the fact that $\imin\rho\rh(G,\C^{\infty,t}_X) \simeq \tho(G,\C^{\infty}_X)$ and $\imin\rho\rh(G,\CW_X) \simeq D'G \wtens \C^\infty_X$ for $G \in D^b_{\rc}(\CC_X)$. Define $F \wtens \tho(G,\C^\infty_X)=\tho(G,F \wtens \C^\infty_X)$ saying that, if $U,V$ are open subanalytic $\CC_U \wtens \tho(\CC_V,\C^\infty_X)=\tho(\CC_V,\CC_U \wtens \C^\infty_X)$ are $\C^\infty$-functions tempered on $V$ and vanishing up to infinity outside $U$. Then we have
\begin{eqnarray*}
\tho(\CC_Z,\C^\infty_X) \otimes \tho(F,\C^\infty_X) & \to & \tho(\CC_Z,\C^\infty_X) \otimes \tho(F,\CC_Z \wtens \C^\infty_X) \\
& \to & \tho(F_Z,\C^\infty_X) \\
& \to & \tho(F,\C^\infty_X)
\end{eqnarray*}
and
\begin{eqnarray*}
\tho(\CC_Z,\C^\infty_X) \otimes D'F \wtens \C^\infty_X & \to & \tho(\CC_Z,\C^\infty_X) \otimes (D'F)_Z \wtens \C^\infty_X \\
& \to & \tho(\CC_Z,D'F \wtens \C^\infty_X) \\
& \to & D'F \wtens \C^\infty_X.
\end{eqnarray*}
The first and the third arrows of the two diagrams are clearly compatible. Let us see the compatibility between the second arrows. Note that $F \in D^b_{\rc}(\CC_X)$ plays no role in these arrows (it denotes a growth conditions which is preserved after the multiplication), so in order to better understand how they are constructed we set $F=\CC_X$. Let $U=X\setminus Z$. Then $\tho(\CC_Z,\C^\infty_X)$ and $\CC_Z \wtens \C^\infty_X$ are represented by the complexes
$$
\begin{array}{rcccl}
0 & \to & \C^\infty_X & \to & \tho(\CC_U,\C^\infty_X) \\
\CC_U \wtens \C^\infty_X & \to & \C^\infty_X & \to & 0
\end{array}
$$
where in both cases $\C^\infty_X$ is the degree zero of the complex.
The morphism is induced by the following diagram, where the vertical arrows are given by multiplication
$$
\xymatrix{
\CC_U \wtens \C^\infty_X \otimes \C^\infty_X \ar[r] \ar[d] & \CC_U \wtens \C^\infty_X \otimes \tho(\CC_U,\C^\infty_X) \oplus \C^\infty_X \otimes \C^\infty_X \ar[r] \ar[d] & \C^\infty_X \otimes \tho(\CC_U,\C^\infty_X) \ar[d] \\
\CC_U \wtens \C^\infty_X \ar[r] & \CC_U \wtens \C^\infty_X \oplus \C^\infty_X \ar[r] & \tho(\CC_U,\C^\infty_X).
}
$$
In the complex in the second line the first arrow is given by $s \mapsto (s,s)$ and the second one by $(u,v) \mapsto u-v$. Computing the cohomology, it is quasi-isomorphic to $\tho(\CC_Z,\C^\infty_X)$.


\end{oss}

\subsection{Microlocal integral transformations}

In the case of contact transformation the hypothesis of properness of the previous section are not satisfied. Hence we are going to define microlocal operations on $\muh^{sa}(\cdot,\ol_X)$ extending those of \cite{KS90} and \cite{An94}.\\

\noindent Let $\Omega \subset T^*X$. Denote by $D^b(X_{sa},\Omega)$ (resp. $D^b(X,\Omega)$, resp. $D^b_{\rc}(X,\Omega)$ the category $D^b(\CC_{X_{sa}})/{\cal N}_\Omega$ (resp. $D^b(\CC_X)/{\cal N}_\Omega$, resp. $D^b_{\rc}(\CC_X)/{\cal N}_\Omega$), where ${\cal N}_\Omega=\{F \in D^b(\CC_{X_{sa}})\;;\; SS(F) \cap \Omega=\varnothing\}$ (resp. $F \in D^b(\CC_X)$, resp. $F \in D^b_{\rc}(\CC_X)$. It follows from Corollary \ref{muhSS(F)G} that the functor
$$
\imin \rho \muh^{sa}: D^b(X_{sa},\Omega)^{op} \times D^b(X_{sa},\Omega) \to D^b(\Omega)
$$
is well defined.

\begin{nt} If there is no risk of confusion we will write for short $\muh(\cdot,\ol_X)$ instead of $\imin \rho \muh^{sa}(\cdot,\ol_X)$.
\end{nt}

Denote by $\mucirc$ the microlocal composition of kernels of \cite{KS90} (and \cite{An94} for $\R$-constructible sheaves). As usual, given $K \in D^b(\CC_{(X \times Y)_{sa}})$ and $F \in D^b(\CC_{Y_{sa}})$ we set $\Phi^\mu_KF=K \mucirc F$.

\begin{prop} \label{microintransf} (i) Let $X,Y$ be two complex analytic manifolds, let $K \in D^b_{\rc}(\CC_{X \times Y})$, $p_X \in T^*X$, $p_Y \in T^*Y$ such that $SS(K) \cap (\{p_X\} \times T^*Y) \subseteq (p_X,p_Y^a)$ in a neighborhood of this point. Then for each $F\in D^b_{\rc}(\CC_Y)$ and $G\in D^b_{\rc}(\CC_X)$ there are morphisms
\begin{eqnarray}
& \muh(K,\OO^{t(0,d_Y)}_{X \times Y})_{(p_X,p_Y^a)}[d_Y]
\otimes
\muh(F,\ot_Y)_{p_Y} & \label{microt} \\
& \to \muh(\Phi^\mu_KF, \ot_X)_{p_X}. & \notag
\end{eqnarray}
\begin{eqnarray}
& \muh(K,\OO^{t(0,d_Y)}_{X \times Y})_{(p_X,p_Y^a)}[d_Y] \otimes \muh(D'(\Phi^\mu_KG),\OWY)_{p_Y} &  \label{microw} \\
& \to \muh(D'G,\OWX)_{p_X}. \notag &
\end{eqnarray}

(ii) Let $Z$ be another complex analytic manifold, let $K_1 \in D^b_{\rc}(\CC_{X \times Y})$ and $K_2 \in D^b_{\rc}(\CC_{Y \times Z})$ be microlocally composable at $(p_X,p_Y,p_Z) \in T^*X \times T^*Y \times T^*Z$, i.e.
$$
(SS(K_1) \times_{T^*Y} SS(K_2)) \cap p^{a-1}_{13}(p_X,p_Z^a) \subseteq \{((p_X,p_Y^a),(p_Y,p_Z^a))\}
$$
in a neighborhood of $((p_X,p_Y^a),(p_Y,p_Z^a))$. Then there is a morphism
\begin{eqnarray*}
& \muh(K_1,\OO_{X \times Y}^{t (0,d_Y)})_{(p_X,p_Y^a)} \otimes
\muh(K_2,\OO_{Y \times Z}^{t (0,d_Z)}))_{(p_Y,p_Z^a)} & \\
& \to \muh(K_1 \mucirc K_2,\OO_{X \times Z}^{t (0,d_Z)})_{(p_X,p_Z^a)}[-d_Y]& \notag
\end{eqnarray*}
\end{prop}
\dim\ \ The result follows thanks to the morphisms defined in the previous section and adapting the proof of Proposition 3.3.12 of \cite{An94}.\\
\qed

\subsection{Contact transformations}

Let $X,Y$ be two complex analytic manifolds of the same complex dimension $n$ and let $\Omega_X \subset T^*X$, $\Omega_Y \subset T^*Y$ be two open subanalytic subsets. Let $\chi$ be a contact transformation from $\Omega_X$ to $\Omega_Y$. We set $\Lambda \subset \Omega_X \times \Omega_Y^a$ be the Lagrangian manifold associated to the graph of $\chi$ (i.e. $(p_X,p_Y^a) \in \Lambda$ if $p_Y=\chi(p_X)$). We denote by $p_1$ and $p_2^a$ the projections from $\Lambda$ to $\Omega_X$ and $\Omega_Y$ respectively.\\

Let $(p_X,p_Y) \in \Omega_X \times \Omega_Y$ and consider $K \in D^b_{\cc}(X \times Y,(p_X,p_Y^a))$ satisfying the following properties:
\begin{equation}\label{microhyp}
\begin{cases}
    SS(K) \subset \Lambda, \\
    \text{$K$ is simple with shift 0 along $\Lambda$.}
  \end{cases}
\end{equation}
In this situation we have the following results of \cite{KS03} and \cite{An94}.

\begin{prop} Let $K \in D^b_{\cc}(X \times Y,(p_X,p_Y^a))$ satisfying \eqref{microhyp}. Set $K^*=r_*\rh(K,\omega_{X \times Y|Y})$, where $r:X \times Y \to Y \times X$ is the canonical map. Then the functor $\Phi^\mu_K:D^b(X_{sa},p_X) \to D^b(Y_{sa},p_Y)$ and the functor $\Phi^\mu_{K^*}:D^b(Y_{sa},p_Y) \to D^b(X_{sa},p_X)$ are equivalence of categories inverse to each other.
\end{prop}

\begin{lem}\label{An5.2.3} Let $K \in D^b_{\cc}(X \times Y,(p_X,p_Y^a))$ satisfying \eqref{microhyp}. Then $\muh(K,\ot_{X \times Y})$ is concentrated in degree zero.
\end{lem}

\begin{prop} \label{An5.2.1} Let $K \in D^b_{\cc}(X \times Y,(p_X,p_Y^a))$ satisfying \eqref{microhyp} and let $s \in \muh(K,\OO_{X \times Y}^{t(0,n)})_{(p_X,p_Y^a)}$.
\begin{itemize}
\item[(i)] For each $F \in D^b_{\rc}(Y,p_Y)$ and $G \in D^b_{\rc}(X,p_X)$ there are morphisms induced by $s$
\begin{eqnarray*}
\varphi_s:\muh(F,\ot_Y)_{p_Y}[n] & \to & \muh(\Phi^\mu_KF,\ot_X)_{p_X}  \\
\psi_s:\muh(D'(\Phi^\mu_KG),\OWY)_{p_Y}[n] & \to & \muh(D'F,\OWX)_{p_X}.
\end{eqnarray*}
\item[(ii)] Let $Z$ be a $n$-dimensional complex analytic manifold, $\Omega_Z \subset T^*Z$ and let $\chi':\Omega_Y \to \Omega_Z$ be a contact transformation. Let $\Lambda'$ be the Lagrangian submanifold associated to the graph of $\chi'$. Let $K' \in D^b_{\cc}(Y \times Z, (p_Y,p_Z^a))$ satisfying \eqref{microhyp} and $s \in \muh(K',\OO_{Y \times Z}^{t(0,n)})$. Then $\varphi_s \circ \varphi'_{s'}=(\varphi \circ \varphi')_{s \circ s'}$ and $\psi_s \circ \psi'_{s'}=(\psi \circ \psi')_{s \circ s'}$, where $s \circ s'$ is the image of $s \otimes s'$ by the morphism
\begin{eqnarray*}
& \muh(K,\OO_{X \times Y}^t)_{(p_x,p_Y^a)} \otimes
\muh(K',\OO_{Y \times Z}^t))_{(p_Y,p_Z^a)} & \\
& \to \muh(K \mucirc K'[n],\OO_{X \times Z}^t)_{(p_X,p_Z^a)}.& \notag
\end{eqnarray*}
\item[(iii)] Let $P \in \E_{X,p_X}^{\R,f}$ and $Q \in \E^{\R,f}_{Y,p_Y}$ such that $Ps=sQ$. Then:
$$
P \circ \varphi_s = \varphi_s \circ Q
$$
(and similarly for $\psi_s$).
\end{itemize}
\end{prop}
\dim\ \ (i) Similar to Proposition 5.2.1 (i) of \cite{An94}. There exists a neighborhood $\Omega$ of $(p_X,p_Y^a)$ such that $s \in \Gamma(\Omega,\muh(K,\OO_{X \times Y}^{t(0,n)})$ and we may suppose that $\Lambda$ is closed in $\Omega$. Set $\K=\muh(K,\OO_{X \times Y}^{t(0,n)})$. Then
\begin{equation}\label{arrs}
s \in \Gamma(\Omega,\K) \simeq \Ho(\CC_\Lambda,\K).
\end{equation}
Moreover we can find a relatively compact neighborhoods $V_Y$ and $V_X$ of $\pi_Y(p_Y)$ and $\pi_X(p_X)$ respectively such that $\Phi^\mu_KF=\Phi_{K_{X \times V_Y}}F=K_{X \times V_Y} \circ F$ and $\Phi^\mu_KG=\Phi_{K_{V_X \times Y}}G=K_{V_X \times Y} \circ G$. Now set
\begin{eqnarray*}
\F_1=\muh(\Phi_{K_{X \times V_Y}}F,\ot_X), && \G_1=\muh(D'G,\OWX), \\
\F_2=\muh(F,\ot_Y)[n], && \G_2=\muh(D'(\Phi_{K_{V_X \times Y}}G),\OWY)[n].
\end{eqnarray*}
Then the morphisms $\varphi_s$ and $\psi_s$ are given by the diagrams
$$
\xymatrix{
\F_2|_{\Omega_Y} \ar[r]^{\sim\hspace{0.6cm}} & (\CC_{\Lambda}^a \circ \F_2)|_{\Omega_X} \ar[r] & (\K^a \circ \F_2)|_{\Omega_X} \ar[r] & \F_1|_{\Omega_X}
}
$$
$$
\xymatrix{
\G_2|_{\Omega_Y} \ar[r]^{\sim\hspace{0.6cm}} & (\CC_{\Lambda}^a \circ \G_2)|_{\Omega_X} \ar[r] & (\K^a \circ \G_2)|_{\Omega_X} \ar[r] & \G_1|_{\Omega_X}
}
$$
where the first arrows are given by \eqref{arrs} and the second ones by \eqref{microt} and \eqref{microw}.

(ii) The arrow follows from (i) and the associativity of the composition.

(iii) See \cite{An94}, Proposition 5.2.1 (iii).

\qed

\begin{teo} \label{An5.2.2} Let $\chi$ be a contact transformation from $\Omega_X$ to $\Omega_Y$ and let $\Lambda$ be the Lagrangian manifold associated to the graph of $\chi$. Then there exists $K \in D^b_{\cc}(X \times Y,(p_X,p_Y^a))$ satisfying \eqref{microhyp} and $s \in \muh(K,\OO_{X \times Y}^{t(0,n)})_{(p_X,p_Y^a)}$ such that:
\begin{itemize}
\item[(i)] the correspondence $\E_{X,p_X} \ni P \mapsto Q \in \E_{Y,p_Y}$ such that $Ps=sQ$ is an isomorphism of rings,
\item[(ii)] for each $F \in D^b_{\rc}(Y,p_Y)$ and $G \in D^b_{\rc}(X,p_X)$ the morphisms induced by $s$
\begin{eqnarray*}
\varphi_s:\muh(F,\ot_Y)_{p_Y}[n] & \to & \muh(\Phi^\mu_KF,\ot_X)_{p_X}  \\
\psi_s:\muh(D'(\Phi^\mu_KG),\OWY)_{p_Y}[n] & \to & \muh(D'G,\OWX)_{p_X}.
\end{eqnarray*}
are isomorphisms compatible with (i).
\end{itemize}
\end{teo}
\dim\ \ The proof is similar to the proof of Proposition 5.2.2 of \cite{An94}.\\
\qed

\begin{oss} Set $G=D'\Phi^\mu_KF$ with $F \in D^b_{\rc}(Y,p_Y)$, then $D'G=\Phi^\mu_KF$ and
$$
D'\Phi^\mu_KG \simeq D'\Phi^\mu_KD'\Phi^\mu_KF \simeq \Psi^\mu_KD'D'\Phi^\mu_KF \simeq \Psi^\mu_K\Phi^\mu_KF \simeq F,
$$
where the last isomorphism follows from Theorem 7.1.2 of \cite{KS90}. Hence we obtain the isomorphism
$$
\muh(F,\OWY)_{p_Y}[n] \iso \muh(\Phi^\mu_KF,\OWX)_{p_X}.
$$
\end{oss}

\section{Cauchy-Kowaleskaya-Kashiwara theorem}

Here we prove the Cauchy-Kowaleskaya-Kashiwara theorem for holomorphic functions with growth conditions.
The idea of the proof is the following: we divide the proof in two parts. In the first one we prove that for holomorphic functions with growth conditions the characteristic variety of a coherent $\D$-module coincides with the microsupport of the complex of solutions.
In the second part we use inverse image formulas for $\ol_Y$ to finish the proof of the theorem.

\subsection{Microsupport and characteristic variety}

We are going to study the relation between microsupport of subanalytic sheaves and the characteristic variety of a $\D$-module.

\begin{lem}\label{D'rho!} Let $F \in D^b_{\rc}(\CC_X)$ and let $G \in D^b(\CC_X)$. Then
$$
D'F \otimes \rho_!G \simeq \rh(F,\rho_!G).
$$
\end{lem}
\dim\ \ It follows from the following isomorphism of \cite{Pr1}: let $F \in D^b_{\rc}(\CC_X)$, $K \in D^b(\CC_{X_{sa}})$ and $G \in D^b(\CC_X)$. Then
$$
\rh(F,K) \otimes \rho_!G \simeq \rh(F,K \otimes \rho_!G).
$$
Setting $K=\CC_X$ we obtain the result.\\
\qed


Let us recall the notion of elliptic pair of \cite{SS94}. Let $\M$ be a coherent
$\D$-module and $F \in D^b_{\rc}(\CC_X)$, then $(F,\M)$ is an
elliptic pair if
$$ SS (F) \cap \vchar(\M) \subseteq T^*_XX.$$
We consider the sheaf $\OO^\lambda_X$, for
$\lambda=\varnothing,t,\mathrm{w},\omega$.

\begin{prop}\label{ellpairFG} Let $(F,\M)$ be an elliptic pair. Let $G \in D^b(\rho_!\D_X)$ such that $H^k\muh(F,G)$ is a $\E_X$-module for each $k \in \Z$. Then we have the
isomorphism
$$
 \rh_{\D_X}(\M,D'F \otimes \imin \rho G) \iso
 \rh_{\D_X}(\M,\imin \rho\rh(F,G)).
$$
\end{prop}
\dim\ \ Let $\delta: \Delta \to X \times X$ be the embedding and
let us consider the Sato's triangle
\begin{eqnarray}\label{Satos}
& \imin \delta\rh(\imin {q_1}F,q_2^!G) \otimes
\omega_{\Delta|X
\times X}  \to  \delta^!\rh(\imin {q_1}F,q_2^!G) & \notag \\
 & \to R\dot{\pi}_*\muh^{sa}(F,G) \stackrel{+}{\to}. &
\end{eqnarray}
We have $\delta^!\rh(\imin {q_1}F,q_2^!G) \simeq
\rh(F,G)$. Moreover
\begin{eqnarray*}
\imin \delta\rh(\imin {q_1}F,q_2^!G) \otimes
\omega_{\Delta|X
\times X} & \simeq & \imin \delta\rh(\imin {q_1}F,\imin {q_2}G) \\
& \simeq & \imin \delta (D'F \boxtimes G) \\
& \simeq & D'F \otimes G
\end{eqnarray*}
where the second isomorphism follows from Lemma \ref{rhD'}. Hence applying $\imin \rho$ we obtain
$$
\dt{D'F \otimes \imin \rho G}{\imin \rho \rh(F,G)}{
R\dot{\pi}_*\muh(F,G)}.
$$
Applying the functor $\rh_{\D_X}(\rho_!\M,\cdot)$ we obtain
\begin{eqnarray*}
& \rh_{\D_X}(\M,D'F \otimes \imin \rho G) \to
 \rh_{\D_X}(\M,\imin \rho\rh(F,G)) & \\
& \to \rh_{\D_X}(\M,R\dot{\pi}_*\muh(F,G))
\stackrel{+}{\to}. &
\end{eqnarray*}
Then it is enough to prove that $\rh_{\D_X}(\M,
R\dot{\pi}_*\muh(F,G))=0$. First remark that we have by adjunction
\begin{eqnarray*}
\lefteqn{\rh_{\D_X}(\M, R\ppi_*\muh(F,G))}\\ &
\simeq & R\ppi_*\rh_{\imin \pi \D_X}(\imin \pi \M,
\muh(F,G)).
\end{eqnarray*}
Let $k \in \Z$. Since $H^k\muh(F,G)$ is a $\E_X$-module for each $k \in \Z$ we have
\begin{eqnarray*}
\lefteqn{R\ppi_*\rh_{\imin \pi \D_X}(\imin \pi \M,H^k
\muh(F,G))} \\
& \simeq & R\ppi_*\rh_{\E_X}((\E_X \otimes_{\imin \pi \D_X} \imin
\pi \M),H^k \muh(F,G)).
\end{eqnarray*}
We have $\supp(H^k\muh(F,G)) \subseteq SS(F)$ for each $k \in \Z$ by Corollary \ref{muhSS(F)G} and
$\supp(\E_X \otimes_{\imin \pi \D_X} \imin \pi \M)=\vchar(\M)$.
Hence
\begin{equation}\label{Hk=0}
R\ppi_*\rh_{\imin \pi \D_X}(\imin \pi \M,H^k
\muh(F,G))=0
\end{equation}
for each $k \in \Z$ since the pair $(F,\M)$ is elliptic. Let us
suppose that the length  of the bounded complex
$\muh(F,G)$ is $n$ and let us argue by induction on the
truncation $\tau^{\leq i}\muh(F,G)$. If $i=0$ the result
follows from \eqref{Hk=0}. Let us consider the distinguish
triangle
$$
\dt{\tau^{\leq n-1} \muh(F,G)}{
\muh(F,G)}{H^n \muh(F,G)}
$$
and apply the functor $R\dot{\pi}_*\rh_{\imin \pi\D_X}(\imin \pi
\M,\cdot)$. The first term becomes zero by the induction hypothesis
and the third one is zero by \eqref{Hk=0}. Hence
$\rh_{\D_X}(\M,R\dot{\pi}_*\muh(F,G))=0$ and
the result follows.\\
\qed

Setting $G=\ol_X$ in Proposition \ref{ellpairFG}, $\lambda=\varnothing,t,{\rm w},\omega$ we obtain the following result.

\begin{teo}\label{ellpair} Let $(F,\M)$ be an elliptic pair. Then we have the
isomorphism
\begin{equation}\label{isoellpair}
 \rh_{\D_X}(\M,D'F \otimes \OO_X) \iso
 \rh_{\D_X}(\M,\imin \rho\rh(F,\OO^\lambda_X)).
\end{equation}
\end{teo}
\dim\ \ If $\lambda=\omega$ this is a consequence of Lemma \ref{D'rho!}. If $\lambda=\varnothing,t,{\rm w}$ then $H^k\muh(F,\ol_X)$ is a $\E_X$-module for each $k \in \Z$ by Corollaries \ref{safreccia1} and \ref{safreccia2} and the result follows from Proposition \ref{ellpair}.\\
\qed

\begin{es}\label{ellsis} Let $M$ be a real analytic manifold and let $X$ be a complexification of $M$. Let $\M$ be an elliptic system on $M$. Then
\begin{eqnarray*}
\rh_{\D_X}(\M,\A_M) & \simeq & \rh_{\D_X}(\M,\C^\infty_M) \\
& \simeq & \rh_{\D_X}(\M,\db_M) \\
& \simeq & \rh_{\D_X}(\M,\B_M).
\end{eqnarray*}
This follows from Theorem \ref{ellpair} setting $F=D'\CC_M$ and applying $\imin \rho$ to the isomorphism \eqref{isoellpair}. In fact
$$ \rh_{\D_X}(\M,\imin\rho\rh(D'\CC_M,\OO^\lambda_X)) \simeq
\begin{cases}
    \text{$\rh_{\D_X}(\M,\A_M)$ if $\lambda=\omega$}, \\
    \text{$\rh_{\D_X}(\M,\C^\infty_M)$ if $\lambda={\rm w}$}, \\
    \text{$\rh_{\D_X}(\M,\db_M)$ if $\lambda=t$}, \\
    \text{$\rh_{\D_X}(\M,\B_M)$ if $\lambda=\varnothing$}
  \end{cases}
$$
and they are all isomorphic since we have $\imin \rho \rh_{\rho_!\D_X}(\rho_!\M,D'D'\CC_M \otimes \OO^\lambda_X) \simeq \rh_{\D_X}(\M,\CC_M \otimes \OO_X) \simeq \rh_{\D_X}(\M,\A_M)$ for $\lambda=\varnothing,t,{\rm w},\omega$.
\end{es}

\begin{es} Let $X$ be a smooth submanifold of $Y$, let $\OO_{Y}\widehat{|}_X$ be the formal completion of $X$ along $Y$ and let $\B_{X|Y}$ be the algebraic cohomology of $\OO_Y$ with support in $X$. Let $\M$ be a $\D_Y$-module such that $T^*_XY \cap \vchar(\M) \subseteq T^*_YY$. Then we have the isomorphisms
\begin{eqnarray*}
\rh_{\D_X}(\M,\OO_{Y|X}) & \simeq & \rh_{\D_X}(\M,\OO_Y\widehat{|}_X), \\
 & \simeq & \rh_{\D_X}(\M,\B_{X|Y}).
\end{eqnarray*}
This follows from Theorem \ref{ellpair} setting $F=D'\CC_Y$, applying $\imin \rho$ to the isomorphism \eqref{isoellpair} and arguing as in Example \ref{ellsis}.
\end{es}


Let $\M$ be a $\D_X$-module and let $\lambda=\varnothing,t,\mathrm{w},\omega$. One sets for short
$$\sol^\lambda(\M):=\rh_{\rho_!\D_X}(\rho_!\M,\ol_X).$$

\begin{cor}\label{ssvchar} Let $\M$ be a coherent $\D_X$-module. Then
$$
SS(\sol^\lambda(\M))=\vchar(\M).$$
\end{cor}
\dim\ \ Recall that $SS(\sol(\M))=\vchar(\M)$.

(i) $\vchar(\M) \subseteq SS(\sol^\lambda(\M))$ follows from the
fact that $\imin \rho \sol^\lambda(\M)=\sol(\M)$ and $SS(\imin
\rho G) \subseteq SS(G)$ for each $G \in D^b(\CC_{X_{sa}})$.

(ii) $\vchar(\M) \supseteq  SS(\sol^\lambda(\M))$. Let $(x,\xi)
\notin \vchar(\M)=SS(\sol(\M))$
and let $U$ be a conic open neighborhood of $(x,\xi)$, $U \cap \vchar(\M)=\varnothing$,
such that for each $F \in D^b_{\rc}(\CC_X)$ with $\supp(F) \subset \subset \pi(U)$ and $SS(F) \subset U \cup T^*_XX$ we have $\Ho_{D^b(\CC_X)}(F,\sol(\M))=0$. By
Theorem \ref{ellpair} the complexes
\begin{eqnarray*}
\rh(F,\rh_{\D_X}(\M,\OO_X)) & \simeq & \imin \rho \rh(F,
\rh_{\rho_!\D_X}(\rho_!\M,R\rho_*\OO_X))
\\
& \simeq & \imin \rho \rh(F,\rh_{\rho_!\D_X}(\rho_!\M,\ol_X))
\end{eqnarray*}
are all quasi-isomorphic for
$\lambda=\varnothing,t,\mathrm{w},\omega$.
Hence $\Ho_{D^b(\CC_{X_{sa}})}(F,\sol^\lambda(\M))=0$ and $(x,\xi) \notin SS(\sol^\lambda(\M))$.\\
\qed

\subsection{Cauchy-Kowaleskaya-Kashiwara theorem}

Now we apply the preceding results to prove the Cauchy-Kowaleskaya-Kashiwara theorem for holomorphic functions with growth conditions $\lambda=\varnothing,t,\mathrm{w},\omega$. We refer to \cite{Ka03} for the statement and proof of the Cauchy-Kowaleskaya-Kashiwara theorem for holomorphic functions. \\

Let $f:X \to Y$ be a morphism of complex manifolds. Set $d=\mathrm{dim}_\CC X-\mathrm{dim}_\CC Y$. We recall the inverse image isomorphisms
\begin{eqnarray}
f^!\ot_Y & \simeq & \rho_!\ddyx \underset{\rho_!\D_X}{\ltens}\ot_X[d], \label{iminvot} \\
f^!\OWY & \simeq & \rh_{\rho_!\D_X}(\rho_!\ddxy,\OWX)[2d]. \label{iminvow}
\end{eqnarray}

\begin{prop}\label{division} Let $\M$ be a coherent $\D_Y$-module, and suppose that $f$ is non characteristic for $\M$.  Then we have the following isomorphism for $\lambda=\varnothing,t,\mathrm{w},\omega$:
$$
f^! \rh_{\rho_!\D_Y}(\rho_!\M,\ol_Y) \simeq \rh_{\rho_!\D_X}(\rho_!\imin {\underline{f}} \M,\ol_X)[2d].
$$
\end{prop}
\dim\ \ (i) Let $\lambda=t$. Recall that if $\M$ is a coherent $\D_Y$-module and $f$ is non characteristic, then $\imin {\underline{f}} \M$ is a coherent $\D_X$-module and
$$
\imin {\underline{f}} \rh_{\D_Y}(\M,\D_Y) \simeq \rh_{\D_X}(\imin {\underline{f}} \M,\D_X)[d].
$$
We have the chain of isomorphisms
\begin{eqnarray*}
\rh_{\D_X}(\rho_!\imin {\underline{f}}\M,\ot_X)[2d] & \simeq &  \rho_!\rh_{\D_X}(\imin {\underline{f}}\M,\D_X) \otimes_{\rho_!\D_X} \ot_X[2d] \\
& \simeq & \rho_!\imin {\underline{f}} \rh_{\D_Y}(\M,\D_Y) \otimes_{\rho_!\D_X} \ot_X[d] \\
& \simeq & \rho_!\imin f \rh_{\D_Y}(\M,\D_Y) \otimes_{\rho_!\imin f\D_Y} f^!\ot_Y \\
& \simeq & f^!(\rho_!\rh_{\D_Y}(\M,\D_Y) \otimes_{\rho_!\D_Y} \ot_Y) \\
& \simeq & f^!\rh_{\rho_!\D_Y}(\rho_!\M,\ot_Y),
\end{eqnarray*}
where the first and the last isomorphisms follow from the coherence of $\imin {\underline{f}} \M$ and $\M$, and the third one follows from \eqref{iminvot}. \\

\noindent (ii) Let $\lambda=\mathrm{w}$. We have the chain of isomorphisms
\begin{eqnarray*}
f^!\rh_{\rho_!\D_Y}(\M,\OWY) & \simeq & \rh_{\rho\imin f \D_X}(\rho_!\imin f \M,f^!\OWY) \\
& \simeq & \rh_{\rho_!\imin f\D_X}(\rho_!\imin f \M,\rh_{\rho_!\D_X}(\rho_!\ddxy,\OWX))[2d] \\
& \simeq & \rh_{\rho_!\D_X}(\rho_!\imin {\underline{f}} \M,\OWX)[2d],
\end{eqnarray*}
where the second isomorphism follows from \eqref{iminvow}.\\

\noindent (iii) Let $\lambda=\varnothing,\omega$. Since $\M$ is coherent and $f$ is non characteristic the result follows from the isomorphism
$$
f^! \rh_{\D_Y}(\M,\OO_Y) \simeq \rh_{\D_X}(\imin {\underline{f}} \M,\OO_X)[2d].
$$
\qed

\begin{teo}\label{ckk} Let $\M$ be a coherent $\D_Y$-module, and suppose that $f$ is non characteristic for $\M$. Then we have the following isomorphism for $\lambda=\varnothing,t,\mathrm{w},\omega$:
$$
\imin f \rh_{\rho_!\D_Y}(\rho_!\M,\OO^\lambda_Y) \simeq \rh_{\rho_!\D_X}(\rho_!\imin {\underline{f}} \M,\ol_X).
$$
\end{teo}
\dim\ \ By Corollary \ref{ssvchar} $f$ is non characteristic for $SS(\sol^\lambda(\M))$. Hence by Proposition \ref{noncharSS(F)}
$$
f^! \rh_{\rho_!\D_Y}(\rho_!\M,\ol_Y) \simeq \imin f \rh_{\rho_!\D_Y}(\rho_!\M,\ol_Y) [2d].
$$
Then the result follows from Proposition \ref{division}.\\
\qed


\appendix

\section{Appendix}

\subsection{Review on subanalytic sets}

We recall briefly some properties of subanalytic subsets.
Reference are made to \cite{BM88} and \cite{Lo93} for the theory
of subanalytic subsets and to \cite{Co00} and \cite{VD98} for the
more general theory of o-minimal structures. Let $X$ be a real
analytic manifold.

\begin{df} Let $A$ be a subset of $X$.

\begin{itemize}

\item[(i)] $A$ is said to be semi-analytic if it is locally
analytic, i.e. each $x \in A$ has a neighborhood $U$ such that $X
\cap U=\cup_{i \in I}\cap_{j \in J} X_{ij}$, where $I,J$ are
finite sets and either $X_{ij}=\{y \in U_x;\,f_{ij}>0\}$ or
$X_{ij}=\{y \in U_x;\,f_{ij}=0\}$ for some analytic function
$f_{ij}$.

\item[(ii)] $A$ is said to be subanalytic if it is locally a
projection of a relatively compact semi-analytic subset, i.e. each
$x \in A$ has a neighborhood $U$ such that there exists a real
analytic manifold $Y$ and a relatively compact semi-analytic
subset $A'\subset X \times Y$ satisfying $X \cap U=\pi(A')$, where
$\pi:X \times Y \to X$ denotes the projection.

\item[(iii)] Let $Y$ be a real analytic manifold. A continuous map
$f:X \to Y$ is subanalytic if its graph is subanalytic in $X
\times Y$.
\end{itemize}
\end{df}

Let us recall some result on subanalytic subsets.

\begin{prop} Let $A,B$ be subanalytic subsets of $X$. Then $A \cup
B$, $A \cap B$, $\overline{A}$, $\partial A$ and $A \setminus B$
are subanalytic.
\end{prop}

\begin{prop} Let $A$ be a subanalytic subsets of $X$. Then the
connected components of $A$ are locally finite.
\end{prop}

\begin{prop} Let $f:X \to Y$ be a subanalytic map. Let $A$ be a
relatively compact subanalytic subset of $X$. Then $f(A)$ is
subanalytic.
\end{prop}

\begin{df} A simplicial complex $(K,\Delta)$ is the data
consisting of a set $K$ and a set $\Delta$ of subsets of $K$
satisfying the following axioms:
\begin{itemize}
\item[S1] any $\sigma \in \Delta$ is a finite and non-empty subset
of $K$,
\item[S2] if $\tau$ is a non-empty subset of an element $\sigma$
of $\Delta$, then $\tau$ belongs to $\Delta$,
\item[S3] for any $p \in K$, $\{p\}$ belongs to $\Delta$,
\item[S4] for any $p \in K$, the set $\{\sigma \in \Delta; p \in
\sigma\}$ is finite.
\end{itemize}
\end{df}

If $(K,\Delta)$ is a simplicial complex, an element of $K$ is
called a vertex. Let $\R^K$ be the set of maps from $K$ to $\R$
equipped with the product topology. To $\sigma \in \Delta$ one
associate $|\sigma|\subset \R^K$ as follows:
$$|\sigma|=\left\{x \in \R^K;\text{ $x(p)=0$ for $p \notin
\sigma$, $x(p)>0$ for $p \in \sigma$ and } \sum_p
x(p)=1\right\}.$$ As usual we set:
$$|K|=\bigcup_{\sigma \in \Delta}|\sigma|,$$
$$U(\sigma)=\bigcup_{\tau \in \Delta,\tau \supset \sigma}|\tau|,$$
and for $x \in |K|$:
$$U(x)=U(\sigma(x)),$$
where $\sigma(x)$ is the unique simplex such that $x \in
|\sigma|$.

\begin{teo} Let $X=\bigsqcup_{i \in I} X_i$ be a locally finite
partition of $X$ con\-sisting of subanalytic subsets. Then there
exists a simplicial complex $(K,\Delta)$ and a subana\-lytic
homeomorphism $\psi:|K|\iso X$ such that
\begin{itemize}
\item[(i)] for any $\sigma \in \Delta$, $\psi(|\sigma |)$ is a
subanalytic submanifold of $X$,
\item[(ii)] for any $\sigma \in \Delta$ there exists $i \in I$
such that $\psi(|\sigma |) \subset X_i$.
\end{itemize}
\end{teo}

Let us recall the definition of a subfamily of the subanalytic
subsets of $\R^n$ which have some very good properties.

\begin{df} A subanalytic subset $A$ of $\R^n$ is said to be globally
subana\-lytic if it is subanalytic in the projective space
$\mathbb{P}^n(\R)$. Here we identify $\R^n$ with a submanifold of
$\mathbb{P}^n(\R)$ via the map $(x_1,\ldots,x_n) \mapsto
(1:x_1:\ldots :x_n)$.
\end{df}

An equivalent way to define globally subanalytic subsets is by
means of the map $\tau_n:\R^n \to \R^n$ given by
$$\tau_n(x_1,\ldots,x_n) := \left({x_1 \over
\sqrt{1+x_1^2}},\ldots,{x_n \over \sqrt{1+x_n^2}}\right).$$ In
particular relatively compact subanalytic subsets are globally
subanalytic.

\begin{df} A map $f:\R^n \to \R^n$ is said to be globally
subanalytic if its graph is globally subanalytic.
\end{df}

\begin{prop} Let $f:\R^n \to \R^n$ be a globally subanalytic map. Let $A$ be a
globally subanalytic subset of $\R^n$. Then $f(A)$ is globally
subanalytic.
\end{prop}

Now we recall the notion of cylindrical cell decomposition, a
useful tool to study the geometry  of a subanalytic subset. We
refer to \cite{Co00} and \cite{VD98} for a complete exposition.\\

A cyindrical cell decomposition (ccd for short) of $\R^n$ is a
finite partition of $\R^n$ into subanalytic subsets,
called the cells of the ccd. It is defined by induction on $n$:

\begin{itemize}

\item[n=1] a ccd of $\R$ is given by a finite subdivision
$a_1<\ldots <a_\ell$ of $\R$. The cells of $\R$ are the points
$\{a_i\}$, $1\leq i \leq \ell$, and the intervals $(a_i,a_{i+1})$,
$0 \leq i \leq \ell$, where $a_0=-\infty$ and
$a_{\ell+1}=+\infty$.

\item[n>1] a ccd of $\R^{n}$ is given by a ccd of $\R^{n-1}$ and,
for each cell $D$ of $\R^{n-1}$, continuous analytic functions
$$\zeta_{D,1}<\ldots <\zeta_{D,\ell_D}:D \to \R.$$
The cells of $\R^n$ are the graphs
$$\{(x,\zeta_{D,i}(x));\; x \in D\},\ \ 1 \leq i \leq \ell_D,$$
and the bands
$$\{(x,y) \in D \times \R;\;
\zeta_{D,i}(x)<y<\zeta_{D,i+1}(x)\}$$ for $0\leq i \leq \ell_D$,
where $\zeta_{D,0}=-\infty$ and $\zeta_{D,\ell_D+1}=+\infty$.
\end{itemize}

\begin{teo} Let $A_1, \ldots A_k$ be globally subanalytic subsets of $\R^n$. There exists a ccd of
$\R^n$ such that each $A_i$ is a union of cells.
\end{teo}

We end this section with the following useful result.

\begin{lem}\label{goodcover} Let $U$ be a globally subanalytic subset of $\R^n$ and denote by $\pi:\R^n \to \R^{n-1}$ the projection. Then $U$ admits a finite open covering $\{U_i\}$ such that each $U_i$ is simply connected and the intersection of each $U_i$ with the fibers of $\pi$ is contractible or empty.
\end{lem}
\dim\ \ Up to take the image of $U$ by the homeomorphism
\begin{eqnarray*}
\varphi:\R^n & \to & (-1,1)^n \\
(x_1,\ldots,x_n) & \mapsto & \left({x_1 \over
\sqrt{1+x_1^2}},\ldots,{x_n \over \sqrt{1+x_n^2}}\right)\end{eqnarray*}
we may assume that $U$ is bounded. Then it follows from a result of \cite{Wi05} that $U$ can be covered by finitely many open cells, and cells satisfy the desired properties.\\
\qed

\subsection{Ind-sheaves and subanalytic sites}\label{indsa}
Let us recall some results of \cite{KS01}.

One denotes by $\I(k_X)$ the category of ind-sheaves of $k$-vector
spaces on $X$, that is $\I(k_X) = \ind(\mod^c(k_X)),$ where
$\mod^c(k_X)$ denotes the full subcategory of $\mod(k_X)$
consisting of sheaves with compact support on $X$. We denote by
$D^b(\I(k_X))$ the bounded derived category of $\I(k_X)$.\\

There are three functors relating ind-sheaves and classical
sheaves:
\begin{eqnarray*}
\iota:\mod(k_X) \to \I(k_X) && F \mapsto \indl {U
\subset\subset X}F_U,\\
\alpha:\I(k_X) \to \mod(k_X) && \indl i F_i \mapsto \lind i
F_i,\\
\beta:\mod(k_X) \to \I(k_X) && \text{left adjoint to $\alpha$}.
\end{eqnarray*}
These functors satisfy the following properties:
\begin{itemize}
\item the functor $\iota$ is fully faithful, exact and commutes with $\Lpro$,
\item the functor $\alpha$ is exact and commutes with $\Lind$ and $\Lpro$,
\item the functor $\beta$ is fully faithful, exact and commute with $\Lind$,
\item $(\alpha,\iota)$ and $(\beta,\alpha)$ are pairs of
adjoint functors.
\end{itemize}
Since $\iota$ is fully faithful and exact we identify $\mod(k_X)$
(resp. $D^b(k_X)$) with a full abelian subcategory of $\I(k_X)$
(resp. $D^b(\I(k_X))$).

The category $\I(k_X)$ admits an internal hom denoted by $\ih$ and
this functor admits a left adjoint, denoted by $\otimes$.
One can also define an external $\ho:\I(k_X) \times \I(k_X) \to
\mod(k_X)$
and one has
$$\ho(F,G) = \alpha \ih(F,G) \ \ \text{and}\ \
\Ho_{\I(k_X)}(F,G)=\Gamma(X;\ho(F,G)).$$

The functor $\otimes$ is exact while $\ih$ and $\ho$ are left
exact and admit right derived functors $\ri$ and $\rh$.

Consider a morphism of real analytic manifolds $f:X \to Y$. One
defines the external operations
\begin{eqnarray*}
\imin f :\I(k_Y) \to \I(k_X) && \indl i G_i \mapsto \indl {i,U\subset\subset X} (\imin f G_i)_U \\
f_*:\I(k_X) \to \I(k_Y) && \indl i F_i \mapsto \lpro {U \subset\subset X} \lind i f_*\Gamma_UF_i \\
f_{!!}:\I(k_X) \to  \I(k_Y)  && \indl i F_i \mapsto \indl i f_!F_i
\end{eqnarray*}
where the notation $f_{!!}$ is chosen to stress the fact that
$f_{!!} \circ \iota \not\simeq \iota \circ f_!$ in general.

While $\imin f$ is exact, the others functors admit right derived
functors. One can show that the functor $Rf_{!!}$ admits a right
adjoint denoted by $f^!$ and we get the usual formalism of the six
Grothendieck operations. Almost all the formulas of the classic
theory of sheaves remain valid for ind-sheaves.\\

There is a strict relations between ind-sheaves and sheaves on the
subanalytic site associated to $X$. Set
$\I_{\rc}(k_X)=\ind(\mod^c_{\rc}(k_X))$ for short.

\begin{teo}\label{indshsa} One has an equivalence of categories
\begin{eqnarray*}
\I_{\rc}(k_X) & \iso & \mod(k_{X_{sa}}) \\
\indl i F_i & \mapsto & \lind i \rho_*F_i.
\end{eqnarray*}
\end{teo}


Let us recall the following functor defined in \cite{KS01}:
\begin{eqnarray*}
I_\T:\mod(k_{X_{sa}}) & \to  & \I(k_X)\\
 \lind i \rho_* F_i & \mapsto & \indl i F_i.
\end{eqnarray*}
It is fully faithful, exact and commutes with $\Lind$ and
$\otimes$. It admits a right adjoint
\begin{eqnarray*}
J_\T:\I(k_X) & \to & \mod(k_{X_{sa}})
\end{eqnarray*}
satisfying, for each $U \in \op(X_{sa})$, $\Gamma(U;J_\T F) =
\Ho_{\I(k_X)}(k_U,F)$. This functor is right exact and commutes
with filtrant inductive limits. Moreover we have $RJ_\T \circ I_\T
\simeq \id$ and
$$RJ_\T\ri(I_\T F,G) \simeq \rh(F,RJ_\T G).$$

We have the following relations:

\begin{eqnarray*}
RJ_\T \circ \iota \simeq R\rho_* & \text{and} & \alpha \simeq
\imin
\rho \circ J_\T \\
\alpha \circ I_\T \simeq \imin \rho & \text{and} & I_\T \circ
\rho_! \simeq \beta
\end{eqnarray*}

 Let $f:X \to Y$ be a morphism of real analytic manifolds and let $U$ be an open subanalytic subset of $X$.

\begin{lem}
Let $F \in D^b(k_{X_{sa}})$ and $G \in D^b(k_{Y_{sa}})$. We have
\begin{itemize}
\item[(i)] $I_{\T} \circ Rf_{!!}F \simeq Rf_{!!} \circ I_{\T}F$,
\item[(ii)] $I_\T \circ \imin fG \simeq \imin f \circ I_\T G$,
\item[(iii)] $I_\T \circ f^!G \simeq f^! \circ I_\T G$,
\item[(iv)] $I_\T F_U \simeq (I_\T F)_U$,
\item[(v)] $I_\T \circ \mathrm{R}\Gamma_UF \simeq R\I\Gamma_U \circ I_\T F$.
\end{itemize}
\end{lem}

\subsection{Inverse image for tempered holomorphic functions}

Let $f:M \to N$ be a morphism of oriented real analytic manifolds of dimension $d_M$ and $d_N$. Set $d=d_N-d_M$. Denote by $\A_M$ the sheaf of analytic functions on $M$.

\begin{lem}\label{lem:uno}
Let $F$ be an $\A_M$-module locally free of finite rank. Then
$R^kf_{!!}(\dbt_M\otimes_{\rho_!\A_M}\rho_! F)=0$ for $k\neq 0$.
\end{lem}
\dim\ \ It is a consequence of the fact that $\dbt_M$ is quasi-injective and Proposition 1.6.5 of \cite{Pr1}.\\
\qed

\begin{lem}\label{lem:due}
Let $M$ and $N$ be orientable real manifolds. There is a natural
morphism of complexes
$$f_{!!}(\dbt_M\underset{\rho_!\A_M}{\otimes}\rho_!\Omega^\bullet_M)[d_M]\rightarrow \dbt_N\underset{\rho_!\A_N}{\otimes}\rho_!\Omega^\bullet_N[d_N]. $$
\end{lem}
\dim\ \ Let $U \in \op^c(N_{sa})$. We have the
chain of morphisms
\begin{eqnarray*}
\Gamma(U;f_{!!}(\dbt_M \underset{\rho_! \A_M}{\otimes}\rho_!\Omega_M^{d_M-i})) & \simeq &
\Gamma(N;\ho(\CC_U,f_{!!}(\dbt_M
\underset{\rho_!\A_M}{\otimes}\rho_!\Omega_M^{d_M-i}))) \\
& \simeq & \Gamma(N;f_!\tho(\imin f \CC_U,\db_M
\underset{\A_M}{\otimes}\Omega_M^{d_M-i})) \\
& \to & \Gamma_c(N;\tho(G,\db_N \underset{\A_N}{\otimes}\Omega_N^{d_N-i})) \\
& \simeq & \Gamma(U;\db_N \underset{\A_N}{\otimes}\Omega_N^{d_N-i}),
\end{eqnarray*}
where the third morphism follows from Proposition 4.3 of
\cite{KS96}.\\ \qed

\begin{prop}\label{mor:imminv}There
is a natural morphism in $D^b(\rho_!\D_M^{op})$:
\begin{equation}\label{mor:phi}
Rf_{!!}(\db^{t\vee}_M\underset{\rho_!\D_M}{\overset{L}{\otimes}}\rho_!\ddmn)\to\db^{t\vee}_N.
\end{equation}
\end{prop}
\dim\ \ The Spencer resolution of $\ddmn$ give rise to the
quasi-isomorphism

$$ \ddmn\stackrel{\sim}{\leftarrow} \D_M\underset{\A_M}{\otimes}\overset{\bullet}{\bigwedge}\Theta_M\underset{\A_M}
{\otimes}\ddmn\simeq\D_M\underset{\A_M
}{\otimes}\overset{\bullet}{\bigwedge}\Theta_M\underset{f^{-1}\A_M
}{\otimes}f^{-1}\D_N
$$
from which we obtain the following quasi-isomorphism for
$\db^{t\vee}_M\underset{\rho_!\D_M}{\otimes}\rho_!\ddmn$ in $ D^b(\rho_!
f^{-1}\D_N^{op}): $

\begin{eqnarray*}
\db^{t\vee}_M\underset{\rho_!\A_M}{\otimes}\rho_!\ddmn  &  \simeq  &
(\dbt_M\underset{\rho_!\D_M}{\otimes}\rho_!\Omega_M)\underset{\rho_!\D_M}{\otimes}(\rho_!(\D_M\underset{\A_M}{\otimes}\overset{\bullet}{\bigwedge}\Theta_M\underset{f^{-1}\A_N}{\otimes}f^{-1}\D_N))
\\
  &  \simeq  &  \dbt_M\underset{\rho_!\A_M}{\otimes}\rho_!(\Omega_M\underset{\A_M}{\otimes}\overset{\bullet}
  {\bigwedge}\Theta_M\underset{f^{-1}\A_N}{\otimes}f^{-1}\D_N)  \\
  &  \simeq  &  \dbt_M\underset{\rho_!\A_M}{\otimes}\rho_!(\Omega_M^\bullet\underset{f^{-1}\A
  N}{\otimes}f^{-1}\D_N)[d_M].
\end{eqnarray*}
Applying $Rf_{!!}$ we obtain:
\begin{eqnarray*}
Rf_{!!}(\db^{t\vee}_M\underset{\rho_!\D_M}{\otimes}\rho_!\ddmn)  &  \simeq &
Rf_{!!}(\dbt_M\underset{\rho_!\A_M}{\otimes}
\rho_!(\Omega_M^\bullet\underset{f^{-1}\A_N}{\otimes}f^{-1}\D_N))[d_M]  \\
  &  \simeq  &  Rf_{!!}(\dbt_M\underset{\rho_!\A_M}{\otimes}\rho_!\Omega_M^\bullet)\underset{\rho_!\A_N}{\otimes}
  \rho_!\D_N[d_M]  \\
  &  \simeq  &  f_{!!}(\dbt_M\underset{\rho_!\A_M}{\otimes}\rho_!\Omega_M^\bullet)\underset{\rho_!\A_N}{\otimes}\rho_!\D_N[d_M]  \\
  &  \to     &  \dbt_N\underset{\rho_!\A_N}{\otimes}\rho_!\Omega_N^\bullet\underset{\rho_!\A_N}{\otimes}\rho_!\D_N[d_N]  \\
  &  \simeq  &  \dbt_N\underset{\rho_!\A_N}{\otimes}\rho_!\Omega_N \\
  &  =       &  \db^{t\vee}_N,
\end{eqnarray*}
where the third isomorphism follows from Lemma \ref{lem:uno} and
the morphism from Lemma \ref{lem:due}. \\
\qed
By adjunction we get a morphism
\begin{equation}\label{mor:psi}
\db^{t\vee}_M\underset{\rho_!\D_M}{\overset{L}{\otimes}}\rho_!\ddmn\to
f^!\db^{t\vee}_N.
\end{equation}

\begin{teo}\label{iso:imminv} The morphism \eqref{mor:psi} is an isomorphism.
\end{teo}
\dim\ \ Let $F \in D^b_{\rc}(\CC_M)$ with compact support. We have the chain of isomorphisms
\begin{eqnarray*}
\Rh(F,f^!\db^{t\vee}_N) & \simeq & \Rh(Rf_{!!}F,\db^{t\vee}_N) \\
& \simeq & \mathrm{R}\Gamma(N,\tho(Rf_!F,\db^\vee_N)) \\
& \simeq & \mathrm{R}\Gamma(N,Rf_!(\tho(F,\db^\vee_M)\underset{\D_M}{\overset{L}{\otimes}}\ddmn)) \\
& \simeq & \mathrm{R}\Gamma(M,\tho(F,\db^\vee_M)\underset{\D_M}{\overset{L}{\otimes}}\ddmn) \\
& \simeq &
\Rh(F,\db^{t\vee}_M\underset{\rho_!\D_M}{\overset{L}{\otimes}}\rho_!\ddmn), \end{eqnarray*}
where the third isomorphism follows from Theorem 4.4 of \cite{KS96}.\\
\qed

By the equivalence between left and right
$\mathcal{D}$-modules, we have an isomorphism
\begin{equation}
\label{eq:DYXDbXDbY}
\rho_!\ddnm\underset{\rho_!\D_M}{\overset{L}{\otimes}}\dbt_M\stackrel{\sim}{\rightarrow}
f^!\dbt_N.
\end{equation}

\begin{cor}\label{smoothdbt} When $f$ is smooth we have an isomorphism
$$
\imin f\dbt_N \iso \rh_{\rho_!\D_M}(\rho_!\ddmn,\dbt_M).
$$
\end{cor}
\dim\ \ The result is
obtained by the following isomorphisms
\begin{eqnarray*}
\rh_{\rho\D_M}(\rho_!\ddmn,\dbt_M) & \simeq & \rho_!\rh(\ddmn,\D_M) \underset{\rho_!\D_M}{\ltens} \dbt_M \\
& \simeq & \rho_!\ddnm \underset{\rho_!\D_M}{\ltens} \dbt_M[d]\\
& \simeq & f^!\dbt_N[d]\\
& \simeq & f^{-1}\dbt_N.
\end{eqnarray*}
The first isomorphism is obtained by replacing $\ddmn$ with its
Koszul complex. The second follows from the smoothness of $f$ and
the isomorphism
$$\rh_{\D_M}(\ddmn,\D_M) \simeq \ddnm[d].$$
The last isomorphism follows since when $f$ is smooth we have the isomorphism $f^!(\cdot)[d]\simeq
f^{-1}$.\\ \qed

From now on $X$ will be a complex manifold, with structure sheaf
$\OO_X$. We denote by $\overline{X}$ the complex conjugate manifold
(with structure sheaf $\mathcal{O}_{\overline{X}}$), and $X_\R$
the underlying real analytic manifold, identified with the
diagonal of $X \times \overline{X}$. Let $\ot_X$ be the sheaf of tempered holomorphic functions on $X$.
We also consider the sheaf $\Omega^t_X \in
D^b(\rho_!\D_X ^{op})$:
$$\Omega^t_X:=\db^{t\vee}_{X_\R} \underset{\rho_!\mathcal{D}_{\overline{X}}}{\ltens} \rho_!\mathcal{O}_{\overline{X}}[-d_X].$$

\begin{prop}\label{iso:imminvcom} Let $f:X \to Y$ be a holomorphic map between complex
manifolds. Then
\begin{equation}\label{mor:imminvcom}
\Omega_X^t \underset{\rho_!\mathcal{D}_{\overline{X}}}{\ltens} \rho_! \ddxy
[d_X] \simeq f^!\Omega^t_Y[d_Y].
\end{equation}
\end{prop}
\dim\ \ 
We have the chain of isomorphisms
\begin{eqnarray*}
f^!(\db^{t\vee}_{Y_\R} \underset{\rho_!\mathcal{D}_{\overline{X}}}{\ltens} \rho_!\OO_{\overline{Y}}) &
\simeq & f^!\db^{t\vee}_{Y_\R} \underset{\rho_! f^{-1}
\mathcal{D}_{\overline{Y}}}{\ltens} \rho_!\imin f
\OO_{\overline{Y}} \\
 & \simeq & \db^{t\vee}_{X_\R} \underset{\rho_!\D_{X_\R}}{\ltens} \rho_! \dd{X_\R}{}{Y_\R} \underset{\rho_!\imin f\mathcal{D}_{\overline{Y}}}
 {\ltens} \rho_! \imin f \OO_{\overline{Y}} \\
 & \simeq & (\db^{t\vee}_{X_\R} \underset{\rho_!\D_{X}}{\ltens} \rho_! \ddxy) \underset{\rho_!\mathcal{D}_{\overline{X}}}
 {\ltens} \rho_! \dd{\overline{X}}{}{\overline{Y}} \underset{\rho_!\imin f \mathcal{D}_{\overline{Y}}}{\ltens}\rho_!\imin f\OO_{\overline{Y}} \\
 & \simeq & (\db^{t\vee}_{X_\R} \underset{\rho_!\D_{X}}{\ltens} \rho_! \ddxy)\underset{\rho_!\mathcal{D}_{\overline{X}}}{\ltens} \rho_! \OO_{\overline{X}} \\
 & \simeq & (\db^{t\vee}_{X_\R} \underset{\rho_!\mathcal{D}_{\overline{X}}}{\ltens} \rho_! \mathcal{O}_{\overline{X}}) \underset{\rho_!\D_X}{\ltens} \rho_! \ddxy \ ,
\end{eqnarray*}
where the second isomorphism follows from Proposition
\ref{iso:imminv}.\\ \qed

By the equivalence between left and right
$\mathcal{D}$-modules, we have an isomorphism
\begin{equation}
\label{eq:DYXotXotY}
\rho_!\ddyx\underset{\rho_!\D_X}{\overset{L}{\otimes}}\ot_X\stackrel{\sim}{\rightarrow}
f^!\ot_Y.
\end{equation}

\begin{cor}\label{smoothot} When $f$ is smooth we have an isomorphism
$$
\imin f\ot_Y \iso \rh_{\rho_!\D_X}(\rho_!\ddxy,\ot_X).
$$
\end{cor}
\dim\ \ The proof is similar to that of Corollary \ref{smoothdbt}.\\
\qed

\subsection{Inverse image for Whitney holomorphic functions}

Let $f:M \to N$ be a morphism of oriented real analytic manifolds of dimension $d_M$ and $d_N$. Set $d=d_N-d_M$.

\begin{lem}\label{degzero} The sheaf $f^!\CWN[d]$ is concentrated in degree zero.
\end{lem}
\dim\ \ If $f$ is smooth, then $f^!(\cdot)[d] \simeq \imin f$, and the result is
clear. Let $f$ be a closed embedding.  Let $F \in
D^b_{\rc}(\CC_M)$. We have the chain of isomorphisms
\begin{eqnarray*}
\Rh(D'F;f^!\CWN)[d]
& \simeq & \Rh(f_!D'F[-d];\CWN) \\
& \simeq & \Rh(D'(f_!F);\CWN ) \\
& \simeq & \mathrm{R}\Gamma(Y;f_!F\wtens \C^\infty_N).
\end{eqnarray*}
The fourth isomorphism follows since $Rf_*DF \simeq D(Rf_!F)$ if
$F \in D^b_{\rc}(\CC_M)$ and $Rf_* \simeq Rf_! \simeq f_!$ since
$f$ is a closed embedding. Let $U \in \op^c(M_{sa})$ be locally
cohomologically trivial. We have $D'\CC_{\overline{U}} \simeq
\CC_U$, and we get
$$R^{k+d}\Gamma(U;f^!\CWN) \simeq R^k\Gamma(X;f_!\CC_{\overline{U}} \wtens
\C^\infty_N)=0$$ if $k \neq 0$ since $f_!\CC_{\overline{U}} \wtens
\C^\infty_N$ is soft. Hence $f^!\CWN[d]$ is concentrated in degree zero on a basis
for the topology of $M_{sa}$ and the result follows. \\
\qed

\begin{prop} There is a natural morphism in $\mod(\CC_{M_{sa}})$
$$f^!\CWN[d] \to \CWM$$
\end{prop}
\dim\ \ Let $U \in \op^c(M_{sa})$ be locally cohomologically
trivial. We have the chain of morphisms
\begin{eqnarray*}
\Gamma(U;f^!\CWN[d])
& \simeq & \mathrm{R}\Gamma(N;Rf_!\CC_{\overline{U}} \wtens
\C_N^\infty) \\
& \to & \mathrm{R}\Gamma(M;\imin fRf_!\CC_{\overline{U}} \wtens
\C_M^\infty) \\
& \to & \mathrm{R}\Gamma(M;\CC_{\overline{U}} \wtens
\C_M^\infty) \\
& \simeq & \Gamma(U;\CWM),
\end{eqnarray*}
where the first isomorphism has been proved in Lemma \ref{degzero}
and the first arrow follows from Theorem 3.3 of \cite{KS96}.\\
\qed

\begin{prop} There is a natural morphism in $D^b(\rho_!\D_M)$:
\begin{equation}\label{admorCW}
\rho_!\ddmn \underset{\rho_!\imin f\D_M}\ltens f^!\CWN[d] \to \CWM.
\end{equation}
\end{prop}
\dim\ \ The Spencer resolution of $\ddmn$ gives rise to the
quasi-isomorphism

$$ \ddmn\stackrel{\sim}{\leftarrow} \D_M\underset{\A_M}{\otimes}\overset{\bullet}{\bigwedge}\Theta_M\underset{\A_M}
{\otimes}\ddmn \simeq \D_M\underset{\A_M
}{\otimes}\overset{\bullet}{\bigwedge}\Theta_M\underset{f^{-1}\A_M
}{\otimes}f^{-1}\D_M
$$
from which we obtain
\begin{eqnarray*}
\lefteqn{\rho_!\ddmn \underset{\rho_!\imin f\D_M}{\ltens} f^!\CWN[d]} \\ & \simeq &
\rho_!\D_M\underset{\rho_!\A_M
}{\otimes}\rho_!\overset{\bullet}{\bigwedge}\Theta_M\underset{\rho_!f^{-1}\A_M
}{\otimes}\rho_!f^{-1}\D_M \underset{\rho_!\imin f\D_M}{\ltens}
f^!\CWN[d] \\
& \simeq & \rho_!\D_M\underset{\rho_!\A_M
}{\otimes}\rho_!\overset{\bullet}{\bigwedge}\Theta_M\underset{\rho_!f^{-1}\A_M
}{\otimes} f^!\CWN[d] \\
& \to & \rho_!\D_M\underset{\rho_!\A_M
}{\otimes}\rho_!\overset{\bullet}{\bigwedge}\Theta_M\underset{\rho_!f^{-1}\A_M
}{\otimes} \CWM \\
& \simeq & \CWM.
\end{eqnarray*}
\qed

By adjunction we get a morphism
\begin{equation}\label{morCW}
 f^!\CWN[d] \to
\rh_{\rho_!\D_M}(\rho_!\ddmn,\CWM).
\end{equation}

\begin{teo}\label{isoOW} The morphism \eqref{morCW} is an isomorphism.
\end{teo}
\dim\ \ Let $F \in D^b_{\rc}(\CC_M)$. We have the chain of
isomorphisms
\begin{eqnarray*}
 \hspace{-0cm}\Rh(D'F,f^!\CWN)[d]
& \simeq & \mathrm{R}\Gamma(Y;Rf_!F \wtens \C^\infty_N) \\
& \simeq & \Rh_{\D_M}(\ddmn,F \wtens \C^\infty_M) \\
& \simeq & \Rh_{\D_M}(\ddmn,\imin \rho \rh(D'F,\CWM)) \\
& \simeq & \Rh_{\rho_!\D_M}(\rho_!\ddmn,\rh(D'F,\CWM))\\
& \simeq & \Rh(D'F,\rh_{\rho_!\D_M}(\rho_!\ddmn,\CWM)),
\end{eqnarray*}
where the second isomorphism follows from Theorem 3.5 of
\cite{KS96}.\\
\qed

\begin{cor}\label{smoothCW} When $f$ is smooth we have an isomorphism
$$
\imin f\C_N^{\infty,{\rm w}} \iso \rh_{\rho_!\D_M}(\rho_!\ddmn,\C_M^{\infty,{\rm w}}).
$$
\end{cor}
\dim\ \ It follows from the fact that $f^!(\cdot)[d] \simeq \imin f$ when $f$ is smooth.\\ \qed

From now on $X$ will be a complex manifold of complex dimension $d_X$, with structure sheaf
$\OO_X$. We denote by $\overline{X}$ the complex conjugate manifold
(with structure sheaf $\mathcal{O}_{\overline{X}}$), and $X_\R$
the underlying real analytic manifold, identified with the
diagonal of $X \times \overline{X}$. Let $\OWX$ be the sheaf of Whitney holomorphic functions on $X$.

\begin{teo}\label{isOW} Let $f:X \to Y$ be a morphism of complex manifolds. Then
\begin{equation}\label{morOW}
f^!\OWY[2d_Y]  \iso \rh_{\rho_!\D_X}(\rho_!\ddxy,\OWX)[2d_X].
\end{equation}
\end{teo}
\dim\ \ Remark that, if $\mathcal{M} \in D^b(\rho_!\D_{X_\R})$ we have
\begin{eqnarray*}
\lefteqn{\rh_{\rho_!\imin f\D_{\overline{Y}}}(\rho_! \imin f
\OO_{\overline{Y}},\rh_{\rho_!\D_{X_\R}}(\rho_!\dd
{X_\R}{}{Y_\R},\mathcal{M}))}\\
& \simeq & \rh_{\rho_!\D_X}(\rho_! \ddxy,\rh_{\rho_!\imin
f\D_{\overline{Y}}}(\rho_! \imin f
\OO_{\overline{Y}},\rh_{\D_{\overline{X}}}(\rho_!\dd
{\overline{X}}{}{\overline{Y}},\mathcal{M}))) \\
& \simeq & \rh_{\rho_!\D_X}(\rho_!
\ddxy,\rh_{\rho_!\D_{\overline{X}}}(\rho_!(\dd
{\overline{X}}{}{\overline{Y}} \underset{_{\rho_!\imin
f\D_{\overline{Y}}}}{\ltens} \imin f
\OO_{\overline{Y}}),\mathcal{M})) \\
& \simeq & \rh_{\rho_!\D_X}(\rho_!
\ddxy,\rh_{\rho_!\D_{\overline{X}}}(\rho_!\OO_{\overline{X}},\mathcal{M})).
\end{eqnarray*}

\noindent We have the chain of isomorphisms
\begin{eqnarray*}
f^! \OWY[2d_Y] & \simeq & f^! \rh_{\rho_!\D_{\overline{Y}}}(
\rho_!\OO_{\overline{Y}},\CWYR)[2d_Y] \\
& \simeq & \rh_{\rho_!\imin f\D_{\overline{Y}}}(\imin f
\rho_!\OO_{\overline{Y}},f^!\CWYR)[2d_Y] \\
& \simeq & \rh_{\rho_!\imin f\D_{\overline{Y}}}(\rho_!\imin
f\OO_{\overline{Y}},\rh_{\rho_!\D_{X_\R}}(\rho_!\dd
{X_\R}{}{Y_\R}, \CWXR))[2d_X] \\
& \simeq & \rh_{\rho_!\D_X}(\rho_!\ddxy,\rh_{\rho_!\D_{\overline{X}}}(\rho_!\OO_{\overline{X}},\CWXR))[2d_X] \\
& \simeq & \rh_{\rho_!\D_X}(\rho_!\ddxy,\OWX)[2d_X].
\end{eqnarray*}
\qed

\begin{cor}\label{smoothOW} When $f$ is smooth we have an isomorphism
$$
\imin f\OW_Y \iso \rh_{\rho_!\D_X}(\rho_!\ddxy,\OW_X).
$$
\end{cor}
\dim\ \ The proof is similar to that of Corollary \ref{smoothCW}.\\
\qed

\addcontentsline{toc}{section}{

\end{document}